%&amstex   
\input amstex\documentstyle {amsppt}  
\pagewidth{12.5 cm}\pageheight{19 cm}\magnification\magstep1
\topmatter
\title Classification of unipotent representations of simple $p$-adic groups,
II\endtitle
\author G. Lusztig\endauthor
\address Department of Mathematics, M.I.T., Cambridge, MA 02139\endaddress
\thanks{Supported in part by the National Science Foundation. This paper was 
written while the author enjoyed the hospitality of the Institut des Hautes 
\'Etudes Scientifiques.}\endthanks

\endtopmatter
\document    
\define\z{\zeta}
\redefine\T{\times}
\redefine\t{\tau}
\redefine\b{\beta}
\define\n{\notin}
\define\r{\rangle}
\define\a{\alpha}
\define\f{\forall}
\define\e{\emptyset}
\redefine\i{^{-1}}
\define\g{\gamma}
\define\h{\frac{\hphantom{aaa}}{\hphantom{aaa}}}
\define\m{\mapsto}
\define\x{\boxed}
\define\s{\star}
\define\k{\kappa}
\define\ta{\vartheta}
\define\sh{\sharp}
\define\fla{\flat}
\define\do{\dots}
\define\us{\underset}
\define\hor{\!\!\!\!\frac{\hphantom{aaa}}{\hphantom{aaa}}\!\!\!\!}
\define\Ra{\Rightarrow}
\define\Lar{\Leftarrow}
\define\dd{\bold d}
\define\hg{\hat\gamma}
\define\hG{\hat G}
\define\hZ{\hat Z}
\define\sqc{\sqcup}
\define\lan{\langle}

\define\lra{\leftrightarrow}

\define\sm{\smallmatrix}
\define\esm{\endsmallmatrix}
\define\sub{\subset}
\define\bxt{\boxtimes}
\define\ti{\tilde}
\define\nl{\newline}

\define\fra{\frac}
\define\un{\underline}

\define\ot{\otimes}

\define\ad{\text{\rm ad}}
\define\Ad{\text{\rm Ad}}
\define\Hom{\text{\rm Hom}}

\define\Irr{\text{\rm Irr}}

\define\Ker{\text{\rm Ker}}
\redefine\Im{\text{\rm Im}}

\define\card{\text{\rm card}}

\define\spa{\spadesuit}
\define\di{\diamond}
\define\opl{\oplus}
\define\sha{\sharp}

\define\de{\delta}

\define\io{\iota}

\define\rh{\rho}
\define\si{\sigma}

\define\th{\theta}

\define\la{\lambda}
\define\ph{\phi}

\define\Ga{\Gamma}

\define\Si{\Sigma}
\define\Th{\Theta}

\define\bod{\bold d}

\define\bc{\bold C}

\define\be{\bold E}
\define\bg{\bold G}

\define\bi{\bold I}
\define\bj{\bold J}
\define\bk{\bold K}

\define\bn{\bold N}

\define\bq{\bold Q}
\define\br{\bold R}

\define\bu{\bold U}

\define\bz{\bold Z}
\define\bx{\bold X}

\define\ca{\Cal A}

\define\cc{\Cal C}

\define\cf{\Cal F}
\define\cg{\Cal G}
\define\ch{\Cal H}
\define\ci{\Cal I}

\define\cl{\Cal L}

\define\cn{\Cal N}
\define\co{\Cal O}
\define\cp{\Cal P}

\define\car{\Cal R}

\define\ct{\Cal T}
\define\cu{\Cal U}

\define\cw{\Cal W}

\define\cx{\Cal X}
\define\cy{\Cal Y}

\define\fg{\frak g}
\define\fh{\frak h}

\define\fl{\frak l}

\define\fo{\frak o}
\define\fp{\frak p}

\define\fs{\frak s}
\define\ft{\frak t}

\define\fz{\frak z}

\define\fJ{\frak J}

\define\fN{\frak N}

\define\fR{\frak R}
\define\fS{\frak S}
\define\fT{\frak T}

\define\bZ{\bar Z}
\define\bac{\bar c}
\define\baz{\bar z}
\define\bam{\bar M}
\define\bah{\bar H}
\define\bas{\bar\Si}
\define\bat{\bar t}
\define\bai{\bar I}
\define\bco{\bar{\co}}
\define\bde{\bar\de}
\define\chr{\check R}
\define\cha{\check\alpha}
\define\chb{\check\beta}
\define\dcc{\dot{\cc}}
\define\dfT{\dot{\fT}}
\define\dby{\dot{\bby}}
\define\hh{\hat H}
\define\hco{\hat{\Cal O}}
\define\hcz{\hat{\Cal Z}}
\define\ug{\un G}
\define\uc{\un c}
\define\ty{\ti y}

\define\tn{\ti n}

\define\tZ{\ti Z}
\define\tig{\ti\gamma}
\define\ts{\ti{\spa}}
\define\utZ{\un{\tZ}}
\define\uct{\un{\ct}}
\define\tcf{\ti{\cf}}

\define\tsi{\ti{\si}}
\define\tbas{\ti{\bas}}

\define\ex{\exp}
\define\ul{\un L}
\define\up{\un P}
\define\uz{\un Z}
\define\uup{\un{U_P}}
\define\Mo{\text{\rm Mod}}
\define\tcr{\ti{\car}}
\define\tih{\ti h}
\define\che{\check}

\define\hah{\hat h}
\define\bby{\bold Y}
\define\KA{K}
\define\KL{KL}
\define\CU{L1}
\define\LG{L2}
\define\II{L3}
\define\IM{L4}
\define\TE{L5}
\define\RE{R}
\define\SE{S}
\define\TI{T}
\define\WA{W}

\head Introduction \endhead
\subhead 0.1\endsubhead
Let $\bk$ be a nonarchimedean local field with a residue field of cardinal $q$.
Let $\bg(\bk)$ be the group of $\bk$-rational points of a connected, adjoint 
simple algebraic group $\bg$ defined over $\bk$ which becomes split over an
unramified extension of $\bk$. Let $\cu(\bg(\bk))$ be the set of isomorphism 
classes of unipotent representations of $\bg(\bk)$ (see \cite{\IM, 1.21}). Let
$G$ be a simply connected almost simple algebraic group over $\bc$ of the type
dual to that of $\bg$ (in the sense of Langlands); let $\ta:G@>>>G$ be the
"graph automorphism" of $G$ associated to the $\bk$-rational structure of $\bg$
as in \cite{\IM, 8.1}. One of the main results of this paper is the
construction of a bijection between $\cu(\bg(\bk))$ and a set of parameters
defined in terms of $G$ and $\ta$. (See 10.11, 10.12.) This result (or rather a
close variant of it) was stated without proof in \cite{\IM, 8.1} and was proved
in \cite{\IM} assuming that $\ta=1$; it supports the Langlands philosophy. See
\cite{\IM, 0.3} for historical remarks concerning this bijection. One of the 
main observations of \cite{\IM} and the present paper is that the various 
affine Hecke algebras which arise in connection with unipotent representations
of $\bg(\bk)$ can be also found in a completely different way, in terms of $G$,
$\ta$ and certain cuspidal local systems. Then the problem reduces to 
classifying the simple modules of these "geometric affine Hecke algebras" with
parameter equal to $\sqrt{q}$. This last problem makes sense in the case where
$\sqrt{q}$ is replaced by any $v_0\in\bc^*$. This problem was solved in
\cite{\IM} assuming that $\ta=1$ and $v_0\in\br_{>0}$. In the present paper we
treat more generally the case where $\ta$ is arbitrary and $v_0$ is either $1$
or is not a root of $1$. Moreover, using results of \cite{\TE}, we determine 
which representations are tempered or square integrable. 

\subhead 0.2. Notation\endsubhead
All algebraic groups are assumed to be affine. All algebraic varieties (in 
particular, all algebraic groups) are assumed to be over $\bc$. If $G$ is an 
algebraic group, $G^0$ denotes the identity component of $G$, $\bar G$ the 
group of components of $G$, $U_G$ the unipotent radical of $G^0$, $Z_G$ the 
centre of $G$, $\ug$ the Lie algebra of $G$. For $x\in\ug$ let $Z_G(x)$ be the
centralizer of $x$ in $G$. For $x,x'\in\ug$ let $Z_G(x,x')=Z_G(x)\cap Z_G(x')$.
If $G'$ is another algebraic groups, let $\Hom(G,G')$ be the set of 
homomorphisms of algebraic groups from $G$ to $G'$. If $G'$ is a subgroup of 
$G$ and $g\in G$, we denote by $Z_{G'}(g)$ the centralizer of $g$ in $G'$; let
$N_G(G')$ be the normalizer of $G'$ in $G$. Let $\ci_G$ be the category of 
finite dimensional rational representations of $G$. If $V\in\ci_G$, then $V$ is
also a $\ug$-module. 

If $\ca$ is a subgroup of $\bc^*$, let $G^\ca$ be the set of all $g\in G$ such
that for any $V\in\ci_G$, any eigenvalue of $g:V@>>>V$ is in $\ca$. If $A$ is a
subgroup of $\bc$, let $\ug_A$ be the set of all $x\in\un G$ such that for any
$V\in\ci_G$, any eigenvalue of $x:V@>>>V$ is in $A$; let $G_A=G^{\ex(A)}$.

If $X$ is an abelian group we write $X_\bq,X_\bc$ instead of $X\ot\bq,X\ot\bc$.

Let $\k=2\pi\sqrt{-1}\in\bc$. 

Let $z=a+\sqrt{-1}b$ where $a,b\in\br$. We say that $z\ge 0$ if either $a>0$, 
or $a=0,b\ge 0$. We say that $z>0$ if either $a>0$, or $a=0,b>0$.

\head Contents\endhead
1. Preliminaries on affine Hecke algebras and graded Hecke algebras.

2. A review of \cite{\LG,\S8,\S9}.

3. Some consequences of the first reduction theorem of \cite{\LG}.

4. Some consequences of the second reduction theorem of \cite{\LG}.

5. Geometric graded Hecke algebras.

6. The subgroups $G_J$.

7. The set $\fR(G\t,G_J,\cc,\cf)$.

8. Geometric affine Hecke algebras.

9. A bijection.

10. The main results.

11. Tables.

Appendix.

\head 1. Preliminaries on affine Hecke algebras and graded Hecke algebras
\endhead
\subhead 1.1\endsubhead
(a) A {\it root system} $(R,\chr,X,Y)$ consists of two finitely generated free
abelian groups $X,Y$, a perfect pairing $\lan,\r:X\T Y@>>>\bz$, finite 
subsets $R\sub X-\{0\},\chr\sub Y-\{0\}$ and a bijection $\a\lra\cha$ between 
$R$ and $\chr$ such that for any $\a\in R$ we have $\lan\a,\cha\r=2$ and 
$s_\a:X@>>>X,x\m x-\lan x,\cha\r\a$ (resp. $s_\a:Y@>>>Y,y\m y-\lan\a,y\r\cha$)
leaves $R$ (resp. $\chr$) stable.

We sometimes write $(R,X)$ instead of $(R,\chr,X,Y)$.

(b) A {\it $\bq$-root system} $(R,\chr,E,E')$ consists of two finite 
dimensional $\bq$-vector spaces $E,E'$, a perfect bilinear pairing 
$\lan,\r:E\T E'@>>>\bq$, finite subsets $R\sub E-\{0\},\chr\sub E'-\{0\}$ 
and a bijection $\a\lra\cha$ between $R$ and $\chr$ such that 
$\lan\a,\cha'\r\in\bz$ for any $\a,\a'\in R$ and for any $\a\in R$ we have 
$\lan\a,\cha\r=2$ and $s_\a:E@>>>E$, $e\m e-\lan e,\cha\r\a$ (resp. 
$s_\a:E'@>>>E'$, $e'\m e'-\lan\a,e'\r\cha$) leaves $R$ (resp. $\chr$) stable.

We sometimes write $(R,E)$ instead of $(R,\chr,E,E')$.

We set $E_\bc=E\ot_\bq\bc,E'_\bc=E'\ot_\bq\bc$. We denote the $\bc$-bilinear
pairing $E_\bc\T E'_\bc@>>>\bc$ defined by $\lan,\r$ again by $\lan,\r$.

Unless otherwise indicated, in both cases (a),(b) it is assumed that 
$\a\in R\implies 2\a\n R$.

In the case (a) (resp. (b)) the Weyl group $W_0$ is defined as the subgroup of
$GL(X)$ or $GL(Y)$ (resp. $GL(E)$ or $GL(E')$) generated by 
$\{s_\a;\a\in R\}$. In both cases one has the standard notion of "basis" (or
"set of simple roots") of $R$ and the corresponding notion of positive roots 
$R^+$ and positive coroots $\chr^+$. A basis of $R$ always exists. If a basis 
of $R$ is given then $W_0$ is naturally a (finite) Coxeter group with length 
function $l:W_0@>>>\bn$.

\subhead 1.2\endsubhead
Assume that we are given a root system $(R,\chr,X,Y)$ and a basis $\Pi$ for it.
A {\it parameter set} consists of a function $\la:\Pi@>>>\bn$ such that 
$\la(\a)=\la(\a')$ whenever $\lan\a,\cha'\r=\lan\a',\cha\r=-1$ together with a
function $\la^*:\{\a\in\Pi;\cha\in 2Y\}@>>>\bn$. If such $(\la,\la^*)$ is 
given, we define $H^{\la,\la^*}_{R,X}$ to be the associative algebra over 
$\bc[v,v\i]$ ($v$ is an indeterminate) defined by the generators $T_w,w\in W_0$
and $\th_x,x\in X$ and by the relations

$T_wT_{w'}=T_{ww'}$ for all $w,w'\in W_0$ such that $l(ww')=l(w)+l(w')$,

$(T_{s_\a}+1)(T_{s_\a}-v^{2\la(\a)})=0$ for all $\a\in\Pi$,

$\th_{x_1}\th_{x_2}=\th_{x_1+x_2}$ for all $x_1,x_2\in X$,

$\th_x(T_{s_\a}+1)-(T_{s_\a}+1)\th_{s_\a(x)}=(\th_x-\th_{s_\a(x)})\cg(\a)$
\nl
for all $x\in X,\a\in\Pi$ where, for $\a\in\Pi$, $\cg(\a)$ equals
$$\fra{\th_\a v^{2\la(\a)-1}}{\th_\a-1} \text{ if $\cha\n 2Y$ and }
\fra{(\th_\a v^{2\la(\a)+\la^*(\a)}-1)(\th_\a v^{2\la(\a)-\la^*(\a)}+1)}
{\th_{2\a}-1} \text{ if } \cha\in 2Y.$$
(Note that $(\th_x-\th_{s_\a(x)})\cg(\a)$ is a $\bz$-linear combination of 
elements $\th_{x_1},x_1\in X$.) Now $\th_0$ is a unit element for
$H^{\la,\la^*}_{R,X}$.

Let $\ct=Y\ot\bc^*$. Let $\co$ be the algebra of regular functions 
$\ct\T\bc^*@>>>\bc$. We may identify $\co$ with the $\bc[v,v\i]$-submodule of
$H^{\la,\la^*}_{R,X}$ spanned by $\{\th_x,x\in X\}$ (a commutative subalgebra):
to $v^n\th_x$ corresponds the regular function $(t,a)\m a^nx(t)$ where 
$x(t)=\prod_n a_n^{\lan x,y_n\r}$ for
$t=\sum_ny_n\ot a_n,y_n\in Y,a_n\in\bc^*$. Now $W_0$ acts naturally on $\co$ 
and the algebra of invariants $\co^{W_0}$ is the centre of 
$H^{\la,\la^*}_{R,X}$. For any $W_0$-orbit $\Si$ on $\ct$ and $v_0\in\bc^*$ let
$J_{\Si,v_0}$ be the maximal ideal of $\co^{W_0}$ consisting of the functions 
in $\co^{W_0}$ that vanish at all points of $\Si\T\{v_0\}$. Let 
$(\co^{W_0})\hat{}$ be the $J_{\Si,v_0}$-adic completion of $\co^{W_0}$ and let
$\hh=H^{\la,\la^*}_{R,X}\ot_{\co^{W_0}}(\co^{W_0})\hat{}$.

\subhead 1.3\endsubhead
For $v_0\in\bc^*$, let $\Mo_{v_0}H^{\la,\la^*}_{R,X}$ be the category of
$H^{\la,\la^*}_{R,X}$-modules that are finite dimensional over $\bc$ and in
which $v$ acts as $v_0$ times $1$. Let $\Irr_{v_0}H^{\la,\la^*}_{R,X}$ be the 
set of isomorphism classes of simple objects of $\Mo_{v_0}H^{\la,\la^*}_{R,X}$.

Let $\Mo_{\Si,v_0}H^{\la,\la^*}_{R,X}$ be the category of
$H^{\la,\la^*}_{R,X}$-modules $M\in\Mo_{v_0}H^{\la,\la^*}_{R,X}$ that satisfy
$J_{\Si,v_0}M=0$. Let $\Irr_{\Si,v_0}H^{\la,\la^*}_{R,X}$ be the set of 
isomorphism classes of simple objects of $\Mo_{\Si,v_0}H^{\la,\la^*}_{R,X}$. We
have 
$$\Irr_{v_0}H^{\la,\la^*}_{R,X}=\sqc_\Si\Irr_{\Si,v_0}H^{\la,\la^*}_{R,X}\tag a
$$
where $\Si$ runs over the $W_0$-orbits in $\ct$.

For $M\in\Mo_{v_0}H^{\la,\la^*}_{R,X}$ and for $t\in\ct$ let $M_t$ be the
subspace of all $m\in M$ such that for any $x\in X$, $m$ is in the generalized
eigenspace of $\th_x:M@>>>M$ corresponding to the eigenvalue $x(t)\in\bc^*$. We
say that $M_t$ is a weight space of $M$. We have $M=\opl_tM_t$ where $t$ runs 
over $\ct$.

Let $\z:\bc^*@>>>\br$ be a group homomorphism such that $\z(v_0)\ne 0$. Let 
$X^+$ be the set of all $x\in X$ such that $\lan x,\cha\r\ge 0$ for all 
$\a\in\Pi$. We say that $M$ (as above) is {\it $\z$-tempered} if the following
holds: for any $t\in\ct$ such that $M_t\ne 0$ and any $x\in X^+$ we have 
$\z(x(t))/\z(v_0)\ge 0$. In the case where $R$ generates a subgroup of finite
index of $X$, we say that $M$ is {\it $\z$-square integrable} if the following
holds: for any $t\in\ct$ such that $M_t\ne 0$ and any $x\in X^+-\{0\}$ we have
$\z(x(t))/\z(v_0)>0$. 

\subhead 1.4\endsubhead
Assume that we are given a $\bq$-root system $(R,\chr,E,E')$ and a basis $\Pi$
for it. A {\it parameter set} is a function $\mu:\Pi@>>>\bz$ such that 
$\mu(\a)=\mu(\a')$ whenever $\lan\a,\cha'\r=\lan\a',\cha\r=-1$. If such $\mu$
is given, we define $\bah^\mu_{R,E}$ to be the associative algebra over 
$\bc[r]$ ($r$ is an indeterminate) defined by the generators $t_w,w\in W_0$ and
$(f), f\in\bco$ (the algebra of regular functions $E'_\bc\opl\bc@>>>\bc$) and 
by the relations 

$t_wt_{w'}=t_{ww'}$ for all $w,w'\in W_0$;

$(f_1)(f_2)=(f_1f_2)$ for all $f_1,f_2\in\bco$;

$(a_1f_1+a_2f_2)=a_1(f_1)+a_2(f_2)$ for $f_1,f_2\in\bco$ and 
$a_1,a_2\in\bc[r]$;

$(f)t_{s_\a}-t_{s_\a}(s_\a(f))=\mu(\a)r\fra{f-s_\a(f)}{\a}$
\nl
for all $f\in\bco,\a\in\Pi$, where $\a$ is regarded as a linear form on 
$E'_\bc\opl\bc$ (zero on the second factor) so that 
$\fra{f-s_\a(f)}{\a}\in\bco$. (We regard $\bco$ as a $\bc[r]$-algebra, by 
identifying $r$ with the second projection $E'_\bc\opl\bc@>>>\bc$.) Now $(0)$ 
is a unit element for $\bar H^\mu_{R,E}$.

We may identify $\bco$ with the $\bc[r]$-submodule of $\bah^\mu_{R,E}$ 
consisting of all $(f)$ with $f\in\bco$ (a commutative subalgebra): to $(f)$ 
corresponds $f$. Now $W_0$ acts naturally on $\bco$ and the algebra of 
invariants $\bco^{W_0}$ is the centre of $\bah^\mu_{R,E}$. For any $W_0$-orbit 
$\bas$ on $E'_\bc$ and $r_0\in\bc$ let $\bar J_{\bas,r_0}$ be the maximal ideal
of $\bco^{W_0}$ consisting of functions in $\bco^{W_0}$ that vanish at all 
points of $\bas\T\{r_0\}$.

\subhead 1.5\endsubhead
For $r_0\in\bc$ let $\Mo_{r_0}\bah^\mu_{R,E}$ be the category of 
$\bah^\mu_{R,E}$-modules that are finite dimensional over $\bc$ and in which 
$r$ acts as $r_0$ times $1$. Let $\Irr_{r_0}\bah^\mu_{R,E}$ be the set of 
isomorphism classes of simple objects of $\Mo_{r_0}\bah^\mu_{R,E}$.

Let $\Mo_{\bas,r_0}\bah^\mu_{R,E}$ be the category of $\bah^\mu_{R,E}$-modules
$\bam\in\Mo_{r_0}\bah^\mu_{R,E}$ that satisfy $J_{\bas,r_0}\bam=0$. Let 
$\Irr_{\bas,r_0}\bah^\mu_{R,E}$ be the set of isomorphism classes of simple 
objects of $\Mo_{\bas,r_0}\bah^\mu_{R,E}$. We have 
$$\Irr_{r_0}\bah^\mu_{R,E}=\sqc_{\bas}\Irr_{\bas,r_0}\bah^\mu_{R,E}\tag a$$
where $\bas$ runs over the $W_0$-orbits in $E'_\bc$.

For $\bam\in\Mo_{r_0}\bah^\mu_{R,E}$ and for $e'\in E'_\bc$ let $\bam_{e'}$ be
the subspace of all $m\in\bam$ such that for any $f\in\bco$, $m$ is in the 
generalized eigenspace of $(f):\bam@>>>\bam$ corresponding to the eigenvalue
$f(e',r_0)$. We say that $\bam_{e'}$ is a weight space of $\bam$. We have 
$\bam=\opl_{e'\in E'_\bc}\bam_{e'}$.

Let $\t:\bc@>>>\br$ be a group homomorphism such that $\t(r_0)\ne 0$. We say
that $\bam$ (as above) is {\it $\t$-tempered} if the following holds: for any
$e'\in E'_\bc$ such that $\bam_{e'}\ne 0$ and any $e\in E$ such that 
$\lan e,\cha\r\ge 0$ for all $\a\in\Pi$ we have 
$\t(\lan e,e'\r)/\t(r_0)\ge 0$. In the case where $R$ generates $E$ as a vector
space, we say that $\bam$ is {\it $\t$-square integrable} if the following 
holds: for any $e'\in E'_\bc$ such that $\bam_{e'}\ne 0$ and any $e\in E-\{0\}$
such that $\lan x,\cha\r\ge 0$ for all $\a\in\Pi$ we have 
$\t(\lan e,e'\r)/\t(r_0)>0$. 

\head 2. A review of \cite{\LG,\S8,\S9}\endhead
\subhead 2.1\endsubhead
Let $\spa$ be a $\bq$-subspace of $\bc$ such that $\k\bq\cap\spa=0$. Let $\ts$
be the image of $\spa$ under $\ex:\bc@>>>\bc^*$. Then $\ex$ restricts to a
group isomorphism $\spa@>\sim>>\ts$.

If $T$ is a torus, we have canonically $T=\cl\ot\bc^*$ where $\cl$ is the free
abelian group $\Hom(\bc^*,T)$. We have $T_\spa=\cl\ot\ts$. If $\ft=\un T$ we
have canonically $\ft=\cl_\bc,\ft_\spa=\cl\ot\spa$ and $\ex:\ft@>>>T$ (denoted
also by $\ex_T$) induces an isomorphism $\ft_\spa@>\sim>>T_\spa$.

\subhead 2.2\endsubhead
In this section we will refer to a subsection of \cite{\LG} such as
\cite{\LG, 8.13} simply as [8.13]. 

Now \cite{\LG,\S6,\S9} gives a method which allows one to reduce a number of 
questions on representations of an affine Hecke algebra to analogous questions
on graded Hecke algebras. Here we shall give a variation of this method. We 
will indicate how to modify \S8 and \S9 of \cite{\LG} (for example [8.13] will
become [8.13]${}'$) to obtain this variation. 

[8.1]${}'$. {\it From now on we assume that 

(a) $Y$ is generated by $\chr\cup(\fra{1}{2}\chr\cap Y)$.
\nl
Assume that a $W_0$-orbit $\Si$ in $\ct$ and an element $v_0\in\ts$ are given.
We define an equivalence relation on $\Si$ as follows: we say that $t,t'\in\Si$
are equivalent if $t't\i\in T_\spa$. Let $\cp$ be the set of equivalence 
classes. Note that $W_0$ acts transitively on $\cp$.

Let $c\in\cp$. We choose $t\in c$ and we define
$$R_c=\{\a\in R;\a(t)\in\ti\spa \text{ if } \cha\n 2Y, \a(t)\in\pm\ti\spa 
\text{ if }\cha\in 2Y\}.$$
This clearly does not depend on the choice of $t$. We set
$\chr_c=\{\cha;\a\in R_c\}$. There is a unique subset $\Pi_c$ of $R_c\cap R^+$
such that $(R_c,\chr_c,X,Y,\Pi_c)$ is a root system. Let $W_0^c$ be the Weyl 
group of this root system (a subgroup of $W_0$). Using (a) and \cite{\IM, 4.5}
we see that $W^c_0=\{w\in W_0;w(c)=c\}$. Now $c$ is a $W_0^c$-orbit in $\ct$.}

[8.2]${}'$ is empty.

[8.3]${}'$ is the same as [8.3] except that the last four lines of [8.3] are 
replaced by:

{\it Let $T_{w,c}\th_x (w\in W_0^c,x\in X)$ be the basis elements of $H_c$ 
analogous to the basis elements $T_w\th_x (w\in W_0,x\in X)$ of $H$.}

[8.4]${}'$ is the same as [8.4] except that the last three lines of [8.4] are 
deleted.

[8.5]${}'$. {\it If $\ca$ is an associative ring with $1$, denote by $\ca_n$ 
the ring of all $n\T n$ matrices with entries in $A$. We have 
$\hcz_c=\hco_c^{W_0^c}$. Thus $\hh_c$ is a $\hco_c^{W_0^c}$-algebra. The 
identity map $\co\to\co$ extends continuously to a ring homomorphism 
$i:\hco@>>>\hco_c$ (since $J_{\Si,v_0}\sub J_{c,v_0}$). This restricts to a 
ring isomorphism $\hco^{W_0}@>\sim>>\hco_c^{W_0^c}$ (since $W_0^c\sub W_0$). 
Via this isomorphism we can regard $\hh_c$ also as a $\hco^{W_0}$-algebra.}

[8.6]${}'$. {\it Theorem. If $c\in\cp$, there exists an isomorphism of 
$\hco^{W_0}$-algebras $\hh\cong(\hh_c)_n$ where $n=\card(\cp)$.}

[8.7]${}'$=[8.7].

[8.8]${}'$ is the same as [8.8] except that the reference to [8.2](b) is 
deleted.

[8.9]${}'$=[8.9].

[8.10]${}'$ remains unchanged except that formula (a) should be replaced by:

(a) $T_w^c=T_{s_{\a_1}}^cT_{s_{\a_2}}^{s_{\a_1}(c)}\do
T_{s_{\a_p}}^{s_{\a_{p-1}}\do s_{\a_2}s_{\a_1}(c)}$.

[8.11]${}'$=[8.11], [8.12]${}'$=[8.12].

[8.13]${}'$ is deleted except for the line (e) and the three lines following it
which are left unchanged.

[8.14]${}'$=[8.14], [8.15]${}'$=[8.15].

[8.16]${}'$. Lines 2,3 of 8.16 are replaced by

{\it For any $c'\in\cp$ let $w\in W_0$ be the unique element of minimal length
such that $w(c)=c'$. If $s_{\a_1}s_{\a_2}\do s_{\a_p}=w$ is a reduced
expression in $W_0$, then }

$c\ne s_{\a_p}(c)\ne s_{\a_{p-1}}s_{\a_p}(c)\ne\do\ne 
s_{\a_1}s_{\a_2}\do s_{\a_p}(c)=c'$.

The rest of [8.16] remains unchanged except that 
$[\Ga(c)],\g,\g\in\Ga(c),T_\g^c$ are deleted.

[8.17]${}'$ is empty.

[9.1]${}'$. {\it We preserve the setup of \S3. Assume that we are given a 
$W_0$-orbit $\Si$ in $\ct$ and an element $v_0\in\ts$ such that for any 
$t\in\Si$ and any $\a\in R$ we have }

$\a(t)\in\ts\text{ if } \cha\n 2Y, \a(t)\in\pm\ts\text{ if }\cha\in 2Y.$

[9.2]${}'$. The text of [9.2] except for the last three lines is replaced by 
the following:

{\it Define $r_0\in\spa$ by $\ex(r_0)=v_0$. Let $\ft=\uct$. We show that there
exists a $W_0$-invariant element $t_0\in\ct$ and a $W_0$-orbit $\bas$ in
$\ft_\spa$ such that $t_0\ex_\ct(\bas)=\Si$.

Choose a $\bq$-subspace $\di$ of $\bc$ complementary to $\spa$, and consider
the subgroup $\ex(\di)$ of $\bc^*$. Then $\bc^*=\ex(\di)\T\ts$ and we have a
$W_0$-invariant decomposition $\ct=\ct_\di\T\ct_\spa$. For any $\a\in R$ we
have $\a(\ct_\di)\in\ex(\di)$ and $\a(\ct_\spa)\in\ts$. Hence if 
$pr_1:\ct@>>>\ct_\di$ is the first projection, we have for any $t\in\Si$}
$$\a(pr_1(t))=1  \text{ if } \cha\n 2Y,\quad
\a(pr_1(t))=\pm 1 \text{ if }\cha\in 2Y.$$
{\it It follows that $pr_1(t)$ is $W_0$-invariant for any $t\in\Si$. Since 
$pr_1(\Si)$ is a single $W_0$-orbit, it follows that $pr_1(\Si)=\{t_0\}$ for 
some $W_0$-invariant $t_0\in\ct_\di$. Let $\bas$ be the unique subset of
$\ft_\spa$ such that $\ex_\ct(\bas)=t_0\i\Si$. Then $t_0,\bas$ are as required.
Another choice for $t_0,\bas$ must be of form $t_0\ex_\ct(\xi_0),\bas-\xi_0$ 
where $\xi_0\in\ft_\spa$ is $W_0$-invariant.}

The last three lines of [9.2] remain unchanged.

[9.3]${}'$=[9.3], [9.4]${}'$=[9.4].

[9.5]${}'$ is the same as [9.5] except that the last three lines of [9.5] are 
replaced by the following.

{\it Assume for example that $\cha\n 2Y$ and

(c) $\a(\bat)+2\la(\a)r_0\in\k\bz-\{0\}$.
\nl
Since $\bat\in\ft_\spa$, we have $\a(\bat)\in\spa$. Since $r_0\in\spa$ and
$\la(\a)\in\bn$, it follows that the left hand side of (c) is contained in
$\spa$. But the right hand side of (c) is not in $\spa$ and we have a 
contradiction. Similarly we see that the other statements (a),(b) hold.}

[9.6]${}'$ is the same as [9.6] except that the reference to [9.2](c) is 
deleted.

[9.7]${}'$ is empty.

\head 3. Some consequences of the first reduction theorem of \cite{\LG}\endhead
\subhead 3.1\endsubhead
We place ourselves in the setup of 2.1 and we fix $v_0\in\ts$. Let 
$$X,Y,R,\chr,\Pi,W_0,\co,\la,\la^*,\ct$$
be as in 1.2. We write $H$ instead of $H^{\la,\la^*}_{R,X}$. We assume that

(a) $Y$ {\it is generated by } $\chr\cup(\fra{1}{2}\chr\cap Y)$.
\nl
Let $\Si$ be a $W_0$-orbit on $\ct$. Let $\cp$ be as in [8.1]${}'$ (see 2.2). 
Let $c\in\cp$. Recall that $c$ is a $W_0^c$-orbit on $\ct$. Let 
$R_c,\chr'_c,\Pi_c,W_0^c$ be as in [8.1]${}'$ (see 2.2). Let $H_c$ be the 
algebra defined in the same way as $H$, but in terms of 
$(X,Y,R_c,\chr_c,\Pi_c)$ instead of $(X,Y,R,\chr,\Pi)$; the parameter set 
$(\la_c,\la^*_c)$ that we use to define $H_c$ is given by $\la_c(\a)=\la(\a')$,
$\la^*_c(\a)=\la^*(\a')$, where $\a\in\Pi_c,\a'\in\Pi$ are in the same 
$W_0$-orbit. (This does not depend on the choice of $\a'$.) Note that $\co$ is
a subalgebra of $H_c$ in the same way as $\co$ is a subalgebra of $H$. 

Let $J_{c,v_0}$ be the maximal ideal of $\co^{W_0^c}$ (the centre of $H_c$) 
consisting of the functions in $\co^{W_0^c}$ that vanish at all points of 
$c\T\{v_0\}$. Let $(\co^{W_0^c})\hat{}$ be the $J_{c,v_0}$-adic completion of
$\co^{W_0^c}$ and let $\hh_c=H_c\ot_{\co^{W_0^c}}(\co^{W_0^c})\hat{}$.

Assume that $M\in\Mo_{c,v_0}H_c$. In particular in the $H_c$-module $M$ we have
$J_{c,v_0}M=0$. Hence $M$ extends naturally to an $\hh_c$-module. The 
$\bc[v,v\i]$-module $M^\cp=M\opl M\opl\do\opl M$ (one summand for each 
$c'\in\cp$) is naturally a module over the algebra of matrices with entries in 
$\hh_c$ indexed by $\cp\T\cp$. The first reduction theorem \cite{\LG, 8.6}, in
the variant [8.6]${}'$ (see 2.2), gives an explicit isomorphism $\io$ of this
algebra of matrices with the algebra $\hh$ (see 1.3). Via this isomorphism, 
$M^\cp$ becomes an $\hh$-module and, by restriction, an $H$-module in 
$\Mo_{\Si,v_0}H$. From the definition we see that, if $f\in\co$ (regarded as an
element of $\hh$), then the $(c',c'')$-entry of $\io\i(f)$ (for $c',c''$ in 
$\cp$) is $0$ if $c'\ne c''$ and is $w'{}\i(f)$ if $c'=c''$; here $w'\in W_0$ 
is the unique element of minimal length of $W_0$ such that $w'(c)=c'$.

\proclaim{Lemma 3.2} The rule $M\m M^\cp$ is a bijection
$\Irr_{c,v_0}H_c@>\sim>>\Irr_{\Si,v_0}H$.
\endproclaim
This is an immediate consequence of the definitions and of \cite{\LG, 8.6}, in
the variant [8.6]${}'$ (see 2.2).

\medpagebreak

Let $S$ be the set of all $w\in W_0$ such that the length of $w$ is minimal in
$wW^c_0$, or equivalenty, such that $w(\cha)\in\chr^+$ for any $\a\in\Pi_c$. 

\proclaim{Lemma 3.3} Let $x\in X$ be such that $\lan x,\cha\r\ge 0$ for all 
$\a\in\Pi_c$. Let $w\in W_0$ and $x'\in X$ be such that $x=w\i(x')$,
$\lan x',\cha'\r\ge 0$ for all $\a'\in\Pi$ and $w\i(\cha')\in\chr^+$ for any 
$\a'\in\Pi$ for which $\lan x',\cha'\r=0$. Then $w\in S$.
\endproclaim
If $w\n S$ then there exists $\a\in\Pi_c$ such that $w(\cha)\in -\chr^+$
that is, $w(\cha)=\sum_{\a'\in\Pi}n_{\a'}\cha'$ where $-n_{\a'}\in\bn$. It 
follows that 
$$0\le\lan x,\cha\r=\lan x',w(\cha)\r=\sum_{\a'}n_{\a'}\lan x',\cha'\r.$$
Since $n_{\a'}\lan x',\cha'\r\le 0$ for all $\a'$, it follows that
$n_{\a'}\lan x',\cha'\r=0$ for all $\a'$. Hence for any $\a'\in\Pi$ such that
$\lan x',\cha'\r\ne 0$ we have $n_{\a'}=0$. In other words,
$$w(\cha)=\sum_{\a'\in\Pi;\lan x',\cha'\r=0}n_{\a'}\cha'.$$
Hence $\cha=\sum_{\a'\in\Pi;\lan x',\cha'\r=0}n_{\a'}w\i(\cha')$. For each 
$\a'$ in the sum, we have $w\i(\cha')\in\chr^+$ and $n_{\a'}\le 0$ hence 
$\cha\in -\chr^+$, a contradiction. The lemma is proved.

\proclaim{Lemma 3.4} Let $M\in\Mo_{c,v_0}H_c$. Assume that $\z:\bc^*@>>>\br$ is
a homomorphism such that $\z(v_0)\ne 0$. The following two conditions are
equivalent:

(i) the $H_c$-module $M$ is $\z$-tempered;

(ii) the $H$-module $M^\cp$ is $\z$-tempered.
\endproclaim
Let $D$ (resp. $D'$) be the set of all $t\in\ct$ such that $M_t\ne 0$ (resp. 
$M^\cp_t\ne 0$). By the description of $\io\i(f)$ given in 3.1, we see that 
$D'=\cup_{w\in S}w(D)$ where $S$ consists of all elements $w\in W_0$ such that
the length of $w$ is minimal in $wW^c_0$, or equivalenty, such that 
$w(\cha)\in\chr^+$ for any $\a\in\Pi_c$. Hence (ii) is equivalent to the 
following condition:

for any $t\in D$, any $w\in S$ and any $x\in X^+$, we have 
$\z(x(w(t)))/\z(v_0)\ge 0$, or equivalently $\z((w\i x)(t))/\z(v_0)\ge 0$.

We see that it is enough to show that the following two conditions for $x\in X$
are equivalent:

(iii) $\lan x,\cha\r\ge 0$ for all $\a\in\Pi_c$;

(iv) there exists $w\in S$ and $x'\in X^+$ such that $x=w\i(x')$.

Assume first that (iv) holds. Write $x=w\i(x')$ as in (iv). Let $\a\in\Pi_c$. 
Since $w\in S$, we have $w(\cha)\in\chr^+$. Using (iv) we deduce that 
$\lan x',w(\cha)\r\ge 0$. Thus $\lan w\i(x'),\cha\r\ge 0$ so that (iii) holds.

Assume next that (iii) holds. We can write uniquely $x=w\i(x')$ where $x'\in X$
satisfies $\lan x',\chb\r\ge 0$ for all $\b\in\Pi$ and $w\in W_0$ is such that
$w\i(\chb)\in\chr^+$ for any $\b\in\Pi$ for which $\lan x',\chb\r=0$. By 3.3
we have $w\in S$. Hence (iv) holds. The lemma is proved.

\proclaim{Lemma 3.5} Assume that $R$ generates a subgroup of finite index in 
$X$. Let $M\in\Mo_{c,v_0}H_c$. Assume that $\z:\bc^*@>>>\br$ is a homomorphism
such that $\z(v_0)\ne 0$. The following two conditions are equivalent:

(i) $R_c$ generates a subgroup of finite index of $X$ and the $H_c$-module $M$
is $\z$-square integrable;

(ii) the $H$-module $M^\cp$ is $\z$-square integrable.
\endproclaim
Assume that (ii) holds but $R_c$ generates a subgroup of infinite index of $X$.
Then $\chr_c$ generates a subgroup of infinite index of $Y$ hence we can find 
$z\in X-\{0\}$ such that $\lan z,\cha\r=0$ for any $\a\in\Pi_c$. We can write
uniquely $z=w'{}\i(x')$ where $x'\in X-\{0\}$ satisfies $\lan x',\chb\r\ge 0$
for all $\b\in\Pi$ and $w'\in W_0$ is such that $w'{}\i(\chb)\in\chr^+$ for any
$\b\in\Pi$ for which $\lan x',\chb\r=0$. By 3.3, we have $w'\in S$. Thus,
$z=w'{}\i(x')$ where $x'\in X-\{0\}$ satisfies $\lan x',\chb\r\ge 0$ for all 
$\b\in\Pi$ and $w'\in S$. The same argument can be applied to $-z$ instead of
$z$. We see that $-z=w''{}\i(x'')$ where $x''\in X-\{0\}$ satisfies
$\lan x'',\chb\r\ge 0$ for all $\b\in\Pi$ and $w''\in S$.

Since (ii) holds, and $D'=\cup_{w\in S}w(D)$ ($D,D'$ as in the proof of 3.4),
we see that, for any $t\in D$, any $w\in S$ and any $x\in X^+-\{0\}$, we have 
$\z(x(w(t)))/\z(v_0)>0$ that is, $\z((w\i x)(t))/\z(v_0)>0$. In particular, for
any $t\in D$ we have 
$$\z((w'{}\i x')(t))/\z(v_0)>0 \text{ and }
\z((w''{}\i x'')(t))/\z(v_0)>0.$$
We have $0=z-z=w'{}\i(x')+w''{}\i(x'')$ hence
$$\align&0<\z((w'{}\i x')(t))/\z(v_0)+\z((w''{}\i x'')(t))/\z(v_0)\\&
=\z((w'{}\i x')(t)(w''{}\i x'')(t))/\z(v_0)\\&
=\z((w'{}\i x'+w''{}\i x'')(t))/\z(v_0)=\z(0(t))/\z(v_0)=\z(1)/\z(v_0)=0.
\endalign$$
This is a contradiction. We see that (ii) implies the first condition in (i). 
Now the proof continues exactly as in 3.4; in particular we see that the 
equivalence of (i) and (ii) follows from the equivalence of 3.4(iii) and 
3.4(iv). The lemma is proved.

\head 4. Some consequences of the second reduction theorem of \cite{\LG}
\endhead
\subhead 4.1\endsubhead
We place ourselves in the setup of 2.1 and we fix $v_0\in\ts$. Define 
$r_0\in\spa$ by $\ex(r_0)=v_0$. Let 
$$X,Y,R,\chr,\Pi,W_0,\la,\la^*,\ct$$
be as in 1.2. We write $H$ instead of $H^{\la,\la^*}_{R,X}$. Let $\ft=\uct$. 
Assume that we are given a $W_0$-orbit $\Si$ in $\ct$ such that for any 
$t\in\Si$ and any $\a\in R$ we have
$$\a(t)\in\ts \text{ if } \cha\n 2Y,\quad
\a(t)\in\pm\ts \text{ if }\cha\in 2Y.$$
As in [9.2]${}'$ (see 2.2) we can find $t_0\in\ct$ ($W_0$-invariant) and a
$W_0$-orbit $\bas$ in $\ft_\spa$ such that $\Si=t_0\ex_\ct(\bas)$. 

Now $(R,\chr,X_\bq,Y_\bq)$ is a $\bq$-root system with basis $\Pi$ and with a 
parameter set $\mu:\Pi@>>>\bz$ defined by
$$\mu(\a)=2\la(\a) \text{ if $\cha\n 2Y$ and }
\mu(\a)=\la(\a)+\a(t_0)\la^*(\a) \text{ if } \cha\in 2Y.$$
(In the last equality we have $\a(t_0)=\pm 1$ since $s_\a(t_0)=t_0$.) Let 
$\bah=\bah^\mu_{R,E}$. Let $\bco\sub\bah$ be as in 1.4. We have 
$Y_\bc=\ft$. Define $\Psi:\ft\opl\bc@>>>\ct\T\bc^*$ by 
$(e',z)\m(t_0\ex_\ct(e'),\ex(z))$.

Let $(\bco^W)\hat{}$ be the $J_{\bas,r_0}$-adic completion of $\bco^W$. Let 
$\bam\in\Mo_{\bas,r_0}\bah$. Since $J_{\bas,r_0}\bam=0$, $\bam$ extends
naturally to a module over $\hat{\bah}=\bah_{\bco^W}(\bco^W)\hat{}$. The second
reduction theorem \cite{\LG, 9.3}, in the variant [9.3]${}'$ (see 2.2), gives 
an explicit algebra isomorphism of $\hh$ (as in 1.2) with $\hat{\bah}$. Via 
this isomorphism, $\bam$ becomes an $\hh$-module and, by restriction, an 
$H$-module $\bam^\dag\in\Mo_{\Si,v_0}H$. Note that $\bam^\dag$ and $\bam$ have
the same underlying $\bc$-vector space. Let $e'\in Y_\bc=\ft$. From the 
definitions we see that the $e'$-weight space $\bam_{e'}$ of $\bam$ is equal to
the $t$-weight space $\bam^\dag_t$ of $\bam^\dag$ where $t\in\ct$ is defined by
$$x(t)=x(t_0)\ex\lan x,e'\r\tag a$$
for all $x\in X$.

\proclaim{Lemma 4.2} The rule $\bam\m\bam^\dag$ is a bijection
$\Irr_{\bas,r_0}\bah@>\sim>>\Irr_{\Si,v_0}H$.
\endproclaim
This is an immediate consequence of the definitions and of \cite{\LG, 9.3}, in
the variant [9.3]${}'$ (see 2.2).

\proclaim{Lemma 4.3}Let $\bam\in\Mo_{\bas,r_0}\bah$. Assume that 
$\z:\bc^*@>>>\br$ is a homomorphism such that $\z(v_0)\ne 0$. Assume that 
$t_0\in\ct^{\Ker\z}$. Define a homomorphism $\t:\bc@>>>\br$ by 
$\t(z)=\z(\ex(z))$. Then $\t(r_0)\ne 0$. The following two conditions are 
equivalent:

(i) the $\bah$-module $\bam$ is $\t$-tempered;

(ii) the $H$-module $\bam^\dag$ is $\z$-tempered.
\endproclaim
In view of 4.1(a) it is enough to show that for $e'\in Y_\bc$, the following
two conditions are equivalent:

(iii) for any $e\in X_\bq$ such that $\lan e,\cha\r\ge 0$ for all $\a\in\Pi$ we
have \linebreak $\z(\ex\lan e,e'\r)/\z(v_0)\ge 0$;

(iv) for any $x\in X^+$ we have $\z(x(t_0)\ex\lan x,e'\r)/\z(v_0)\ge 0$.

Since $t_0\in\ct^{\Ker\z}$, for any $x\in X$ we have $\z(x(t_0))=0$ and in
(iv) we have 
$$\z(x(t_0)\ex\lan x,e'\r)=\z(\ex\lan x,e'\r).\tag a$$
Since $X\sub X_\bq$ it follows that, if (iii) holds then (iv) holds.

Assume now that (iv) holds. Let $e\in X_\bq$ be such that $\lan e,\cha\r\ge 0$
for all $\a\in\Pi$. We can find $n\in\bn-\{0\}$ such that $ne\in X^+$. Since 
(iv) holds, it follows that $\z(\ex\lan ne,e'\r)/\z(v_0)\ge 0$. (We use (a).)
Hence $n\z(\ex\lan e,e'\r)/\z(v_0)\ge 0$ so that
$\z(\ex\lan e,e'\r)/\z(v_0)\ge 0$. Thus (iii) holds. The lemma is proved.

\proclaim{Lemma 4.4} Let $\bam\in\Mo_{\bas,r_0}\bah$. Assume that 
$\z:\bc^*@>>>\br$ is a homomorphism such that $\z(v_0)\ne 0$. Define a 
homomorphism $\t:\bc@>>>\br$ by $\t(z)=\z(\ex(z))$. Then $\t(r_0)\ne 0$. Assume
that $R$ generates $X_\bq$ as a $\bq$-vector space. The following two 
conditions are equivalent:

(i) the $\bah$-module $\bam$ is $\t$-square integrable;

(ii) the $H$-module $\bam^\dag$ is $\z$-square integrable.
\endproclaim
In view of 4.1(a) it is enough to show that for $e'\in Y_\bc$, the following
two conditions are equivalent:

(iii) for any $e\in X_\bq-\{0\}$ such that $\lan e,\cha\r\ge 0$ for all 
$\a\in\Pi$ we have \linebreak $\z(\ex\lan e,e'\r)/\z(v_0)>0$;

(iv) for any $x\in X^+-\{0\}$ we have $\z(x(t_0)\ex\lan x,e'\r)/\z(v_0)>0$.

This is shown in the same way as in the proof of 4.3. (In this case we have 
automatically $\z(x(t_0))=0$. Indeed, since $t_0$ is $W_0$-invariant and $R$ 
generates $X_\bq$, we see that $t_0$ has finite order in $\ct$ hence 
$\z(x(t_0))$ has finite order in $\br$ hence $\z(x(t_0))=0$.) The lemma is 
proved.

\head 5. Geometric graded Hecke algebras \endhead
\subhead 5.1\endsubhead
If $G$ is an algebraic group, the exponential map $\ex:\ug@>>>G$ restricts to a
bijection $\ug_\spa@>\sim>>G_\spa$.

\subhead 5.2\endsubhead
In this section we review some results of \cite{\CU}, \cite{\TE} and give some
variants of them.

Assume that $G$ is a connected reductive algebraic group. Let $\fg=\ug$. Let
$L$ be the Levi subgroup of some parabolic subgroup of $G$. Let $\cc$ be a 
nilpotent $L$-orbit in $\ul$ and let $\cf$ be an irreducible $L$-equivariant 
cuspidal local system (over $\bc$) on $\cc$. Let $T=Z^0_L$. Let $\ft=\un T$. We
have $\fg=\opl_{\a\in\ft^*}\fg^\a$ where 
$$\fg^\a=\{x\in\fg;[y,x]=\a(y)x\quad\f y\in\ft\}.$$
Let $R'=\{\a\in\ft^*-\{0\};\fg^\a\ne 0\},R=\{\a\in R';\a/2\n R'\}$. The group
$W=N(T)/L$, where $N(T)$ is the normalizer of $T$ in $G$, acts naturally on 
$\ft$ and $\ft^*$. For any $\a\in R$ there is a unique element $s_\a\in W$
which acts on $\ft^*$ as a reflection sending $\a$ to $-\a$; there is a unique
element $\cha\in\ft$ such that $s_\a(x)=x-x(\cha)\a$ for all $x\in\ft^*$. Let 
$\chr=\{\cha;\a\in R\}$. We have canonically $T=Y\ot\bc^*$, $\ft=Y_\bc$ where
$Y$ is the group of all one-parameter subgroups of $T$. Let
$$\ft_\bq=Y_\bq,\ft^*_\bq=\{x\in\ft^*;x(z)\in\bq\quad\f z\in\ft_\bq\}.$$
Then $(R,\chr,\ft^*_\bq,\ft_\bq)$ is a $\bq$-root system. Let $\Pi$ a basis for
it. Let $y_0\in\cc$. For any $\a\in\Pi$ we denote by $\uc(\a)$ the integer 
$\ge 2$ such that 
$\ad(y)^{\uc(\a)-2}:\fg^\a\opl\fg^{2\a}@>>>\fg^\a\opl\fg^{2\a}$ is $\ne 0$ and 
$\ad(y)^{\uc(\a)-1}:\fg^\a\opl\fg^{2\a}@>>>\fg^\a\opl\fg^{2\a}$ is $0$. (This 
is independent of the choice $y_0$.) Then $\a\m\uc(\a)$ is a parameter set (see
1.4) for our $\bq$-root system. The corresponding algebra 
$\bah_{R,\ft^*_\bq}^{\uc}$ (see 1.4) is denoted by $\bah(G,L,\cc,\cf)$.

\subhead 5.3\endsubhead
Let $r_0\in\bc$. Let $\si$ be a semisimple element of $\fg$ and let $y$ be a 
nilpotent element of $\fg$ such that $[\si,y]=2r_0y$. Let $P$ be a parabolic 
subgroup $P$ of $G$ with Levi subgroup $L$. Let
$$\bx_{\si,y}=\{g\in G;\Ad(g\i)y\in\cc+\uup,\Ad(g\i)\si\in\up\}.\tag a$$
We have an obvious map $\bx_{\si,y}@>>>\cc$ which takes $g$ to the image of
$\Ad(g\i)y$ under $\cc+\uup@>>>\cc,a+b\m a$. The inverse image of $\cf$ under
this map is denoted again by $\cf$. On $\bx_{\si,y}$ we have a free $P$-action
by right translation and $\cf$ is $P$-equivariant hence it descends to a local
system $\tcf$ on $\bx_{\si,y}/P$. The group $Z_G(\si,y)$ acts on 
$\bx_{\si,y}/P$ by left translation and $\tcf$ is naturally a 
$Z_G(\si,y)$-equivariant local system. Then $\bZ_G(\si,y)$ acts naturally on 
the cohomology 
$$\opl_nH^n_c(\bx_{\si,y}/P,\tcf).\tag b$$
The set of irreducible representations (up to isomorphism) of $\bZ_G(\si,y)$
which appear in the representation (b) is denoted by $\Irr_0\bZ_G(\si,y)$.
(This set is independent of the choice of $P$ since another choice of $P$ is
of the form $nPn\i$ where $nLn\i=L,\Ad(n)\cc=\cc$.)

Let $\fS(G,L,\cc,\cf,r_0)$ be the set consisting of all triples $(\si,y,\rh)$ 
(modulo the natural action of $G$) where $\si,y$ are as above and 
$\rh\in\Irr_0\bZ_G(\si,y)$.

In \cite{\TE} a canonical bijection
$$\Irr_{r_0}\bah(G,L,\cc,\cf)\lra\fS(G,L,\cc,\cf,r_0)\tag c$$
is established using geometric methods (equivariant homology).

\subhead 5.4\endsubhead
We fix elements $e^0,h^0,f^0$ in $\ul$ which satisfy the standard relations of
$\fs\fl_2$ and $e^0\in\cc$. Let 
$$\align&Z=\{g\in G;\Ad(g)e^0=e^0,\Ad(g)h^0=h^0,\Ad(g)f^0=f^0\},\\&\tZ=\{(g,a)
\in G\T\bc^*;\Ad(g)e^0=a^2e^0,\Ad(g)h^0=h^0,\Ad(g)f^0=a^{-2}f^0\}.\endalign$$
We have 
$$\align&\uz=\{x\in G;[x,e^0]=0,[x,h^0]=0,[x,f^0]=0\},\\&
\utZ=\{(x,a)\in\fg\T\bc;[x,e^0]=2ae^0,[x,h^0]=0,[x,f^0]=-2af^0\}.\endalign$$
There is a unique isomorphism of algebraic groups $\io:Z^0\T\bc^*@>\sim>>\tZ^0$
such that the induced Lie algebra isomorphism $\uz\opl\bc@>\sim>>\utZ$ is given
by $(x,a)\m x+ah^0$.

\proclaim{Lemma 5.5} The inclusion $\tZ^0@>>>G\T\bc^*$ induces an injective map
from the set of semisimple $\tZ^0$-orbits in $\utZ$ to the set of semisimple 
$G\T\bc^*$-orbits in $\fg\opl\bc$.
\endproclaim
This has been stated without proof and used in \cite{\CU, 14.3(a)},
\cite{\II, 8.13}. The proof is given in the appendix.

\subhead 5.6\endsubhead
Let $\bas$ be a $W$-orbit on $\ft$ and let
$\bam\in\Irr_{\bas,r_0}\bah(G,L,\cc,\cf)$. Assume that $\bam$ corresponds to
$\si,y,\rh$ under 5.3(c).

The centre $\bco^W$ of $\bam\in\Irr_{\bas,r_0}\bah(G,L,\cc,\cf)$ acts on $\bam$
via a character which may be identified with the $W$-orbit $\bas\T\{r_0\}$ on
$\ft\T\bc$. By \cite{\II, 8.13}, the centre $\bco^W$ (as in 1.4) of $\bah$ is
identified with the equivariant cohomology $H^*_{Z^0\T\bc^*}(point)$ which via
$\io$ is identified with $H^*_{\tZ^0}(point)$; from the definitions, the 
natural action of $H^*_{\tZ^0}(point)$ on $\bam$ is via a character that may be
identified with a semisimple orbit in $\utZ$ which is contained in the 
$\Ad(G\T\bc^*)$-orbit of $(\si,r_0)$ in $\fg\opl\bc$.

Thus, $\si$ is related to $\bas$ as follows: there exists $(\si',r_0)\in\utZ$ 
such that $\si,\si'$ are conjugate under $G$ and 
$\io\i(\si',r_0)=(\si'-r_0h_0,r_0)\in\bas\T\{r_0\}$.

We define a map 
$$\fS(G,L,\cc,\cf,r_0)@>>>\ft/W.\tag a$$
Consider an element $u$ of $\fS(G,L,\cc,\cf,r_0)$ represented by $(\si,y,\rh)$.
We can find

(b) $\si'$ in the $G$-orbit of $\si$ such that $\si'-r_0h^0\in\ft$. 
\nl
Indeed, since $\bx_{\si,y}\ne\e$, there exists $g\in G$ such that 
$\si_1=\Ad(g\i)\si\in\ul,y'=\Ad(g\i)y\in e^0+\uup$. From $[\si_1,y']=2r_0y'$ we
deduce that $[\si_1,e^0]=2r_0e^0$. Now we can find $l\in L$ such that 
$\Ad(l)e^0\in\bc^*e^0$ and such that $\si'=\Ad(l)\si_1$ satisfies $\si'\in\ul$,
$[\si',e^0]=2r_0e^0,[\si',h^0]=0,[\si',f^0]=-2r_0f^0$. Since $e^0$ is 
distinguished in $\ul$, it follows that $\si'-r_0h^0\in\ft$, as required.

Let $\si'$ be as in (b) and let $\bas$ be the $W$-orbit of $\si'-r_0h^0$ in
$\ft$. Let $\tsi'$ be another element like $\si'$ and let $\tbas$ be the
$W$-orbit of $\tsi'-r_0h^0$ in $\ft$. Then $(\si',r_0),(\tsi',r_0)$ are
semisimple elements of $\utZ$ and are in the same $G\T\bc^*$-orbit in
$\fg\opl\bc$ hence, by 5.5, there exists $(g',a)\in\tZ^0$ such that
$\Ad(g')\si'=\tsi'$. Since $h^0$ is central in $\utZ$, we have
$\Ad(g')\si'=\tsi'$. Since $(h^0,1)$ is central in $\utZ$, we have
$\Ad(g')h^0=h^0$. Hence $\Ad(g')(\si'-r_0h^0)=\tsi'-r_0h^0$. Since 
$\si'-r_0h^0,\tsi'-r_0h^0$ belong to $\ft$ hence to the Cartan subalgebra 
$\ft\opl\bc h^0$, they are in the same orbit of the Weyl group of $\tZ^0$ with
respect to that Cartan subalgebra, which may be identified with $W$. It follows
that $\bas=\tbas$. Thus we have a map $u\m\bas$ as in (a). 

For any $W$-orbit $\bas$ in $\ft$ let $\fS(G,L,\cc,\cf,\bas,r_0)$ be the 
inverse image of $\bas$ under the map (a). The previous discussion yields the 
following result.

\proclaim{Lemma 5.7} Let $\bas$ be a $W$-orbit in $\ft$. The bijection 5.3(c) 
restricts to a bijection
$\Irr_{\bas,r_0}\bah(G,L,\cc,\cf)\lra\fS(G,L,\cc,\cf,\bas,r_0)$.
\endproclaim

\proclaim{Lemma 5.8} In the setup of 2.1, assume that $r_0\in\spa$. The 
bijection 5.3(c) restricts to a bijection
$$\Irr_{r_0}^\spa\bah(G,L,\cc,\cf)\lra\fS^\spa(G,L,\cc,\cf,r_0)$$
where $\Irr_{r_0}^\spa\bah(G,L,\cc,\cf)=\sqc\Irr_{\bas;r_0}\bah(G,L,\cc,\cf)$
(union over all $W$-orbits $\bas$ in $\ft_\spa$) and 
$\fS^\spa(G,L,\cc,\cf,r_0)$
consists of all $(\si,y,\rh)$ in $\fS(G,L,\cc,\cf,r_0)$ such that 
$\si\in\fg_\spa$.
\endproclaim
Using 5.7 we see that it is enough to show that, if
$(\si,y,\rh)\in\fS(G,L,\cc,\cf,r_0)$ and $\bas\in\ft/W$ correspond to each
other under 5.6(a), then we have $\si\in\fg_\spa$ if and only if 
$\bas\sub\ft_\spa$. We may assume that $\si-r_0h^0\in\bas$. Now $r_0h^0$ has 
eigenvalues in $\bz r_0$ in any $V\in\ci_G$ (a property of $\fs\fl_2$) hence 
$r_0h^0\in\fg_\spa$ (since $r_0\in\spa$). If we have two elements of $\fg_\spa$
that commute, then their sum is again in $\fg_\spa$. Applying this to the 
commuting elements $\si,-r_0h^0$ and to the commuting elements 
$\si-r_0h^0,r_0h^0$ we deduce that $\si\in\fg_\spa$ if and only if 
$\si-r_0h^0\in\fg_\spa$ that is, if and only if $\bas\sub\fg_\spa$. It remains
to observe that $\fg_\spa\cap\ft=\ft_\spa$. The lemma is proved.

The following result is closely related to \cite{\TE, 1.21}.
\proclaim{Lemma 5.9} Assume that $\t:\bc@>>>\br$ is a group homomorphism such 
that $\t(r_0)\ne 0$. Let $\bam\in\Irr_{r_0}\bah(G,L,\cc,\cf)$ and let 
$(\si,y,\rh)$ correspond to $\bam$ under 5.3(c). The following two conditions 
are equivalent.

(i) $\bam$ is $\t$-tempered;

(ii) there exists $h\in\fg,\ty\in\fg$ such that
$[y,\ty]=h,[h,y]=2y,[h,\ty]=-2\ty$, $[\si,h]=0,[\si,\ty]=-2r_0\ty$
and such that $\si-r_0h\in\fg^{\Ker\t}$.
\endproclaim
Assume that the lemma holds for $(G_i,L_i,\cc_i,\cf_i)$, $i=1,2$, (two data 
like $(G,L,\cc,\cf)$). Then one checks easily that it also holds also for
$(G_1\T G_2,L_1\T L_2,\cc_1\T\cc_2,\cf_1\bxt\cf_2)$. Without loss of 
generality we can assume that $G=C\T G'$ where $G'$ is semisimple and $C$ is a
torus. Hence it is enough to prove the lemma assuming that $G$ is either
semisimple or a torus.

Assume first that $G$ is a torus. Then $G=L=T$. In this case, we must verify
that for any $e'\in\ft$, the following two conditions are equivalent.

(iii) For any $e\in\ft^*_\bq$ we have $\t(\lan e,e'\r)/\t(r_0)\ge 0$; 

(iv) $e'\in\ft^{\Ker\t}=Y\ot\Ker\t$. 

Now (iii) is equivalent to the condition that, for any $e\in\ft^*_\bq$, we have
$\t(\lan e,e'\r)=0$. One checks easily that this is equivalent to (iv).

Next, assume that $G$ is semisimple. The condition that $\bam$ is 
$\t$-tempered (as given in 1.5) is in this case equivalent to the condition 
that $\bam$ is $\t$-tempered according to the definition in \cite{\TE, 1.20}.
(This follows easily from \cite{\TE, 3.6}.) Using \cite{\TE, 1.21} we are
therefore reduced to verifying the following statement.

If $x\in\fg$ and any eigenvalue $\nu$ of $\ad(x):\fg@>>>\fg$ satisfies 
$\t(\nu)=0$, then for any $V\in\ci_G$, 

(a) {\it any eigenvalue $\nu$ of $x:V@>>>V$ satisfies $\t(\nu)=0$ (that is, }
$x\in\fg^{\Ker\t}$). 

Clearly, if (a) holds for $V$ then it also holds for $V^{\ot n}$ for any 
$n\ge 0$. Since (a) holds for the adjoint representation it also holds for
tensor powers of the adjoint representations, hence for any direct summand of 
such a tensor power, hence for any irreducible $V$ on which $Z_G$ 
acts trivially, hence for any $V$ on which $Z_G$ acts trivially. If
$V\in\ci_G, n\ge 1$ and (a) holds for $V^{\ot n}$ then it also holds for $V$.
(Indeed, if $\nu$ is an eigenvalue of $x:V@>>>V$ then $n\nu$ is an eigenvalue 
of $x:V^{\ot n}@>>>V^{\ot n}$ hence $\t(n\nu)=0$ hence $n\t(\nu)=0$ hence 
$\t(\nu)=0$.) If $V\in\ci_G$ then there exists $n\ge 1$ such that $Z_G$ acts 
trivially on $V^{\ot n}$. As we have seen earlier, (a) holds for
$V^{\ot n}$ hence it holds for $V$. Thus (a) holds in general. The lemma is
proved.

The following result is closely related to \cite{\TE, 1.22}.
\proclaim{Lemma 5.10} Assume that $G$ is semisimple. Assume that 
$\t:\bc@>>>\br$ is a group homomorphism such that $\t(r_0)\ne 0$. Let 
$\bam\in\Irr_{r_0}\bah(G,L,\cc,\cf)$ and let $(\si,y,\rh)$ correspond to $\bam$
under 5.3(c). The following two conditions are equivalent.

(i) $\bam$ is $\t$-square integrable;

(ii) there exists $h\in\fg,\ty\in\fg$ such that
$[y,\ty]=h,[h,y]=2y,[h,\ty]=-2\ty$, $\si=r_0h$; moreover, $y$ is distinguished.
\endproclaim
Using \cite{\TE, 3.6}, we see that the condition that $\bam$ is $\t$-square 
integrable (as given in 1.5) is in our case equivalent to the condition that 
$\bam$ is $\t$-square integrable according to the definition in 
\cite{\TE, 1.20}. Hence the lemma follows from \cite{\TE, 1.22}.

\subhead 5.11\endsubhead
For any $\xi\in\fg$ let
$$\bby_\xi=\{g\in G;\Ad(g\i)\xi\in\cc+\ft+\uup\}.\tag a$$
We have an obvious map $\bby_\xi@>>>\cc$ which takes $g$ to the image of
$\Ad(g\i)y$ under $\cc+\ft+\uup@>>>\cc,a+b+c\m a$. The inverse image of $\cf$
under this map is denoted again by $\cf$. On $\bby_\xi$ we have a free
$P$-action by right translation and $\cf$ is $P$-equivariant hence it descends
to a local system $\tcf$ on $\bby_\xi/P$. The group $Z_G(\xi)$ acts on 
$\bby_\xi/P$ by left translation and $\tcf$ is naturally a
$Z_G(\xi)$-equivariant local system. Then $\bZ_G(\xi)$ acts naturally on the 
cohomology 
$$\opl_nH^n_c(\bby_\xi/P,\tcf).\tag b$$
The set of irreducible representations (up to isomorphism) of $\bZ_G(\xi)$
which appear in the representation (b) is denoted by $\Irr_0\bZ_G(\xi)$.

\subhead 5.12\endsubhead
Let $\si,y$ be two elements of $\fg$ such that $\si$ is semisimple, $y$ is 
nilpotent and $[\si,y]=2r_0y$. We choose (as we may) elements $h,\ty$ in $\fg$
such that
$$[y,\ty]=h,[h,y]=2y,[h,\ty]=-2\ty,[\si,h]=0,[\si,\ty]=-2r_0\ty\tag a$$
and we set $\xi_1=\si-r_0h$. This is a semisimple element of $\fg$ (since
$\si,h$ are commuting semisimple elements) and it commutes with $y,h,\ty$. We
set $\xi=\xi_1+y$. If we make another choice $h',\ty'$ instead of $h,\ty$ then,
as it is known, there exists an element $g\in Z_G(\si,y)$ such that
$\Ad(g)h=h',\Ad(g)\ty=\ty'$. Let $\xi'_1=\si-r_0h'$. We have
$\Ad(g)\xi_1=\xi'_1$. Let $\xi'=\xi'_1+y$. We have $\Ad(g)\xi=\xi'$. Thus, the
$G$-orbit of $\xi$ is well defined by $\si,y$ (in fact, it depends only on the
$G$-orbit of $(\si,y)$.

Conversely, assume that $\xi\in\fg$ is given. We can write uniquely
$\xi=\xi_1+y$ where $\xi_1\in\fg$ is semisimple, $y\in\fg$ is nilpotent and
$[\xi_1,y]=0$. We choose (as we may) elements $h,\ty$ in $\fg$ such that
$$[y,\ty]=h,[h,y]=2y,[h,\ty]=-2\ty,[\xi_1,h]=0,[\xi_1,\ty]=0.\tag b$$
Let $\si=\xi_1+r_0h$. Then $\si$ is a semisimple element (since $\xi_1,h$ are 
commuting semisimple elements) and 
$[\si,y]=2r_0y,[\si,h]=0,[\si,\ty]=-2r_0\ty$. If we make another choice
$h',\ty'$ instead of $h,\ty$ then, as it is known, there exists 
$g'\in Z_G(\xi_1,y)$ such that $\Ad(g')h=h'$, $\Ad(g')\ty=\ty'$. Let 
$\si'=\xi_1+r_0h'$. We have $\Ad(g')\si=\si'$. Thus the $G$-orbit of $(\si,y)$
is well defined by $\xi$ (in fact, it depends only on the $G$-orbit of $\xi$.) 
Thus we have defined a bijection $\si,y\lra\xi$ between the set of $G$-orbits 
of pairs $(\si,y)$ of elements of $\fg$ such that $\si$ is semisimple, $y$ is 
nilpotent and $[\si,y]=2r_0y$, on the one hand, and the set of $G$-orbits in 
$\fg$, on the other hand.

\proclaim{Lemma 5.13} Assume that $(\si,y)$ corresponds to $\xi$ as above. Then
the groups $\bZ_G(\si,y)$ and $\bZ_G(\xi)$ may be naturally identified so that
$\Irr_0\bZ_G(\si,y)=\Irr_0\bZ_G(\xi)$.
\endproclaim
We may assume that $\si,y,\xi$ are related as follows: there exist $h,\ty$ in
$\fg$ so that 5.12(a) holds and $\xi=\xi_1+y$ where $\xi_1=\si-r_0h$. It is 
known that $Z'=\{g\in G;\Ad(g)y=y,\Ad(g)h=h,\Ad(g)\ty=\ty,\Ad(g)\si=\si\}$ is a
maximal reductive subgroup of $Z_G(\si,y)$ hence it has the same group of
components as $Z_G(\si,y)$. Similarly, since $Z_G(\xi_1)$ is connected,
reductive, $Z''=\{g\in G;\Ad(g)y=y,\Ad(g)h=h,\Ad(g)\ty=\ty,\Ad(g)\xi_1=\xi_1\}$
is a maximal reductive subgroup of 
$\{g\in G;\Ad(g)y=y,\Ad(g)\xi_1=\xi_1\}=Z_G(\xi)$ hence $Z''$ has the same 
group of components as $Z_G(\xi)$. Now $Z'=Z''$. It follows that $Z_G(\si,y)$ 
and $Z_G(\xi)$ have the same group of components. Note that
$$\bby_\xi=\{g\in G;\Ad(g\i)\xi_1\in\ft+\uup,\Ad(g\i)y\in\cc+\uup\}.\tag a$$
Now in the presence of the condition $\Ad(g\i)y\in\cc+\uup$, the conditions

(b) $\Ad(g\i)\xi_1\in\ft+\uup$,

(c) $\Ad(g\i)\xi_1\in\up$,
\nl
are equivalent. Indeed, it is clear that if (b) holds then (c) holds.
Conversely, assume that (c) holds. Then $\Ad(g\i)\xi_1=l\mod\uup$ where 
$l\in\ul$. By our assumption we have $\Ad(g\i)y=y_0\mod\uup$ where $y\in\cc$. 
Since $[\xi_1,y]=0$, we have $[\Ad(g\i)x_1,\Ad(g\i)y]=0$ and taking images 
under $\up@>>>\ul$ we deduce $[l,y_0]=0$. Since $y_0$ is distinguished in 
$\ul$, its centralizer in $\ul$ is $\ft$. Thus, $l\in\ft$ so that (b) holds. We
see that (a) can be rewritten as follows:
$$\bby_\xi=\{g\in G;\Ad(g\i)\xi_1\in\up,\Ad(g\i)y\in\cc+\uup\}.$$
Let $s\in\Hom(\bc^*,G)$ be such that the tangent map $\bc@>>>\fg$ of $s$ 
carries $1$ to $h$. We define a $\bc^*$-action of $\bby_\xi$ by 
$a:g\m s(a)g$. This induces a $\bc^*$-action on $\bby_\xi/P$ whose fixed point
set is
$$\bby'=\{g\in G;\Ad(g\i)\xi_1\in\up,\Ad(g\i)y\in\cc+\uup,\Ad(g\i)h\in\up\}/P.
$$
Similarly, we define a $\bc^*$-action on $\bx_{\si,y}$ by $a:g\m s(a)g$. This
induces a $\bc^*$-action on $\bx_{\si,y}/P$ whose fixed point set is
$$\bx'=\{g\in G;\Ad(g\i)y\in\cc+\uup,\Ad(g\i)\si\in\up,\Ad(g\i)h\in\up\}/P.$$
In the presence of the condition $\Ad(g\i)h\in\up$, the conditions
$\Ad(g\i)\si\in\up$ and $\Ad(g\i)\xi_1\in\up$ are equivalent (since 
$\xi_1=\si-r_0h$). It follows that $\bby'=\bx'$. Thus the $\bc^*$-actions on 
$\bby_\xi/P$ and on $\bx_{\si,y}/P$ have the same fixed point set. They also
have the same action of $Z'=Z''$.

The restriction of the local system $\tcf$ on $\bby_\xi/P$ (see 5.3) to this 
fixed point set is the same as the restriction of the local system $\tcf$ on 
$\bx_{\si,y}/P$ (see 5.11) to this fixed point set; the restriction is denoted
again by $\tcf$. By the principle of conservation of Euler characteristics by 
passage to the fixed point set of a torus action, we have
$$\align&\sum_n(-1)^nH^n_c(\bx_{\si,y}/P,\tcf)=\sum_n(-1)^nH^n_c(\bx',\tcf)\\&
=\sum_n(-1)^nH^n_c(\bby',\tcf)=\sum_n(-1)^nH^n_c(\bby_\xi/P,\tcf)\endalign$$
as virtual representations of $\bZ_G(\si,y)=\bZ_G(\xi)$. Hence for an 
irreducible representation $\rh$ of $\bZ_G(\si,y)=\bZ_G(\xi)$, the conditions

$\rh$ appears in $\sum_n(-1)^nH^n_c(\bx_{\si,y}/P,\tcf)$,

$\rh$ appears in $\sum_n(-1)^nH^n_c(\bby_\xi/P,\tcf)$
\nl
are equivalent. For $n$ odd we have $H^n_c(\bx_{\si,y}/P,\tcf)=0$ and 
$H^n_c(\bby_\xi/P,\tcf)=0$. (Indeed, both $\bx_{\si,y}/P$ and $\bby_\xi/P$ can
be regarded as fixed point sets of torus actions on the variety denoted by 
$\cp_y$ in \cite{\CU}, and it suffices to use the odd vanishing theorem 
\cite{\CU, 8.6} for $\cp_y$ together with \cite{\II, 4.4}). It follows that the
conditions 

$\rh$ appears in $\sum_n H^n_c(\bx_{\si,y}/P,\tcf)$,

$\rh$ appears in $\sum_nH^n_c(\bby_\xi/P,\tcf)$ 
\nl
are equivalent. The lemma is proved.

\subhead 5.14\endsubhead
Let $\fT(G,L,\cc,\cf)$ be the set consisting of all pairs $(\xi,\rh)$ (modulo
the natural action of $G$) where $\xi\in\fg$ and $\rh\in\Irr_0\bZ_G(\xi)$. By
5.12 and 5.13 we have a canonical bijection
$$\fS(G,L,\cc,\cf,r_0)\lra\fT(G,L,\cc,\cf).\tag a$$
Composing this with the bijection 5.3(c) we obtain a bijection
$$\Irr_{r_0}\bah(G,L,\cc,\cf)\lra\fT(G,L,\cc,\cf).\tag b$$
In the setup of 2.1 and assuming that $r_0\in\spa$, let 
$\fT^\spa(G,L,\cc,\cf)$ be the set of all $(\xi,\rh)$ in $\fT(G,L,\cc,\cf)$ 
such that $\xi\in\fg_\spa$. Then (a) restricts to a bijection
$$\fS^\spa(G,L,\cc,\cf,r_0)\lra\fT^\spa(G,L,\cc,\cf).\tag c$$
Indeed, let $\xi\in\fg$ and write $\xi=\xi_1+y$ where $\xi_1\in\fg$ is 
semisimple, $y\in\fg$ is nilpotent and $[\xi_1,y]=0$. Let $h,\ty$ in $\fg$ such
that 5.12(b) holds and let $\si=\xi_1+r_0h$. We must show that
$\xi\in\fg_\spa$ if and only if $\si\in\fg_\spa$. Clearly, $\xi\in\fg_\spa$ if
and only if $\xi_1\in\fg_\spa$. As in 5.8, we have $r_0h\in\fg_\spa$. (The 
eigenvalues of $h$ in any $V\in\ci_G$ are integers.) As in 5.8, if we have two
elements of $\fg_\spa$ that commute then their sum is again in $\fg_\spa$. 
Applying this to the commuting elements $\xi_1,r_0h^0$ and to the commuting 
elements $\si,-r_0h^0$ we deduce that $\si\in\fg_\spa$ if and only if 
$\xi_1\in\fg_\spa$. This yields (c).

Composing (c) with the bijection in 5.8 we obtain a bijection
$$\Irr_{r_0}^\spa\bah(G,L,\cc,\cf)\lra\fT^\spa(G,L,\cc,\cf).\tag d$$
(Notation of 5.8.) We can now reformulate Lemmas 5.9 and 5.10 as follows.

\proclaim{Lemma 5.15} Assume that $\t:\bc@>>>\br$ is a group homomorphism such
that $\t(r_0)\ne 0$. Let $\bam\in\Irr_{r_0}\bah(G,L,\cc,\cf)$ and let 
$(\xi,\rh)$ correspond to $\bam$ under 5.14(b). The following two conditions 
are equivalent:

(i) $\bam$ is $\t$-tempered;

(ii) $\xi\in\fg^{\Ker\t}$.
\endproclaim
By 5.9, condition (i) is equivalent to the condition that the semisimple part
$\xi_1$ of $\xi$ satisfies $\xi_1\in\fg^{\Ker\t}$. But this is clearly
equivalent to condition (ii).

\proclaim{Lemma 5.16} Assume that $G$ is semisimple. Assume that 
$\t:\bc@>>>\br$ is a group homomorphism such that $\t(r_0)\ne 0$. Let 
$\bam\in\Irr_{r_0}\bah(G,L,\cc,\cf)$ and let $(\xi,\rh)$ correspond to $\bam$
under 5.14(b). The following two conditions are equivalent:

(i) $\bam$ is $\t$-square integrable;

(ii) $\xi$ is a distinguished nilpotent element.
\endproclaim

\subhead 5.17\endsubhead
The local system on $\ex(\cc)$ (a unipotent class in $L$) that corresponds to 
$\cf$ under $\ex:\un L@>>>L$ is denoted again by $\cf$. For any $f\in G$, let
$$\dby_f=\{g\in G;g\i fg\in\ex(\cc)TU_P\}.\tag a$$
Consider the map $\dby_f@>>>\ex(\cc)$ which takes $g$ to the image of $g\i fg$
under $\ex(\cc)TU_P@>>>\ex(\cc),abc\m a$. The inverse image of $\cf$
under this map is denoted again by $\cf$. On $\dby_f$ we have a free $P$-action
by right translation and $\cf$ is $P$-equivariant hence it descends to a local 
system $\tcf$ on $\dby_f/P$. The group $Z_G(f)$ acts on $\dby_f/P$ by left 
translation and $\tcf$ is naturally a $Z_G(f)$-equivariant local system. Then 
$\bZ_G(f)$ acts naturally on the cohomology 
$$\opl_nH^n_c(\dby_f/P,\tcf).\tag b$$
The set of irreducible representations (up to isomorphism) of $\bZ_G(f)$ which
appear in the representation (b) is denoted by $\Irr_0\bZ_G(f)$.

In the setup of 2.1, let $\dfT^\spa(G,L,\cc,\cf)$ be the set of all $(f,\rh)$ 
(modulo the natural action of $G$) where $f\in G_\spa$ and 
$\rh\in\Irr_0\bZ_G(f)$.

If $f\in G_\spa$ corresponds to $\xi\in\fg_\spa$ under the bijection 5.1, we
have
$Z_G(f)=Z_G(\xi)$ hence $\bZ_G(f)=\bZ_G(\xi)$; we also have $\dby_f=\dby_\xi$.
(We use that $\ex:\up@>>>P$ restricts to a bijection
$\cc+\ft_\spa+\uup@>\sim>>\ex(\cc)T_\spa U_P$.) It follows that 
$\Irr_0\bZ_G(f)=\Irr_0\bZ_G(\xi)$. We see that $(\xi,\rh)\m(\ex(\xi),\rh)$ 
defines a bijection
$$\fT^\spa(G,L,\cc,\cf)@>\sim>>\dfT^\spa(G,L,\cc,\cf).\tag c$$
Composing this with the bijection 5.14(d) we obtain a bijection
$$\Irr_{r_0}^\spa\bah(G,L,\cc,\cf)\lra\dfT^\spa(G,L,\cc,\cf).\tag d$$

We can now state the following variant of Lemma 5.15.

\proclaim{Lemma 5.18} Assume that $\z:\bc^*@>>>\br$ is a group homomorphism 
such that $\z(\ex(r_0))\ne 0$. Let $\t=\z\ex:\bc@>>>\br$. Let 
$\bam\in\Irr_{r_0}^\spa\bah(G,L,\cc,\cf)$ and let $(f,\rh)$ correspond to 
$\bam$ under 5.17(d). The following two conditions are equivalent.

(i) $\bam$ is $\t$-tempered;

(ii) $f\in G^{\Ker\z}$.
\endproclaim
Using 5.15 we see that it is enough to verify the following statement: for
$\xi\in\fg$, we have $\xi\in\fg^{\Ker\t}$ if and only if 
$\ex(\xi)\in G^{\Ker\z}$. This is immediate. The lemma is proved.

\head 6. The subgroups $G_J$\endhead
\subhead 6.1\endsubhead
We fix an algebraic group $\hG$ such that $\hG^0$ is simply connected, almost
simple. We set $G=\hG^0$. We assume that we are given an element $\ta\in\hG$ of
finite order $d$ such that $G\T\bz/d\bz@>>>\hG$, $(g,j)\m g\ta^j$ is a 
bijection and such that the following holds: there exists a set of Chevalley 
generators $\{e_{i'},h_{i'},f_{i'}; i'\in I'\}$ for $\fg=\un{\hG}=\ug$ (with
standard notation) and a bijection $I'@>\sim>>I',i'\m{}^\ta i'$ of order $d$,
such that 

$\Ad(\ta)(e_{i'})=e_{{}^\ta i'},\Ad(\ta)(h_{i'})=h_{{}^\ta i'},
\Ad(\ta)(f_{i'})=f_{{}^\ta i'}$
\nl
for all $i'\in I'$. It follows that $\hG/G$ is a cyclic group of order $d$
generated by the image of $\ta$. Let $G^1$ be the connected component of $\hG$
that contains $\ta$. The subspace $\ft'$ of $\fg$ spanned by 
$\{h_{i'};i'\in I'\}$ is $\un T'$ for a maximal torus $T'$ of $G$. Let 
$\cy'=\Hom(\bc^*,T')$, $\cx'=\Hom(T',\bc^*)$. We have canonically 
$T'=\cy'\ot\bc^*$, $\ft'=\cy'_\bc$, $\ft'{}^*=\cx'_\bc$. Hence we may identify
$\cy'$ with a subgroup of $\ft'$ and $\cx'$ with a subgroup of $\ft'{}^*$. Let
$\chr'\sub\cy'\sub\ft'$ (resp. $R'\sub\cx'\sub\ft'{}^*$) be the set of coroots
(resp. roots) of $G$ with respect to $T'$. For $\a\in R'$ let $h_\a$ be the 
corresponding coroot and let $\fg_\a$ be the corresponding root subpace of 
$\fg$. For $i'\in I'$ define $\a_{i'}\in R'$ by $\fg_{\a_{i'}}=\bc e_{i'}$. 
Then $(R',\chr',\cx',\cy')$ is a root system with basis $\{\a_{i'};i'\in I'\}$.
Now $\ta$ normalizes $T'$. For any $j\in[0,d-1]$ we have a bijection 
$R'@>>>R',\a\m{}^{\ta^j}\a$ given by ${}^{\ta^j}\a(\ta^j t\ta^{-j})=\a(t)$ for
$t\in T,\a\in R'$, that is, $\Ad(\ta^j)\fg_\a=\fg_{{}^{\ta^j}\a}$ for 
$\a\in R'$.

\subhead 6.2\endsubhead
Let $G^\ta=Z_G(\ta),T=Z_{T'}(\ta)=T'\cap G^\ta$,
$\ft=\un T=\{x\in\ft';\Ad(\ta)x=x\}$. It is known that

(a) {\it $G^\ta$ is connected and $T$ is a maximal torus of $G^\ta$. Moreover,
$T'=Z_G(T)$; in particular, $N_G(T)\sub N_G(T')$.}

Let $W'=N_G(T')/T'$. Conjugation by $\ta$ induces an isomorphism $W'@>\sim>>W'$
whose fixed point set is denoted by $W'{}^\ta$. Let $W=N_{G^\ta}(T)/T$. It is 
known that

(b) {\it the obvious maps $N_G(T)/T'@<<<W@>>>W'{}^\ta$ are isomorphisms.}

Let ${}'\cy=\Hom(\bc^*,T),{}'\cx=\Hom(T,\bc^*)$. For $\b\in{}'\cx$ we set 
$\fg_\b=\{x\in\fg;\Ad(t)x=\b(t)x\quad\f t\in T\}$. Then 
$\fg=\opl_{\b\in{}'\cx}\fg_\b$ and $\fg_0=\ft'$. Let 
${}'R=\{\b\in{}'\cx-\{0\};\fg_\b\ne 0\}$. There is a unique subset ${}'\chr$ of
${}'\cy-\{0\}$, in bijection ${}'h_\b\lra\b$ with ${}'R$, such that 
$({}'R,{}'\chr,{}'\cx,{}'\cy)$ is a (not necessarily reduced) root system whose
associated Weyl group is $W$. We have canonically $T={}'\cy\ot\bc^*$, 
$\ft={}'\cy_\bc$, $\ft^*={}'\cx_\bc$. Hence we may identify ${}'\cy$ with a 
subgroup of $\ft$ and ${}'\cx$ with a subgroup of $\ft^*$ and we may regard 
${}'R\sub\ft^*,{}'\chr\sub\ft$.

\proclaim{Lemma 6.3} Define $\psi:R'@>>>\ft'{}^*$ by 
$\a\m\a+{}^\ta\a+\do+{}^{\ta^{d-1}}\a$. If $\a,\a'\in R'$ satisfy 
$\psi(\a)=\psi(\a')$ then $\a'={}^{\ta^j}\a$ for some $j\in[0,d-1]$.
\endproclaim
Let $R'_0$ be a set of representatives for the orbits of bijection 
$R'@>>>R',\a\m{}^\ta\a$. It is enough to show that $\a\m\psi(\a)$ is an
injective map $R'_0@>>>\ft'{}^*$. This can be easily checked in every case (we
may assume that $d\ge 2$).

\subhead 6.4\endsubhead
If $\a\in R'\sub\ft'{}^*$, then $\a|_\ft\in{}'R$. We thus obtain a map 
$R'@>>>{}'R,\a\m\a|_\ft$ which is constant on the orbits of $\a\m{}^\ta\a$. In
fact, using 6.3, we see that this map induces a bijection from the set of 
orbits of $\a\m{}^\ta\a$ on $R'$ onto ${}'R$. 

For $\b\in{}'R$ let $d'_\b$ be the cardinal of the corresponding orbit in $R'$;
thus $d'_\b=\dim\fg_\b$. For $\b\in{}'R$ we set $d''_\b=2$ if either
$2\b\in{}'R$ or $\fra{1}{2}\b\in{}'R$ and we set $d''_\b=1$ if 
$2\b\n{}'R,\fra{1}{2}\b\n{}'R$. We also set $d_\b=d'_\b d''_\b$. 

If $\a\in R'$ and $\b=\a|_\ft$, we have

${}'h_\b=h_\a$ if $\b\in{}'R-{}'R_*$ or if $d'_\b=1$,

${}'h_\b=h_\a+h_{{}^\ta\a}$ if $d'_\b=2,d''_\b=1$,

${}'h_\b=2h_\a+2h_{{}^\ta\a}$ if $d_\b=4$,

${}'h_\b=h_\a+h_{{}^\ta\a}+h_{{}^{\ta^2\a}}$ if $d_\b=3$.

Let $\bai'$ be the set of orbits of the bijection $I'@>>>I',i'\m{}^\ta i'$. For
$i\in\bai'$ let $\b_i=\a_{i'}|_\ft$ where $i'$ is any element of the orbit
$i$. Then $\{\b_i;i\in\bai'\}$ is a basis of the root system
$({}'R,{}'\chr,{}'\cx,{}'\cy)$.

Let $R$ be the subset of $\ft^*$ consisting of the vectors $d_\b\b$ for various
$\b\in{}'R$. For $\g\in R$ we set $h_\g=\fra{1}{d_\b}{}'h_\b$ where
$\g=d_\b\b,\b\in{}'R$. Let $\chr$ be the subset of $\ft$ consisting of the 
vectors $h_\g$ for various $\g\in R$. 

For $i\in\bai'$ let $d'_i=d'_{\b_i},d''_i=d''_{\b_i},d_i=d'_id''_i=d_{\b_i}$
and let $\g_i=d_i\b_i$. Let $\cy$ be the subgroup of $\ft$ generated by 
$\{h_{g_i};i\in\bai'\}$. Let $\cx$ be the set of all $\xi\in\ft^*$ that take
integer values on $\cy$. Then $R\sub\cx$, $\chr\sub\cy$ and $(R,\chr,\cx,\cy)$
is a (reduced) root system with Weyl group $W$ and with basis 
$\{\g_i;i\in\bai'\}$. It is also irreducible. (If $d=1$ we have $R'={}'R=R$. If
$d=2$ and $R'$ is of type $A_{2n-1}$, then ${}'R,R$ are of type $C_n,B_n$. If
$d=2$ and $R'$ is of type $A_{2n}$, then ${}'R,R$ are of type $BC_n,C_n$. If
$d=2$ and $R'$ is of type $D_n$, then ${}'R,R$ are of type $B_{n-1},C_{n-1}$. 
If $d=2$ and $R'$ is of type $E_6$, then ${}'R,R$ are of type $F_4,F_4$. If 
$d=3$ and $R'$ is of type $D_4$, then ${}'R,R$ are of type $G_2,G_2$.)

Let $\g_0\in R$ be the negative of the highest root of $R$ relative to 
$\{\g_i;i\in\bai'\}$. Then $\g_0=d_0\b_0$ for a unique $\b_0\in{}'R$ such 
that $2\b_0\n{}'R$. Here $d_0=d_{\b_0}=d$. Setting 

$I=\bai'\sqc\{0\}$, 
\nl
there are unique integers $n_i\in\bz_{>0} (i\in I)$ with $n_0=1$, such that
$$\sum_{i\in I}n_i\g_i=0$$
and $\sum_{i\in I}n_i$ is maximum possible (the Coxeter number of $W$).

For $i\in I$, we set $h_i=h_{\g_i}$. 

Let $V$ be a $\bc$-vector space with basis $\{b_i;i\in I\}$. Let $V'$ be the 
dual vector space with dual basis $\{b'_i;i\in I\}$. The canonical pairing
$V\T V'@>>>\bc$ is denoted by $x,x'\m x(x')$. We imbedd $\ft$ into $V'$ by
$y\m\sum_{i\in I}\g_i(y)b'_i$; we identify $\ft$ with its image, the subspace
$\{\sum_ic_ib'_i;c_i\in\bc,\sum_in_ic_i=0\}$ of $V'$. In particular, we may 
regard $h_i$ as a vector in $V'$ with $b_i(h_i)=2$. We have $b_i(x')=\g_i(x')$
for any $x'\in\ft,i\in I$. Let 
$$\ft^1=\{\sum_ic_ib'_i;c_i\in\bc,\sum_in_ic_i=1\}.$$
For $i\in I$ define $s_i:V@>>>V$ by $s_i(x)=x-x(h_i)b_i$ and its contragredient
$s_i:V'@>>>V'$ by $s_i(x')=x'-b_i(x')h_i$. Let $W^a$ be the subgroup of $GL(V)$
or $GL(V')$ generated by $\{s_i;i\in I\}$ (an affine Weyl group). Note that
$\ft,\ft^1$ are $W^a$-stable subsets of $V'$. We obtain a
homomorphism $W^a@>>>GL(\ft)$ whose image coincides with $W$.

\subhead 6.5\endsubhead
For any $S\sub I,S\ne\e$ let 
$$C_S=\{x'\in\ft^1; x'=\sum_{i\in S}c_ib'_i \text{ with } c_i\in\bc,c_i>0\quad
\f i\in S\}.$$ 
The sets $C_S$ are disjoint. Let 
$$C'=\cup_{S\sub I;S\ne\e}C_S.$$
For $J\sub I,J\ne I$, let $W_J$ be the subgroup of $W^a$ generated by 
$\{s_i;i\in J\}$ (a finite Coxeter group). For $S$ as above and $x'\in C_S$, we
have

(a) $\{w\in W^a;w(x')=x'\}=W_{I-S}$.

\proclaim{Lemma 6.6} Let $x'\in\ft^1$. The $W^a$-orbit $W^ax'$ meets $C'$ in
exactly one point.
\endproclaim
Let 
$V'_\br=\sum_{i\in I}\br b'_i,\ft_R=\ft\cap V'_\br,\ft^1_\br=\ft^1\cap V'_\br$.
 
The following $\br$-analogue of the lemma is well known.

(a) {\it Let $x'_1\in\ft^1_\br$. The $W^a$-orbit $W^ax'_1$ meets 
$C'\cap\ft^1_\br$ in exactly one point.}

We can write $x'=x'_1+\sqrt{-1}x'_2$ where $x'_1\in\ft^1_\br,x'_2\in\ft_\br$.
Using (a) we can find $w\in W^a$ and $S\sub I,S\ne\e$ such that 
$w(x')=x'_3+\sqrt{-1}x'_4$ where $x'_3\in C_S\cap\ft^1_\br,x'_4\in\ft_\br$. By
a well known property of Weyl chambers applied to $W_{I-S}$, we can find 
$w'\in W_{I-S}$ such that $w'(x_4)(b_i)\in\br_{\ge 0}$ for all $i\in I-S$. By 
6.5(a) we have $w'(x'_3)=x'_3$. Thus, $w'(w(x'))(b_i)$ is in 
$\br_{>0}+\sqrt{-1}\br$ if $i\in S$ and is in $\sqrt{-1}\br_{\ge 0}$ if 
$i\in I-S$. Hence $w'(w(x'))\in C_{S'}$ for some $S'$, $S\sub S'\sub I$. We see
that $W^ax'\cap C'\ne\e$.

Now let $z,z'$ be two points of $C'$ that $w(z)=z'$ for some $w\in W^a$. We can
write $z=z_1+\sqrt{-1}z_2$, $z'=z'_1+\sqrt{-1}z'_2$ where 
$z_1,z'_1\in C'\cap\ft^1_\br$, $z_2,z'_2\in\ft_\br$. Since $z_1,z'_1$ are in 
the same $W^a$-orbit, we see using (a) that $z_1=z'_1$. We can find
$S\sub I,S\ne\e$ such that $z_1=z'_1\in C_S\cap\ft^1_\br$. Since 
$z,z'\in C_S$, we have $z_2(b_i)\ge 0, z'_2(b_i)\ge 0$ for all $i\in I-S$. 
Moreover, we have $w(z_1)=z'_1=z_1$ (hence $w\in W_{I-S}$, see 6.5(a)) and 
$w(z_2)=z'_2$. Using a well known property of Weyl chambers applied to 
$W_{I-S}$, we deduce that $z_2=z'_2$. Thus, $z=z'$. The lemma is proved.

\subhead 6.7\endsubhead
Let $\cn=\{g\in G;gT\ta g\i=T\ta\}$, $\cn'=\cn\cap T'$. Clearly,
$N_{G^\ta}(T)\cap\cn'=T$.

\proclaim{Lemma 6.8} $\cn=N_{G^\ta}(T)\cn'=\cn'N_{G^\ta}(T)$.
\endproclaim
If $g\in\cn$ then $g$ normalizes the subgroup generated by $T\ta$ hence it also
normalizes the identity component $T$ of that subgroup. Since $T'=Z_G(T)$, it 
follows that $g$ normalizes $T'$. Thus, $\cn\sub N_G(T')$. It also follows that
$\cn'$ is normal in $\cn$ hence $N_{G^\ta}(T)\cn'=\cn'N_{G^\ta}(T)$.

Let $g\in\cn$. We set $a=g\i\ta g\ta\i,a'=\ta^2g\i\ta\i g\ta\i$. For $t\in T$
we have $\ta gt\ta g\i=gt\ta g\i\ta$ hence $ata'=t$. Taking $t=1$, we get 
$aa'=1$. Hence $ata\i=t$ for all $t\in T$. Thus, $a\in Z_G(T)=T'$. Thus,

(a) $\ta g\ta\i=ga$ for some $a\in T'$. 
\nl
Let $\bar g$ be the image of $g$ in $W'=N_G(T')/T'$. Then $\bar g\in W'{}^\ta$
(see (a) and 6.2). Since $W@>\sim>>W'{}^\ta$ (see 6.2(b)), there exists 
$g'\in N_{G^\ta}(T),t'\in T'$ such that $g=g't'$. We have $t'\in T'\cap\cn$ 
hence $t'\in\cn'$. Thus, $\cn=N_{G^\ta}(T)\cn'$. The lemma is proved.

\proclaim{Lemma 6.9}We have ${}'\cy\sub\cy$ and $\cy/{}'\cy$ may be identified
with $\cn'/T$.
\endproclaim
Let $F=\prod_{i\in\bai'}\bz/d'_i\bz$. Let $F'$ be the subgroup of 
$(\bc^*)^{\bai'}$ consisting of all $(a_i)$ such that $a_i^{d'_i}=1$ for all
$i\in\bai'$. Define $F@>\sim>>F'$ by $(l_i)\m(\exp(\k l_i/d'_i))$.

From the definitions we see that ${}'\cy$ has a $\bz$-basis
$\{\fra{1}{d''_i}{}'h_{\b_i};i\in\bai'\}$ and $\cy$ has a $\bz$-basis
$\{h_i;i\in\bai'\}$. Recall that $h_i=\fra{1}{d_i}h_{\b_i}$. It follows that
${}'\cy\sub\cy$ and we have 

$F@>\sim>>\cy/{}'\cy,(l_i)\m\text{ ${}'\cy$-coset of } \sum_{i\in\bai'}l_ih_i$.
\nl
By definition, $\cn'=\{t\in T';t\ta t\i\ta\i\in T\}$. The homomorphism 
$\chi:\cn'@>>>T,t\m t\ta t\i\ta\i$ with kernel $T$ induces an isomorphism 
$\cn'/T@>\sim>>\Im(\chi)$. Define $(\bc^*)^{\bai'}@>\sim>>T$ by 
$(a_i)\m\prod_{i\in\bai'}{}'h_{\b_i}(a_i)$. Via this isomorphism, $\Im(\chi)$
corresponds to $F'$. Combining the isomorphisms above yields the lemma.

\subhead 6.10\endsubhead
Let $p$ be the composition $\ft^1@>>>\ft@>>>T\ta$ where the first map is
$x\m x-b'_0$ and the second map is $x\m\exp_T(\k x)\ta$.

\proclaim{Lemma 6.11} The map $x'\m p(x')$ defines a bijection between $C'$ and
a set of representatives for the orbits of the $\cn$-action on $T\ta$ (by
conjugation).
\endproclaim
Now $x\m\exp_T(\k x)\ta$ induces $\ft/{}'\cy@>\sim>>T\ta$. Using 6.8, 6.9, we
see that via this isomorphism the action of $\cn$ on $T\ta$ corresponds to the
action of the obvious semidirect product of $W$ and $\cy/{}'\cy$ on 
$\ft/{}'\cy$ (with $\cy/{}'\cy$ normal) where the action of $W$ is the obvious
one and the action of $\cy/{}'\cy$ is by translation. It follows that we have 
an induced bijection

(a) $\{W-\text{orbits on } \ft/\cy\}\lra\{\cn-\text{orbits on } T\ta\}$.
\nl
It is well known that one may regard $\cy$ as a normal subgroup of $W^a$ in 
such a way that an element $y\in\cy$ acts on $\ft^1$ (as part of the 
$W^a$-action) in the same way as $y$ acts on $\ft^1$ by $x'\m x'+y$. Hence we
have an obvious bijection 

(b) 
$\{W^a-\text{orbits on } \ft^1\}\lra\{W^a/\cy-\text{orbits on } \ft^1/\cy\}$.
\nl
Now $\cy$ also acts on $\ft$ by translation and $x'\m x'-b'_0$ induces a 
bijection $\ft^1/\cy\lra\ft/\cy$ which is compatible with the action of 
$W^a/\cy$ on $\ft^1/\cy,\ft/\cy$. (This is because for $w\in W^a$ we have 
$w(b'_0)-b'_0\in\cy$.) Thus we have an obvious bijection

(c) $\{W^a/\cy-\text{orbits on } \ft^1/\cy\}\lra
\{W^a/\cy-\text{orbits on } \ft/\cy\}$.
\nl
By the last sentence in 6.4 we have

(d) $\{W^a/\cy-\text{orbits on } \ft/\cy\}=\{W-\text{orbits on } \ft/\cy\}$.
\nl
Combining (a)-(d), we obtain a bijection

$\{W^a-\text{orbits on } \ft^1\}\lra\{\cn-\text{orbits on } T\ta\}$.
\nl
This is induced by $x'\m p(x')$. We now use 6.6. The lemma follows.

\proclaim{Lemma 6.12}Let $Z$ be a semisimple $G$-conjugacy class in $G^1$. Then
$Z\cap(T\ta)$ is exactly one $\cn$-orbit in $T\ta$.
\endproclaim
This is classical when $d=1$. This is also known when $d>1$. It can be deduced
for example from \cite{\SE} (this reference deals with compact groups but our
case can be treated in a similar way).

Combining 6.11, 6.12, we have the following result.

\proclaim{Proposition 6.13} The map $x'\m p(x')$ defines a bijection between
$C'$ and a set of representatives for the $G$-conjugacy classes of semisimple
elements in $G^1$.
\endproclaim

\proclaim{Lemma 6.14} Let $\g\in R$ and $n\in\bz$. Let 
$H=\{y'\in\ft^1;\g(y'-b'_0)=n\}$. Let $S\sub I,S\ne\e$ and let $x'\in C_S$. If
$x'\in H$ then $\fra{1}{n_k}b'_k\in H$ for any $k\in S$.
\endproclaim
Let 

$V'_\br=\sum_{i\in I}\br b'_i,\ft_R=\ft\cap V'_\br,
\ft^1_R=\ft^1\cap V'_\br$,

$H_\br=H\cap\ft^1_R, H'=\{y'\in\ft_\br;\g(y')=0\}$.
\nl
The following $\br$-analogue of the lemma is well known.

(a) {\it Let $x'_1\in C_S\cap\ft^1_\br$. If $x'_1\in H_\br$ then 
$\fra{1}{n_k}b'_k\in H_\br$ for any $k\in S$.}
\nl
(This follows from the fact that $C_S\cap\ft^1_\br$ is a facet of a
configuration of reflection hyperplanes in $\ft^1_\br$ (one of which is
$H_\br$) and that $\fra{1}{n_k}b'_k (k\in S)$ are the vertices of that facet.)

We now write $x'=x'_1+\sqrt{-1}x'_2$ where $x'_1\in H_\br,x'_2\in H'$. We can 
find $S'\sub S$ such that

$x'_1=\sum_{i\in S'}c_{1,i}b'_i$, $x'_2=\sum_{i\in S}c_{2,i}b'_i$,
\nl
where 

(b) $c_{1,i}\in\br_{>0}$ for $i\in S'$, $c_{2,i}\in\br$ for $i\in S$,
$c_{2,i}\in\br_{>0}$ for $i\in S-S'$.
\nl
Using (a) for $x'_1$ we see that 

(c) $\fra{1}{n_k}b'_k\in H_\br$ for any $k\in S'$.
\nl
If $S'=S$, we are done. Assume now that $S-S'\ne\e$. From (c) we see that, for
$k\in S'$, we have $\g(b'_k-n_kb'_0)=nn_k$. Hence 

$\g(\sum_{i\in S'}c_{2,i}b'_i-\sum_{i\in S'}c_{2,i}n_ib'_0)=
\sum_{i\in S'}c_{2,i}nn_i$.

Since $x'_2\in H'$, we have $\sum_{i\in S}c_{2,i}n_i=0$ and 
$\g(\sum_{i\in S}c_{2,i}b'_i)=0$. It follows that
$$\align&\g(\sum_{i\in S'}c_{2,i}b'_i-\sum_{i\in S'}c_{2,i}n_ib'_0)\\&=
\g(-\sum_{i\in S-S'}c_{2,i}b'_i+\sum_{i\in S-S'}c_{2,i}n_ib'_0)
=-\sum_{i\in S-S'}c_{2,i}nn_i.\tag d\endalign$$
Set $c=\sum_{i\in S-S'}c_{2,i}n_i$. Since $S-S'\ne\e$ we see using (b) that
$c\in\br_{>0}$. From (d) we deduce

$\g(\sum_{i\in S-S'}c_{2,i}c\i b'_i-b'_0)=n$
\nl
so that $\sum_{i\in S-S'}c_{2,i}c\i b'_i\in H_\br\cap C_{S-S'}$. Using again 
(a) we deduce that $\fra{1}{n_k}b'_k\in H_\br$ for any $k\in S-S'$. Combining 
this and (c) we see that the lemma is proved.

\subhead 6.15\endsubhead
Let $\fN$ be the set of all pairs $(\b,j)$ where $\b\in{}'R,j\in[0,d-1]$ and

$j=0$ if $d'_\b=1,d''_\b=1$, 

$j=1$ if $d'_\b=1,d''_\b=2$.
\nl
We regard $\fN$ as a subset of the group $\ft^*\T\bz/d\bz$ by identifying
$[0,d-1]$ and $\bz/d\bz$ in the obvious way.

For $i\in I$ define $p_i$ by $p_i=0$ for $i\in I-\{0\}$, $p_0=1$.

\proclaim{Lemma 6.16} Let $S\sub I,S\ne\e$ and let $x'\in C_S$. The following 
two conditions for $(\b,j)\in\fN$ are equivalent:

(i) $\b(x'-b'_0)+\fra{j}{d}\in\bz$;

(ii) $(\b,j)=\sum_{i\in I-S}c_i(\b_i,p_i)$ with $c_i\in\bz$.
\endproclaim
Assume that $(\b,j)$ satisfies (ii). To show that it satisfies (i), we may
assume that $(\b,j)=(\b_i,p_i)$ for some $i\in I-S$. Then we have 
$\b_i(x'-b'_0)=d_i\i b_i(x'-b'_0)$ and this is $0$ if $i\ne 0$ and is 
$-d_0\i=-d\i$ if $i=0$. If $i\ne 0$ we have $p_i=0$ hence
$\b_i(x'-b'_0)+\fra{p_i}{d}=0$. If $i=0$ we have $p_i=1$ hence
$\b_i(x'-b'_0)+\fra{p_i}{d}=-\fra{1}{d}+\fra{1}{d}=0$, so that (i) holds.

Conversely, assume that $(\b,j)$ satisfies (i). Thus, we have
$\b(x'-b'_0)+\fra{j}{d}=n$ for some $n\in\bz$. We can write uniquely 
$\b=\sum_{k\in I-0}f_k\b_k$ where $f_k\in\bz$.

Let $H$ be the affine hyperplane $\{y'\in\ft^1;\b(y'-b'_0)+\fra{j}{d}=n\}$ in
$\ft^1$. Note that $H$ is of the form $\{y'\in\ft^1;\g(y'-b'_0)=n'\}$ for some
$\g\in R$ and $n'\in\bz$. Indeed, if $\b\in{}'R,d''_\b=1$, we can take 
$\g=d_\b\b$,$n'=d_\b n-j$; if $\b\in{}'R,2\b\in{}'R$, we can take 
$\g=4\b,n'=4n-2j$; if $\b\in{}'R,\fra{1}{2}\b\in{}'R$, we can take 
$\g=2\b,n'=2n-j$. Since $x'\in H$ and $x'\in C_S$ we have 
$\fra{1}{n_k}b'_k\in H$ for any $k\in S$. (See 6.14.) Thus, 
$\b(\fra{1}{n_k}b'_k-b'_0)+\fra{j}{d}=n$ for any $k\in S$.

If $k\in S,k\ne 0$, we have $\b(\fra{1}{n_k}b'_k-b'_0)=\fra{f_k}{d_kn_k}$.
Thus, $\fra{f_k}{d_kn_k}+\fra{j}{d}=n$. Since $\fra{jd_k}{d}\in\bz$, it follows
that $f_k=n_kg_k$ where $g_k\in\bz$, $g_k+\fra{jd_k}{d}=d_kn$.

If $0\in S$, we have $\fra{j}{d}=\b(0)+\fra{j}{d}=n$ so that $j=0$ and $n=0$.
In this case we deduce that for $k\in S,k\ne 0$ we have $g_k=0$ hence $f_k=0$ 
so that $\b=\sum_{k\in I-S}f_k\b_k$. Moreover, 
$(\b,j)=\sum_{k\in I-S}f_k(\b_k,p_k)$ since $j=0$.

Assume now that $0\in I-S$. Then
$$\align&\b=\sum_{k\in I-S-\{0\}}f_k\b_k
+\sum_{k\in S}n_k(d_kn-\fra{jd_k}{d})\b_k\\&
=\sum_{k\in I-S-\{0\}}f_k\b_k+\sum_{j\in I-0}n_k(d_kn-\fra{jd_k}{d})\b_k
-\sum_{k\in I-S-\{0\}}n_k(d_kn-\fra{jd_k}{d})\b_k\\&
=\sum_{k\in I-S-\{0\}}(f_k-n_k(d_kn-\fra{jd_k}{d}))\b_k-(nd-j)\b_0,\endalign$$
$$(\b,j)=\sum_{k\in I-S-\{0\}}(f_k-n_k(d_kn-\fra{jd_k}{d}))(\b_k,p_k)
-(nd-j)(\b_0,p_0)$$
since $j=-(nd-j)\mod d\bz$. The lemma is proved.

\subhead 6.17\endsubhead
Let $\b\in{}'R\cup\{0\}$. Then $\fg_\b$ is stable under $\Ad(\ta):\fg@>>>\fg$.
For any $j\in[0,d-1]$ we set 
$\fg_{\b,j}=\{x\in\fg_\b;\Ad(\ta)x=\exp(\k j/d)x\}$. Clearly,
$\fg_\b=\opl_{j\in[0,d-1]}\fg_{\b,j}$.

\proclaim{Lemma 6.18} Let $\b\in{}'R$ and let $j\in[0,d-1]$. We have 
$\dim\fg_{\b,j}=1$ if $(\b,j)\in\fN$ and $\fg_{\b,j}=0$ if $(\b,j)\n\fN$.
\endproclaim
We can assume that $d>1$. Assume first that $\dim\fg_\b=d$. Since $\fg_\b$ is
the direct sum of $d$ (one dimensional) root spaces of $\fg$ with respect to 
$T'$ which are cyclically permuted by $\Ad(\ta)$, it follows that 
$\Ad(\ta):\fg_\b@>>>\fg_\b$ has order $d$ and its $n$-th power has trace $0$
for $1\le n<d$. It follows that its $\exp(\k j/d)$-eigenspace is one 
dimensional for $0\le j<d$. Thus the lemma is proved in this case.

Next assume that $\dim\fg_\b=1$. Then $\b=\a|_\ft$ for a unique $\a\in R'$. We
may assume that $\a\in\sum_{i'\in I'}\bn\a_{i'}$. We can find a unique 
partition $I'=I'_1\sqc I'_2$ such that whenever $i'_1\in I'_1,i'_2\in I'_2$, 
the vertices $i'_1,i'_2$ of the Coxeter graph of $G$ are not joined. Let 
$x_0\in\fg_{\g_0}-\{0\}$. Let $B$ be the canonical basis of $\fg$ (as a left 
$\fg$-module) such that $x_0\in B$. Then $c_{i'}e_{i'}\in B$ for well defined 
$c_{i'}\in\bc^*$. Moreover, it is not difficult to check that there exist 
$a_1,a_2\in\bc^*$ such that $c_{i'}=a_1$ for all $i'\in I'_1$, $c_{i'}=a_2$ for
all $i'\in I'_2$ and $a_1+a_2=0$. Now $\Ad(\ta)(B)=uB$ for some $u\in\bc^*$. 
Hence $c_{i'}e_{{}^\ta i'}=uc_{{}^\ta i'}e_{{}^\ta i'}$ for all $i'$. We see
that $c_{i'}=uc_{{}^\ta i'}$. We consider two cases.

{\it Case 1.} Both $I'_1,I'_2$ are stable under $i'\m{}^\ta i'$. Then 
$a_1=ua_1$ and $u=1$. It follows that, if $x\in B\cap\fg_\b$, then 
$\Ad(\ta)x=ux=x$. Thus, if $\fg_{\b,j}\ne 0$ then $j=0$.

{\it Case 2.} $I'_1,I'_2$ are interchanged by $i'\m{}^\ta i'$ (hence $d=2$).
Then $ua_1=a_2=-a_1$ and $u=-1$. It follows that, if $x\in B\cap\fg_\b$, then
$\Ad(\ta)x=ux=-x$. Thus, if $\fg_{\b,j}\ne 0$ then $j=1$.

The lemma follows.

\subhead 6.19\endsubhead
Let $J\sub I,J\ne I$. Let 
$$\fg_J=\ft\opl\opl_{\b,j}\fg_{\b,j}$$
where $(\b,j)\in\fN$ is subject to the condition
$(\b,j)\in\sum_{i\in J}\bz(\b_i,p_i)$.

\proclaim{Lemma 6.20} There is a unique closed connected reductive subgroup 
$G_J$ of $G$ with Lie algebra $\fg_J$. If $x'\in C_{I-J}$, then 
$G_J=Z_G(p(x'))$.
\endproclaim
Recall that $p(x')=\exp_T(\k(x'-b'_0))\ta\in T\ta$. Now $Z_G(p(x'))$ is a 
closed connected reductive subgroup of $G$ whose Lie algebra is 
$\fh=\{x\in\fg;\Ad(p(x'))x=x\}$. Clearly, $\fh$ is stable under the $\Ad$ 
action of $T$ and that of $\ta$ on $\fg$. Hence $\fh$ is the sum of its 
intersections with the various $\fg_{\b,j}$ where $\b\in{}'R\cup\{0\}$. 
Clearly, $\fh\cap\fg_{0,j}$ is $0$ if $j\ne 0$ and is $\ft$ if $j=0$. Using 
6.18 we deduce that

$\fh=\ft\opl\opl_{\b,j}\fg_{\b,j}$
\nl
where $(\b,j)\in\fN$ is subject to $\b(\exp_T(\k(x'-b'_0)))\exp(\k j/d)=1$
(with $\b$ regarded as a character $T@>>>\bc^*$) or equivalently, to
$\b(x'-b'_0)+\fra{j}{d}\in\bz$ (with $\b$ regarded as a form $\ft@>>>\bc$). 
Using now 6.16 we see that $\fh=\fg_J$. The lemma is proved.

\subhead 6.21\endsubhead
For $i_1,i_2\in I$ let 

$a_{i_1,i_2}=\g_{i_2}(h_{i_1}), 
{}'a_{i_1,i_2}=\b_{i_2}({}'h_{\b_{i_1}})=\fra{d_{i_1}}{d_{i_2}}a_{i_1,i_2}$.
\nl
Then $(a_{i_1,i_2})$ is an untwisted affine Cartan matrix and 
$({}'a_{i_1,i_2})$ is a possibly twisted affine Cartan matrix. 

Let $J\sub I,J\ne I$. Let ${}'R_J$ the set of all $\b\in{}'R$ such that
$\b=w\b_i$ for some $i\in J$ and some $w\in W_J$. Let ${}'\chr_J$ be the set of
all ${}'h_\b$ where $\b\in{}'R_J$. 

\proclaim{Lemma 6.22} ${}'R_J$ (resp. ${}'\chr_J$) is exactly the set of roots 
(resp. coroots) of $G_J$ with respect to $T$.
\endproclaim
For $i\in J$ we have $(\b_i,p_i)\in\fN$ hence by 6.18, 6.19, $\b_i$ is a root
of $G_J$. For $i\in J$ there exists $g\in N_{G_J}(T)$ such that 
$\Ad(g):\ft@>>>\ft$ is a reflection that takes $\b_i$ to $-\b_i$. By 6.2(b),
there exists $g'\in N_{G^\ta}(T)$ and $t'\in T'$ such that $g=g't'$. Hence 
$\Ad(g):T@>>>T$ coincides with $\Ad(g'):T@>>>T$. Now there is a unique element
in $W$ that acts on $T$ as a reflection taking $\b_i$ to $-\b_i$, namely 
$s_i$. It follows that $\Ad(g)=s_i:\ft@>>>\ft$. Hence if $H_i\in\ft$ is the 
coroot of $G_J$ corresponding to $\b_i$, we have
$x-\b_i(x)H_i=x-\b_i(x){}'h_{\b_i}$ for all $x\in\ft$, hence 

(a) $H_i={}'h_{\b_i}$.
\nl
By 6.19, any root of $G_J$ is of the form $\b$ where $\b\in{}'R$ satisfies 
$\b\in\sum_{i\in J}\bz\b_i$. Thus, $(\b_i)_{i\in J}$ is a set of simple roots 
for $G_J$. By the first part of the argument, $W_J$ coincides with the Weyl 
group of $G_J$ (both are subgroups of $Aut(T)$); it follows that $R_J$ is
exactly the set of roots of $G_J$. The claim that ${}'\chr_J$ is exactly the 
set of coroots of $G_J$ follows from (a). The lemma is proved.

\subhead 6.23\endsubhead
From 6.22 we see that the Cartan matrix of $G_J$ is
$({}'a_{i_1,i_2})_{i_1,i_2\in J}$.

\subhead 6.24\endsubhead
Let $\hG_J=\cup_{j\in[0,d-1]}G_J\ta^j$. Since $\ta$ normalizes $G_J$, $\hG_J$ 
is a (closed) subgroup of $\hG$, with identity component $G_J$.
Let $Z_J$ (resp. $\hZ_J$) be the centre of $G_J$ (resp. of $\hG_J$). Let 
$\fz_J=\un{Z_J}$ (a subspace of $\ft$ hence a subspace of $V'$).

Let $V_J$ be the subspace of $V$ spanned by $\{b_i;i\in J\}$. Let $K=I-J$. Let
$V'_K$ be the subspace of $V'$ spanned by $\{b'_i;i\in K\}$, that is, the 
annihilator of $V_J$ in $V'$. We have $\fz_J=\ft\cap V'_K$ (a hyperplane in 
$\ft$). Let $\fz^1_J=\ft^1\cap V'_K$ (an affine hyperplane in $V'_K$).

Let $X_J=p(\fz^1_J)$. Note that $C_K\sub\fz^1_J$ hence $p(C_K)\sub X_J$.

\proclaim{Lemma 6.25} (a) We have $Z_J\sub\hZ_J$ and $Z_J^0=\hZ_J^0$. 

(b) $X_J$ is a connected component of $\hZ_J$.
\endproclaim
We prove (a). If $g\in Z_J$ then $g\in T$ hence $\ta g=g\ta$ hence $g\in\hZ_J$.
Thus, $Z_J\sub\hZ_J$. It follows that $Z_J^0\sub\hZ_J^0$. Now $\hZ_J^0\sub G_J$
hence $\hZ_J^0\sub Z_J$. Thus, $\hZ_J^0=G_J$.

We prove (b). Let $x'\in\fz_J^1$. We set $x^0=x'-b'_0\in\fz_J$. Let $i\in J$. 
Since $x^0\in V'_K$, we have $b_i(x^0)=0$. Now $p(x')$ acts on $\fg_{\b_i,p_i}$
by the scalar
$$\align&\b_i(\ex_T(\k x^0))\ex(\k d\i p_i)=
\ex(\k\b_i(x^0))\ex(\k d\i p_i)\\&=\ex(\k d_i\i\g_i(x^0)+\k d\i p_i)
=\ex(\k d_i\i b_i(x^0)+\k d\i p_i)\\&=\ex(-\k d_i\i b_i(b'_0)+\k d\i p_i)
=\ex(-\k d_i\i\de_{0,i}+\k d\i p_i))=\ex(0)=1.\endalign$$
It follows that $p(x')$ centralizes $G_J$. Since $p(x')\in T\ta$, it 
centralizes $\ta$ hence also $\hG_J$. Thus, $X_J\sub\hZ_J$. Clearly, $X_J$ is
connected. If $z\in Z^0_J$ we have $z=\ex_T(\k(x'_0))$ for some $x'_0\in\fz_J$
hence for $x'$ as above, $p(x'+x'_0)=p(x')z$. Thus, $X_J$ is stable by 
multiplication by $Z^0_J$. The lemma is proved.

\subhead 6.26\endsubhead
Let $G_1$ be the simply connected almost simple algebraic group corresponding
to the root system $(R,\chr,\cx,cy)$ (see 6.4). By 6.13 we have a natural
bijection between the set of $G_1$-conjugacy classes of semisimple elements in
$G_1$ and the set of $G$-conjugacy classes of semisimple elements in $G\ta$.

Our discussion of semisimple $G$-conjugacy classes in $G\ta$ has been 
influenced by \cite{\KA} where a connection between the elements of finite 
order in $G\ta$ and (possibly twisted) affine Lie algebras is given.

\head 7. The set $\fR(G\ta,G_J,\cc,\cf)$\endhead
\subhead 7.1\endsubhead
We preserve the setup of \S6. Let $J\sub I,J\ne I$ and let $K=I-J$. Let 
$G_{(J)}$ be the centralizer of $Z^0_J$ in $G$. Now $G_{(J)}$ is the subgroup 
of $G$ generated by $T'$ and by the root subgroups of $G$ corresponding to 
various $\a\in R'$ such that $\a|_\ft\in\sum_{i\in J}\bq\b_i$. We have 
$Z_{G_{(J)}}^0\cap G^\ta=Z_J^0$.

\proclaim{Lemma 7.2} Let $g\in Z_J$. Then we have $g=g_1g_2\ta^{-n}$ (in $\hG$)
where $g_2\in Z_G\cap G^\ta$, $n$ is an integer and $g_1$ is either the $n$-th 
power of an element in $X_J$ (if $n\ne 0$) or is an element of $Z_J^0$ (if 
$n=0$).
\endproclaim
We have $g\in T$ and $\b_i(g)=1$ for all $i\in J$. Hence $g=\ex_T(\k x)$ where
$x\in\ft$ satisfies $\b_i(x)\in\bz$ for all $i\in J$ hence
$x'=\sum_{i\in I}c_ib'_i$ with $c_i\in d_i\bz$ for all $i\in J$, $c_i\in\bc$ 
for all $i\in K$. We set $n=\sum_{k\in K}n_kc_k$. We have $n\in\bz$ since 
$n=-\sum_{i\in J}n_ic_i$. Let

$x'=\sum_{k\in K}c_kb'_k-nb'_0\in\ft,\quad
x''=\sum_{i\in J}c_ib'_i+nb'_0\in\ft$.
\nl
Then $x=x'+x''$. For any $i\in I-\{0\}$ we have
$\b_i(x'')=d_i\i b_i(\sum_{i_1\in J}c_{i_1}b'_{i_1}+nb'_0)$.

If $0\in J$ then $\b_i(x'')$ equals $d_i\i c_i$ if $i\in J-\{0\}$ and equals 
$0$ if $i\in K$. If $0\in K$ then $\b_i(x'')$ equals $d_i\i c_i$ if $i\in J$ 
and equals $0$ if $i\in K-\{0\}$. In any case case, $\b_i(x'')\in\bz$ for 
$i\in I-\{0\}$. It follows that $\ex_T(\k x'')$ is in the kernel of 
$\b_i:T@>>>\bc^*$ for any $i\in I-\{0\}$. Hence 
$g''=\ex_T(\k x'')\in Z_G\cap G^\ta$.

Assume first that $0\in J$. If $n=0$ we have $x'\in\fz_J$ hence 
$g'=\ex_T(\k x')\in Z_J^0$ and $g=g'g''$. If $n\ne 0$, we set 
$g'=\ex_T(\k\fra{1}{n}x')\ta$. Then $g'\in X_J$ and $g=g'{}^ng''\ta^{-n}$.

Assume next that $0\in K$. In this case we have $x'\in\fz_J$ hence 
$g'=\ex_T(\k x')\in Z_J^0$ and $g=g'g''$. The lemma is proved.

\proclaim{Lemma 7.3}If $g\in X_J$ then $Z_{G_{(J)}}(g)=Z_J$.
\endproclaim
Let $n=\sum_{k\in K}n_k$. We have $n\in\bz_{>0}$. If $x'=n\i\sum_{k\in K}b'_k$,
then $x'\in\ft^1$ and $z=p(x')\in X_J$. Since $x'\in C_{I-J}$, we have 
$Z_G(z)=G_J$ (see 6.20). Let $g\in X_J$. We have $g=zt$ where 
$t\in Z_J^0=Z_{G_{(J)}}^0$. Hence 
$Z_{G_{(J)}}(g)=Z_{G_{(J)}}(z)=Z_G(z)\cap G_{(J)}=G_J\cap G_{(J)}=G_J$. The 
lemma is proved.

\subhead 7.4\endsubhead
We fix $J,K$ as above, a nilpotent $G_J$-orbit $\cc$ in $\fg_J$ and an 
irreducible $G_J$-equivariant cuspidal local system $\cf$ (over $\bc$) on 
$\cc$. The local system on $\dcc=\ex(\cc)$ (a unipotent class in $G_J$) that
corresponds to $\cf$ under $\ex:\fg_J@>>>G_J$ is denoted again by $\cf$. Let
$$X_{(J)}=\cup_{g_1\in G_{(J)}}g_1X_J\dcc g_1\i,$$
(a locally closed subset of $G\ta$, stable under conjugacy by $G_{(J)}$).

(a) {\it Let $C$ be an $\Ad(G_{(J)})$-orbit in $G\ta$ that is contained in 
$X_{(J)}$. There exists an irreducible $G_{(J)}$-equivariant local system 
$\cf'$ on $C$ (unique up to isomorphism) such that the following holds: for any
$x\in X_J$ such that $x\dcc\sub C$, the restriction $\cf'|_{x\dcc}$ is the 
local system  obtained from $\cf$ via $\dcc@>\sim>>x\dcc$ (multiplication by 
$x$).}

This is shown as follows. Let $x\in X_J$ be such that $x\dcc\sub C$. There 
exists an irreducible $G_{(J)}$-equivariant local system $\cf(x)$ on $C$ 
(unique up to isomorphism) such that $\cf(x)|_{x\dcc}$ is the local system 
obtained from $\cf$ via $\dcc@>\sim>>x\dcc$ (multiplication by $x$). (We use 
the fact that, if $u\in\dcc$, then $Z_{G_{(J)}}(xu)=Z_{G_J}(u)$, see 7.3). We 
must only show that the isomorphism class of $\cf(x)$ is independent of $x$. If
$\card K=1$ then $X_J$ is a point and there is nothing to prove. Thus we may
assume that $\card K\ge 2$. 

Let $x'\in X_J$ be a second element such that $x'\dcc\sub C$. We must show that
$\cf(x),\cf(x')$ are isomorphic. It is enough to show that, if $f\in G_{(J)}$ 
is such that $fxf\i=x'$ (so that $\Ad(f)Z_J=Z_J$, as we see from 7.3), then 
$\Ad(f)$ carries $(\dcc,\cf)$ to $(\dcc,\cf)$. If $d\ge 2$, then $(\dcc,\cf)$ 
is uniquely determined by $G_J$ for $\card K\ge 2$ (see the tables in \S11) and
we are done. Assume now that $d=1$. Now $Z_J$ acts on $\cf$ through a character
$\chi:Z_J@>>>\bc^*$. Since there is at most one pair consisting of a unipotent
class of $G_J$ and an irreducible $G_J$-equivariant cuspidal local system on it
with prescribed action of $Z_J$, it is enough to show that 
$\chi(z)=\chi(fzf\i)$ for any $z\in Z_J$. Since $\chi$ is trivial on $Z^0_J$ it
is enough to show that $fzf\i=z\mod Z^0_J$ for any $z\in Z_J$. Using 7.2, we 
can write $z$ in the form $z=z_0z_1^nz_2$ where $z_0\in Z_J^0,z_1\in X_J$,
$z_2\in Z_G,n\in\bz$. (In this case the power of $\ta$ in 7.2 is $1$.) It 
suffices to show that $fz_0z_1^nz_2f\i=z_0z_1^nz_2\mod Z_J^0$ or that
$fz_0z_1^nf\i=z_0z_1^n\mod Z_J^0$, or that $fz_1^nf\i=z_1^n\mod Z_J^0$. (We use
that $fz_0=z_0f$.) We have $z_1=xa$ where $a\in Z^0_J$ hence

$fz_1^nf\i=fx^na^nf\i=fx^nf\i fa^nf\i=x'{}^na^n=z_1^n(x\i x')^n=z_1^n\mod
Z^0_J$
\nl
since $x\i x'\in Z_J^0$. Thus, (a) is established.

\subhead 7.5\endsubhead
We can find a parabolic subgroup $P$ of $G$ which has $G_{(J)}$ as a Levi 
subgroup and satisfies $\ta P\ta\i=P$. (We can choose a general enough
$y\in\Hom(\bc^*,Z^0_J)$ such that, setting 
$\fg(n)=\{x\in\fg;\Ad(y(a))x=a^nx\quad\f a\in\bc^*\}$ for $n\in\bz$, we have
$\fg(0)=\un{G_{(J)}}$. Then $\opl_{n\in\bn}\fg(n)=\un P$ for a well defined $P$
which satisfies our requirements.) For any $f\in G\ta$ we set
$$\bu_f=\{g\in G;g\i fg\in X_{(J)}U_P\}.$$
We have an obvious map $\pi:\bu_f@>>>X_{(J)}$ which takes $g$ to the image of
$g\i fg$ under $X_{(J)}U_P@>>>X_{(J)},ab\m a$.

The image of $\pi$ is a disjoint union of finitely many $\Ad(G_{(J)}$-orbits in
$X_{(J)}$ (since the semisimple part of a point in this image is contained in a
fixed $\Ad(G)$-orbit, namely that of the semisimple part of $f$). This image 
carries a $G_{(J)}$-equivariant local system (on each connected component we 
take the local system in 7.4(a)). Taking inverse image under 
$\bu_f@>>>\pi(\bu_f)$ of this local system we obtain a local system on $\bu_f$
which is $P$-equivariant for the free $P$-action on $\bu_f$ given by right 
translation, hence it descends to a local system on $\bu_f/P$ denoted by
$\tcf$. Now $Z_G(f)$ acts on $\bu_f/P$ by left translation and $\tcf$ is 
naturally a $Z_G(f)$-equivariant local system. Then $\bZ_G(f)$ acts naturally 
on the cohomology 
$$\opl_nH^n_c(\bu_f/P,\tcf).\tag a$$
The set of irreducible representations (up to isomorphism) of $\bZ_G(f)$ which
appear in the representation (a) is denoted by $\Irr_1\bZ_G(f)$.

Let $\fR(G\ta,G_J,\cc,\cf)$ be the set of all $(f,\rh)$ (modulo the 
$\Ad$-action of $G$) where $f\in G\ta$ and $\rh\in\Irr_1\bZ_G(f)$. 

\proclaim{Lemma 7.6} Assume that $S\sub K,S\ne\e$. Let $P'$ be a parabolic 
subgroup of $G$ which has $G_{(J)}$ as a Levi subgroup and satisfies 
$\ta P'\ta\i=P'$. Then

(a) $G_{I-S}\cap G_{(J)}=G_J$;

(b) $G_{I-S}\cap P'$ is a parabolic subgroup of $G_{I-S}$ with Levi subgroup
$G_J$;

(c) $G_{I-S}\cap U_{P'}=U_{G_{I-S}\cap P'}$.
\endproclaim
The proof is routine. It will be omitted.

\head 8. Geometric affine Hecke algebras\endhead
\subhead 8.1\endsubhead
We preserve the setup of 7.4. In this and the next subsection we assume that 
$\card(K)\ge 2$.

As in \cite{\IM, 5.6}, for any $J'\sub I$ such that $J\sub J'\ne I$, 
conjugation by the longest element $w_0^{J'}$ of $W_{J'}$ leaves stable 
$\{s_i;i\in J\}$. Hence for any $k\in K$ we have 
$w_0^{J\cup k}w_0^Jw_0^{J\cup k}=w_0^J$ and 
$\si_k=w_0^{J\cup k}w_0^J=w_0^Jw_0^{J\cup k}$ is an involution. Now $\si_k$
preserves the subspace $V_J$ of $V$ hence also the subspace $V'_K$ of $V'_K$.
Hence the subgroup $W^*$ of $W^a$ generated by $\{\si_k;k\in K\}$ acts on 
$V'_K$. As in \cite{\IM, 2.11}, $W^*$ is a Coxeter group (an affine Weyl 
group).

Let $y_0\in\cc$. For $k\in K$, $\ad(y_0):\fg@>>>\fg$ induces a nilpotent
endomorphism of $\fg_{J\cup k}/\fg_J$. Let $\uc_k$ be the largest integer 
$\ge 2$ such that the $(\uc_k-2)$-th power of this nilpotent endomorphism is
non-zero. (This does not depend on the choice of $y_0$.)

The $W^*$-action on $V'_K$ leaves stable the subset $\fz_J^1$ of $V'_K$. It 
also leaves stable the subspace $\fz_J$ of $V'_K$ where it acts through a 
finite quotient $\cw\sub GL(\fz_J)$. Let $\cl'$ be the set of all $x\in\fz_J$ 
such that the translation $z\m z+x$ of $\fz_J^1$ coincides with the 
automorphism $x\m w(x)$ of $\fz_J^1$ for some $w\in W^*$. Then $\cl'$ is a 
subgroup of $\fz_J$ such that $\cl'_\bc=\fz_J$. For $k\in K$ there exist 
non-zero vectors $\tih_k\in\fz_J$ and $\tig_k\in\Hom(\fz_J,\bc)$ such that
$\si_k(x)=x-\tig_k(x)\tih_k$ for all $x\in\fz_G$ and $\tig_k(\tih_k)=2$. These
vectors are uniquely determined if we require that $\tig_k(x)=z_kx(b_k)$ for 
all $x\in\fz_J$, $\tih_k\in\cl$ and $z_k\in\bz_{>0}$ is maximum possible (see 
\cite{\IM, 2.11}). We have $z_k\in\{1,2,3,4\}$. Let 
$\cl=\{x\in\fz_J^*;x(\cl')\in\bz\}$. Then $\tig_k\in\cl$. Clearly, $\cw$ acts 
naturally on $\cl',\cl$. Let $\tcr$ (resp. $\tcr\che{}$) be the set of vectors
in $\fz_J^*$ (resp. $\fz_J$) that are of the form $w(\tig_k)$ (resp. 
$w(\tih_k)$) for some $w\in\cw$ and some $k\in K$. Then 

(a) $(\tcr,\tcr\che{},\cl,\cl')$ {\it is an irreducible root system with Weyl 
group $\cw$.}
\nl
See \cite{\IM, 2.11}. Moreover, $\{\tih_k;k\in K\}$ generates $\cl'$ (see
\cite{\IM, 2.14}). Now  $\{\tig_k;k\in K\}$ spans $\cl_\bc$ over $\bc$ with a 
single relation $\sum_{k\in K}\tn_k\tig_k=0$ where $\tn_k\in\bz_{>0}$ for all 
$k$ (at least one $\tn_k$ is $1$) and $\sum_{k\in K}\tn_k$ is the Coxeter 
number of (a). For $k\in K$ we have 
$$z_k=\fra{n_k}{\tn_k}.$$ 
We define a subset $K^\fla$ of $K$ as follows. If $W^*$ is a Coxeter group of
type $\ti C_n,n\ge 1$, and $k,k'$ correspond to the two ends of the Coxeter 
graph then $K^\fla=\{k,k'\}$. In any other case, $K^\fla=\e$. We set 
$K^\sh=K-K^\fla$.

For $k\in K$ we define $\baz_k$ by $\baz_k=z_k/2$ if $k\in K^\fla$ and
$\baz_k=z_k$ if $k\in K^\sh$.

We have $\baz_k\in\{\fra{1}{2},1,2,3\}$. We set 
$$\hg_k=(\baz_k/z_k)\tig_k,\hah_k=(z_k/\baz_k)\tih_k.$$
We have $\hg_k\in\cl$. Let $\car$ (resp. $\car\che{}$) be the set of vectors in
$\fz_J^*$ (resp. $\fz_J$) that are of the form $w(\hg_k)$ (resp. $w(\hah_k)$) 
for some $w\in\cw$ and some $k\in K$. 

In the case where $K^\fla\ne\e$, we have $\{k\in K;\tn_k=1\}=K^\fla$ and we
choose $k_0\in K$ so that $k_0\in K^\fla$ and 
$\uc_{k_0}\baz_{k_0}d_{k_0}\le\uc_{k'}\baz_{k'}d_{k'}$ where 
$K^\fla=\{k_0,k'\}$. In the case where $K^\fla=\e$ we choose $k_0\in K$ such
that $\tn_{k_0}=1$. One can verify that 
$$\dim\fg_{I-\{k_0\}}\ge\dim\fg_{I-\{k\}} \text{ for any } k\in K.$$
Now $\{\tig_k;k\in K-\{k_0\}\}$ is a basis for 
$(\tcr,\tcr\che{},\cl,\cl')$ and $\tig_{k_0}$ is the negative of the highest
root of $(\tcr,\tcr\che{},\cl,\cl')$ with respect to this basis. Moreover,

(d) $(\car,\car\che{},\cl,\cl')$ {\it is a root system and 
$\Pi=\{\hg_k;k\in K-\{k_0\}\}$ is a basis for it. Its Weyl group is $\cw$.}

See \cite{\IM, 2.15}. Since $\{\tih_k;k\in K\}$ generates $\cl'$  we see that 
$\car\che{}\cup((1/2)\car\che{}\cap\cl')$ generates $\cl'$.

\subhead 8.2\endsubhead
If $k\in K^\fla$ and $K^\fla-\{k\}=\{k'\}$ we set
$$\la(\hg_k)=(\uc_k\baz_kd_k+\uc_{k'}\baz_{k'}d_{k'})/2,\quad
\la^*(\hg_k)=|\uc_k\baz_kd_k-\uc_{k'}\baz_{k'}d_{k'}|/2.\tag a$$
If $k\in K^\sh$, we set 
$$\la(\hg_k)=\uc_k\baz_kd_k/2.\tag b$$
Restricting, we get functions $\la:\Pi@>>>\bn$,
$\la^*:\{\hg_k\in\Pi;\hah_k\in 2\cl'\}@>>>\bn$. (One can check that $\la$ and
$\la^*$ have indeed values in $\bn$.) Then $(\la,\la^*)$ is a parameter set for
the root system 7.6(d) with its basis $\Pi$. Hence the $\bc[v,v\i]$-algebra
$H^{\la,\la^*}_{\car,\cl}$ is well defined as in 1.2. We denote this algebra by
$H(G\ta,G_J,\cc,\cf)$; we call it a {\it geometric affine Hecke algebra}.

\subhead 8.3\endsubhead
Assume now that $\card(K)=1$. Then $\fz_J=0$. We set
$(\car,\car\che{},\cl,\cl')=(0,0,\e,\e)$. (A root system.) Then
$H^{\la,\la^*}_{\car,\cl}=\bc[v,v\i]$ is again denoted by \linebreak
$H(G\ta,G_J,\cc,\cf)$.

\subhead 8.4\endsubhead
The definition of $K^\fla$ given in 8.1 differs slightly from the one given in
\cite{\IM, 2.13}. (There are only three cases where the definitions differ:
those in \cite{\IM, 7.16, 7.47, 7.56}.) In the context of \cite{\IM} it does 
not matter which of the two definitions we adopt. They both lead to the same 
$H(G\ta,G_J,\cc,\cf)$. However, in the more general context of this paper, the
present definition should be adopted (it diverges from the definition in
\cite{\IM}). 

\head 9. A bijection\endhead
\subhead 9.1\endsubhead
We preserve the setup of \S6,7,8. In particular we fix $J,K$ as in 7.1,
$\cc,\cf$ as in 7.4, $P$ as in 7.5, $k_0$ as in 8.1. Let $v_0\in\bc^*$ be such
that either $v_0=1$ or $v_0$ is not a root of $1$. We choose $r_0\in\bc$ such 
that $\ex(r_0)=v_0$; if $v_0=1$ we choose $r_0=0$. Let $\spa$ be the 
$\bq$-subspace of $\bc$ spanned by $r_0$. We have $\spa\cap\k\bq=0$. We can 
choose a $\bq$-subspace $\di$ of $\bc$ such that
$$\k\bq\sub\di,\di\opl\spa=\bc.\tag a$$
We have $\bc^* =\ex(\di)\T\ex(\spa)$.

If $G'$ is an algebraic group, any $g\in G'$ can be written uniquely in the 
form $g=g_\di g_\spa=g_\spa g_\di$ where $g_\di\in G'_\di$ is semisimple and 
$g_\spa\in G'_\spa$.

In this section we will define a bijection
$$\Irr_{v_0}H(G\ta,G_J,\cc,\cf)\lra\fR(G\ta,G_J,\cc,\cf)\tag b$$
in terms of $r_0$ and $\di$ as above.

If $\card(K)=1$ then both sides of (b) consist of one element, hence there is a
unique bijection between them. In the remainder of this section we assume that
$\card(K)\ge 2$. 

\subhead 9.2\endsubhead
Let $V'{}^\di=\{x'\in V';x'=\sum_{i\in I}c_ib'_i, c_i\in\k\i\di\}$.

For any $S\sub I,S\ne\e$ let $C_S^\di=C_S\cap V'{}^\di$. Let 
$C'{}^\di=\sqc_{S\sub I,S\ne\e}C_S^\di=C'\cap V'{}^\di$. 
Let $D=\sqc_{S\sub K;S\ne\e}C_S$, $D^\di=\sqc_{S\sub K;S\ne\e}C_S^\di$

(a) {\it The map $x'\m p(x')$ defines a bijection between $C'{}^\di$ and a set
of representatives for the orbits of the $\cn$ action on 
$T_\di\ta=\hG_\di\cap T\ta$ (by conjugation).}
\nl
This is an immediate consequence of 6.11 and its proof.

We set 
$$\ct=\cl'\ot\bc^*.$$
(A torus.) Then $\cl'_\bc=\fz_J=\uct$ hence $\ex_\ct:\fz_J@>>>\ct$ is defined.
Let $p'$ be the composition $\fz_J^1@>>>\fz_J@>>>\ct$ where the first map is 
$x\m x-\fra{1}{n_{k_0}}b'_{k_0}$ and the second map is $x\m\ex_\ct(\k x)$.

(b) {\it The map $x'\m p'(x')$ defines a bijection between $D$ and a set of 
representatives for the $\cw$-orbits in $\ct$ and also a bijection between 
$D^\di$ and a set of representatives for the $\cw$-orbits in $\ct_\di$.}

Just like 6.11 was proved using 6.6, the proof of (b) is based on the following
analogue of 6.6. 

{\it Let $x'\in\fz_J^1$. The $W^*$-orbit $W^*x'$ meets $D$ in exactly one 
point.} 
\nl This is proved exactly like 6.6, by replacing $W^a$ by $W^*$.

For $\dd\in D^\di$, let $\cw_\dd$ be the stabilizer of $p'(\dd)$ in $\cw$. Let
$\Th$ be the set consisting of all pairs $(\dd,\de)$ where $\dd\in D^\di$ and 
$\de$ is a $\cw_\dd$-orbit in $\ct_\spa$. Using the decomposition 
$\ct=\ct_\di\ct_\spa$ we see that there is a bijection $\Th@>\sim>>\ct/\cw$ 
which associates to $(\dd,\de)\in\Th$ the $\cw$-orbit of $p'(\dd)z$ where 
$z\in\de$. (This $\cw$-orbit is denoted by $\Si_{\dd,\de}$.) 

\subhead 9.3\endsubhead
The partition 1.3(a) becomes in our case
$$\Irr_{v_0}H^{\la,\la^*}_{\car,\cl}=
\sqc_{S\sub K,S\ne\e}\sqc_{\dd\in C_S^\di}\sqc_{\de\in\ct_\spa/\cw_\dd}
\Irr_{\Si_{\dd,\de},v_0}H^{\la,\la^*}_{\car,\cl}.\tag a$$
Let $(\dd,\de)\in\Th$ where $\dd\in C_S^\di$. Then $p'(\dd)\de$ is a 
$\cw_\dd$-orbit contained in $\Si_{\dd,\de}$. Moreover, 
$$\car_{p'(\dd)\de},\car\che{}_{p'(\dd)\de},\Pi_{p'(\dd)\de},\la_{p'(\dd)\de},
\la^*_{p'(\dd)\de}$$
are defined in terms of $\cl,\cl',\car,\car\che{},\Pi,p'(\dd)\de,\la,\la^*$ in 
the same way as \linebreak $R_c,\chr_c,\Pi_c,\la_c,\la^*_c$ were defined in 
[8.1]${}'$ (see 2.2) and 3.1 in terms of \linebreak
$X,Y,R,\chr,\Pi,c,\la,\la^*$. Then, by 3.2 we have a natural bijection 
$$\Irr_{\Si_{\dd,\de},v_0}H^{\la,\la^*}_{\car,\cl}\lra\Irr_{p'(\dd)\de,v_0}
H^{\la_{p'(\dd)\de},\la^*_{p'(\dd)\de}}_{\car_{p'(\dd)\de},\cl}.\tag b$$
Let $\cw_{K-S}$ be the subgroup of $\cw$ generated by the image of
$\{\si_k;k\in K-S\}$ in $\cw$. Let $\car_{K-S}$ (resp. $\car\che{}_{K-S}$) be 
the set of vectors of $\cl$ (resp. $\cl'$) of the form $w(\hg_k)$ (resp.
$w(\hah_k)$) for some $k\in K-S,w\in\cw_{K-S}$. Let 
$\Pi_{K-S}=\{\hg_k;k\in K-S\}$. Then
$$(\cl,\cl',\car_{p'(\dd)\de},\car\che{}_{p'(\dd)\de},\Pi_{p'(\dd)\de},\cw_\dd)
=(\cl,\cl',\car_{K-S},\car\che{}_{K-S},\Pi_{K-S},\cw_{K-S}).$$
(The main assertion here is that $\car_{p'(\dd)\de}=\car_{K-S}$. In other 
words, for $\a\in\car$, with corresponding $\check\a\in\car\check{}$, the 
condition that $\a\in\car_{K-S}$ is equivalent to the condition that 
$\a(\dd-\fra{1}{n_{k_0}}b'_{k_0})$ is in $\bz$ if $\check\a\n 2\cl'$ and is in
$\fra{1}{2}\bz$ if $\check\a\in 2\cl'$. This is an assertion of the same type 
as 6.16 and has a similar proof. See also \cite{\IM, 3.9, 3.10}.)

Thus, the bijection (b) can be rewritten as
$$\Irr_{\Si_{\dd,\de},v_0}H^{\la,\la^*}_{\car,\cl}\lra
\Irr_{p'(\dd)\de,v_0}H^{\la,\la^*}_{\car_{K-S},\cl};$$
the exponents $\la,\la^*$ in both sides are restrictions of the function given
by 8.2(a),(b).

Taking union of all $d,\de$ and composing with (a) we obtain a bijection
$$\Irr_{v_0}H^{\la,\la^*}_{\car,\cl}\lra
\sqc_{S\sub K,S\ne\e}\sqc_{\dd\in C_S^\di}\sqc_{\de\in\ct_\spa/\cw_\dd}
\Irr_{p'(\dd)\de,v_0}H^{\la,\la^*}_{\car_{K-S},\cl}.\tag c$$
Let $(\dd,\de)\in\Th$. Then $(\cl_\bq,\cl'_\bq,\car_{K-S},\car\che{}_{K-S})$ is
a $\bq$-root system with basis $\Pi_{K-S}$. 

For $k\in K$ we regard $\tig_k$ as a character $\ct@>>>\bc^*$ given by
$l'\ot a\m a^{\tig_k(l')}$ where $l'\in\cl',a\in\bc^*$. Then for $z\in\fz_J$ 
we have $\g_k(\ex_\ct(z))=\ex(\g_k(z))$. We show that, if
$k\in K-S,k\in K^\fla$ then:

$$\hg_k(p'(\dd))=-1 \text{ if } k=k_0 \text{ and } \hg_k(p'(\dd))=1 \text{ if }
k\ne k_0.$$
Indeed, $\hg_k(\ex_\ct(\k(\bod-\fra{1}{n_{k_0}}b'_{k_0})))=
\ex(\k\baz_k b_k(\bod-\fra{1}{z_{k_0}}b'_{k_0}))=\ex(-\k\de_{k,k_0}/2).$

It follows that, if we define $\mu:\Pi_{K-S}@>>>\bz$ in terms of $\la,\la^*$ as
in 4.1 (with $t_0=p'(\dd)$), then
$$\mu(\hg_k)=d_k\baz_k\uc_k.$$   
Define a $\cw_{K-S}$-orbit $\bde$ in $\uct_\spa$ by $\ex_\ct(\bde)=\de$. By 4.2
we have a bijection
$$\Irr_{\bde,r_0}\bah^\mu_{\car_{K-S},\cl_\bq}\lra
\Irr_{p'(\dd)\de,v_0}H^{\la,\la^*}_{\car_{K-S},\cl}.$$
Taking union over all $\dd,\de$ and composing with (c) we obtain a bijection
$$\Irr_{v_0}H^{\la,\la^*}_{\car,\cl}\lra
\sqc_{S\sub K,S\ne\e}\sqc_{\dd\in C_S^\di}\sqc_{\bde\in\uct_\spa/\cw_{K-S}}
\Irr_{\bde,r_0}\bah^\mu_{\car_{K-S},\cl_\bq}.\tag d$$

For $k\in K$ we set ${}^*\g_k=(1/d_k\baz_k)\hg_k\in\fz_J^*$ (that is, 
${}^*\g_k$ is the restriction of $\b_k$ to $\fz_J$ where $\b_k$ is regarded as
an element of $\ft^*$). We set ${}^*h_k=d_k\baz_k\hah_k\in\fz_J$. 

Let $S\sub K,S\ne\e$. Let ${}^*\car_{K-S}$ (resp. ${}^*\car\che{}_{K-S}$) be 
the set of vectors of $\cl_\bq$ (resp. $\cl'_\bq$) of the form $w({}^*\g_k)$ 
(resp. $w({}^*h_k)$) for some $k\in K-S,w\in\cw_{K-S}$. Let 
${}^*\Pi_{K-S}=\{{}^*\g_k;k\in K-S\}$. Then
$(\cl_\bq,\cl'_\bq,\car_{K-S},\car\che{}_{K-S},{}^*\Pi_{K-S})$ is a $\bq$-root
system. Define $\mu':{}^*\Pi_{K-S}@>>>\bn$ by $\mu'({}^*\g_k)=\uc_k$. There is
an algebra isomorphism
$$\bah^\mu_{\car_{K-S},\cl_\bq}@>\sim>>\bah^{\mu'}_{{}^*\car_{K-S},\cl_\bq}$$
which is the identity on the generators. (We use that 
$\mu(\hg_k)/\hg_k=\mu'({}^*\g_k)/{}^*\g_k$.) This induces a bijection
$$\Irr_{\bde,r_0}\bah^\mu_{\car_{K-S},\cl_\bq}\lra
\Irr_{\bde,r_0}\bah^{\mu'}_{{}^*\car_{K-S},\cl_\bq}$$
for any $\bde\in\uct_\spa/\cw_{K-S}$. Taking union over all $\dd,\bde$ and 
composing with (d) we obtain a bijection
$$\Irr_{v_0}H^{\la,\la^*}_{\car,\cl}\lra
\sqc_{S\sub K,S\ne\e}\sqc_{\dd\in C_S^\di}\sqc_{\bde\in\uct_\spa/\cw_{K-S}}
\Irr_{\bde,r_0}\bah^{\mu'}_{{}^*\car_{K-S},\cl_\bq}.\tag e$$
Let $S\sub K,S\ne\e$. If we apply the definitions of 5.2 to $G_{I-S},G_J$ 
(instead of $G,L$) then $R,W$ of 5.2 become ${}^*\car_{K-S},\cw_{K-S}$. We see
that $\bah^{\mu'}_{{}^*\car_{K-S},\cl_\bq}$ (as above) may be interpreted as
$\bah(G_{I-S},G_J,\cc,\cf)$. Moreover, if $\dd\in C_S^\di$ we have
$G_{I-S}=Z_G(p(\dd))$ (see 6.20). Hence we may rewrite (e) as
$$\Irr_{v_0}H^{\la,\la^*}_{\car,\cl}\lra
\sqc_{S\sub K,S\ne\e}\sqc_{\dd\in C_S^\di}\sqc_{\bde\in(\fz_J)_\spa/\cw_{K-S}}
\Irr_{\bde,r_0}\bah(Z_G(p(\dd)),G_J,\cc,\cf)\tag f$$
(We have used that $\uct=\fz_J$. Moreover $(\fz_J)_\spa$ defined in terms of
$Z^0_J$ coincides with $\uct_\spa$ defined in terms of $\ct$.) By 5.17(d), for
any $S\sub K,S\ne\e$ and $\dd\in C_S^\di$ we have a bijection
$$\sqc_{\bde\in(\fz_J)_\spa/\cw_{K-S}}\Irr_{\bde,r_0}
\bah(Z_G(p(\dd)),G_J,\cc,\cf)\lra\dfT^\spa(Z_G(p(\dd)),G_J,\cc,\cf).$$
(To define $\dfT^\spa(Z_G(p(\dd)),G_J,\cc,\cf)$ we use the parabolic subgroup
$P\cap G_{I-S}$ of $Z_G(p(\dd))=G_{I-S}$ with Levi subgroup $G_J$, see 7.6.) 
Taking union over all $\dd$ and composing with (f) we obtain a bijection
$$\Irr_{v_0}H^{\la,\la^*}_{\car,\cl}\lra\sqc_{S\sub K,S\ne\e}
\sqc_{\dd\in C_S^\di}\dfT^\spa(Z_G(p(\dd)),G_J,\cc,\cf).\tag g$$

\proclaim{Lemma 9.4} Let $f\in G\ta$. The following three conditions are 
equivalent:

(i) there exists $g\in G$ such that, for some $S\sub K,S\ne\e$ we have
$g\i f_\di g\in p(C_S^\di)$ and 
$g\i f_\spa g\in\dcc(Z_J^0)_\spa(U_P\cap G_{I-S})$;

(ii) there exists $g\in G$ such that for some $S\sub K,S\ne\e$ we have
$g\i fg\in\dcc p(C_S^\di)(Z_J^0)_\spa(U_P\cap G_{I-S})$;

(iii) $g\i fg\in X_{(J)}U_P$.
\endproclaim
It is clear that (i) and (ii) are equivalent and that (ii) implies (iii). Now
assume that (iii) holds. We show that (ii) holds. We may assume that 
$f\in X_{(J)}U_P$. By replacing $f$ by a $P$-conjugate, we may assume that
$f\in X_J\dcc U_P$. Let $f_s$ and $f_u$ be the semisimple and unipotent part of
$f$. Then $f_s$ is $P$-conjugate to an element of $X_J$. Hence, replacing $f$ 
by a $P$-conjugate, we may assume that $f_s\in X_J,f_u\in\dcc U_P$. Let 
$(\fz_J^1)_\di$ be the set of all $x'\in\fz_J^1$ such that
$x'=\sum_{k\in K}c_kb'_k$ with $c_k\in\k\i\di$. Then
$X_J=p((\fz_J^1)_\di)(Z_J^0)_\spa$. Now the $W^*$-action on $\fz_J^1$ restricts
to a $W^*$-action on $(\fz_J^1)_\di$ which has $\cup_{S\sub K;S\ne\e}C_S^\di$ 
as a fundamental domain. Hence by replacing $f$ by $nfn\i$ for some 
$n\in G^\ta$ such that 
$$nG_Jn\i=G_J,nX_Jn\i=X_J,n\dcc n\i=\dcc,n G_{(J)}n\i=G_{(J)},nPn\i=P'\tag a$$
($P'$ is another parabolic subgroup normalized by $\ta$ which has $G_{(J)}$ as
a Levi subgroup), we see that we may assume that $f_\di\in p(C_S^\di)$, 
$f_\spa\in(Z_J^0)_\spa\dcc U_{P'}$. Since $f_\spa\in Z_G(f_\di)$ and $f_\di$ as
well as $(Z_J^0)_\spa\dcc$ are contained in the Levi subgroup $G_{(J)}$ of 
$P'$, we have automatically $f_\spa\in(Z_J^0)_\spa\dcc(Z_G(f_\di)\cap U_{P'})$.
Since $Z_G(f_\di)=G_{I-S}$ (see 6.20), we have
$$f\in p(C_S^\di)(Z_J^0)_\spa\dcc(G_{I-S}\cap U_{P'}).$$
Now the parabolic subgroups $G_{I-S}\cap P'$ and $G_{I-S}\cap P$ of $G_{I-S}$
(both with Levi subgroup $G_J$, see 7.6) are conjugate under an element
$z\in G_{I-S}$ which normalizes $G_J$. Conjugation by $z$ carries
$G_{I-S}\cap U_{P'}$ to $G_{I-S}\cap U_P$, $p(C_S^\di)$ to $p(C_S^\di)$, 
$Z^0_J$ to $Z^0_J$ and $\dcc$ to $\dcc$ hence it carries $f$ to
$zfz\i\in p(C_S^\di)(Z_J^0)_\spa\dcc(U_P\cap G_{I-S})$. Thus, $f$ satisfies 
(ii). The lemma is proved.

\subhead 9.5\endsubhead
Let $f\in G\ta$. Assume that $f_\di\in p(C_S^\di)$ where $S\sub K,S\ne\e$. Then
$f_\spa\in Z_G(f_\di)=G_{I-S}$. We want to compare the varieties:
$$\align& A=\{g\in G;g\i fg\in G_{(J)}U_P\}/P,\\&
A'=\{g'\in G_{I-S};g'{}\i f_\spa g'\in Z_J^0\dcc(U_P\cap G_{I-S})\}
/(P\cap G_{I-S}).\endalign$$
Let $P'$ be any parabolic subgroup of $G$ such that $\ta P'\ta\i=P'$, $G_{(J)}$
is a Levi subgroup of $P'$ and $P'\cap G_{I-S}=P\cap G_{I-S}$. Let 
$A(P')=\{g\in G;g\i fg\in G_{(J)}U_{P'}\}/P'$. Define $f_{P'}:A'@>>>A(P')$ by
$g'(P\cap G_{I-S})\m g'P'$. This is clearly injective. We can find 
$n\in G_{I-S}$ such that 9.4(a) holds. Define $F_n:A(P')@>\sim>>A$ by 
$gP'\m gnP$. The composition $F_nf_{P'}:A'@>>>A(P')$ is injective. Its image
$A_{P'}$ depends only on $P'$, not on $n$. By the argument in the proof of 9.4
we see that $A$ is the disjoint union of finitely many subvarieties $A_{P'}$ 
(for the various $P'$ as above) and each $A_{P'}$ is isomorphic to $A'$. It 
follows easily that 
$$\Irr_1\bZ_G(f)=\Irr_0\bZ_{Z_G(f_\di)}(f_\spa).\tag a$$
(Note that $Z_G(f)=Z_{Z_G(f_\di)}(f_\spa)$.) Using (a) and 9.4 we see that we
have a bijection
$$\sqc_{S\sub K,S\ne\e}\sqc_{\dd\in C_S^\di}\dfT^\spa(Z_G(p(\dd)),G_J,\cc,\cf)
\lra\fR(G\ta,G_J,\cc,\cf)\tag b$$
given by $(\dd,(f,\rh))\m(p(\dd)f,\rh)$. Here $f\in(Z_G(p(\dd)))_\spa$. 
Composing (b) with 9.3(g), we obtain a bijection
$\Irr_{v_0}H^{\la,\la^*}_{\car,\cl}\lra\fR(G\ta,G_J,\cc,\cf)$. This is, by 
definition, the bijection 9.1(b).

\head 10. The main results\endhead
\subhead 10.1\endsubhead
In 10.2-10.7 we preserve the setup of 9.1.

\proclaim{Lemma 10.2}Assume that $\card(K)\ge 2$. There is a canonical 
injective map $\io:\ct/\cw@>>>(T\ta)/\cn$ whose image is exactly the image of 
$X_J$ in $(T\ta)/\cn$.
\endproclaim
By definition, $\io$ sends the $\cw$-orbit of $t\in\ct$ to the $\cn$-orbit of
$p(x')$ where $x'\in\fz_J^1$ is such that $p'(x')=t$. Assume that $x''$ is an 
element of $\fz_J^1$ such that $p'(x'')=w_1(t)$ where $w_1\in\cw$. Then there 
exists $w\in W^*$ such that $p'(w(x''))=t$. Since $p'(w(x'')-x')=1$, we have 
$w(x'')-x'\in\cl'$. Hence $x'=w'w(x'')$ for some $w'\in W^*$. In particular, 
$x'=\ti w(x'')$ for some $w\in W^a$. Hence $p(x')=np(x'')n\i$ for some 
$n\in\cn$. Thus, $\io$ is well defined.

We show that $\io$ is injective. Now $C'$ is a set of representatives for the 
$W^a$-orbits on $\ft^1$. Similarly, $D$ is a set of representatives for the 
$W^*$-orbits on $\ft^1_J$. Let $t_1,t_2\in\ct$ be such that the $\cw$-orbit of
$t_1$ and the $\cw$-orbit of $t_2$ have the same image under $\io$. Let 
$x_1,x_2\in\ft^1_J$ be such that $p'(x_1)=t_1,p'(x_2)=t_2$. We may assume that
$x_1\in D,x_2\in D$. Since $D\sub C'$, we have $x_1\in C',x_2\in C'$. By 
assumption we have $p(x_1)=np(x_2)n\i$ for some $n\in\cn$. Hence $x_1=w(x_2)$ 
for some $w\in W^a$. Since $x_1,x_2\in C'$, it follows that $x_1=x_2$. Hence 
$t_1=t_2$. This shows that $\io$ is injective. The fact that the image of $\io$
is exactly the image of $X_J$ in $(T\ta)/\cn$ is obvious. The lemma is proved.

\subhead 10.3\endsubhead
We show that the bijection 9.1(b) does not depend on the choice of $\di$ as in
9.1(a). When $v_0=1$ this is obvious: we have $\di=\bc$. Assume now that
$v_0\ne 1$. It is enough to show that one can define a map
$$\Irr_{v_0}H(G\ta,G_J,\cc,\cf)@>>>\fR(G\ta,G_J,\cc,\cf)$$
purely in terms of $\spa$ and which coincides with the map defined in \S9 in 
terms of any given $\di$. We can assume that $\card(K)\ge 2$.

Let $M\in\Irr_{v_0}H(G\ta,G_J,\cc,\cf)$. We want to attach to $M$ a pair 
$(f,\rh)$ (up to $G$-conjugacy) where $f\in G\ta$ and $\rh\in\Irr_1\bZ_G(f)$. 
We will only indicate the definition of the $G$-conjugacy class of $f$. (A 
similar definition applies to $\rh$.) 

By 1.3, we have $M\in\Irr_{\Si,v_0}H^{\la,\la^*}_{\car,\cl}$ for a well defined
$\cw$-orbit $\Si$ on $\ct$. Let $c$ be a fibre of $\Si@>>>\ct/\ct_\spa$
(restriction of $\ct@>>>\ct/\ct_\spa$). Define $\car_c,\cw^c,H_c$ in terms of 
$\ct,\car$ in the same way as $R_c,W^c_0,H_c$ were defined in [8.1]${}'$ and 
[8.3]${}'$ (see 2.2) in terms of $\ct,R$. By 3.2, to $M$ corresponds an object 
$M'\in\Irr_{c,v_0}H_c$. As in [9.2]${}'$, we can find an element $t_0\in\ct$ 
whose stabilizer in $\cw$ equals $\cw^c$ and a $\cw^c$-orbit $\bac$ in 
$\ft_\spa$ such that $t_0\ex_\ct(\bac)=c$. By 4.2 (for this $t_0$), $M'$ 
corresponds to an object $M''\in\Irr_{r_0}^\spa\bah$ where $\bah$ is attached 
to $H_c$ as in 4.1. Let $\ti t_0$ be an element of $X_J$ such that $\io(\cw t)$
is the $\cn$-orbit of $\ti t_0$ (see 10.2). Now $\bah$ may be interpreted as 
the algebra $\bah(G',G_J,\cc,\cf)$ where $G'=Z_G(\ti t_0)$. (Note that $G_J$ is
the Levi subgroup of some parabolic subgroup of $G'$.) Under 5.17(d), to $M''$
corresponds a pair $(f',\rh')$ where $f'\in G'_\spa$ is well defined up to 
conjugation in $G'$ and $\rh'\in\Irr_0\bZ_{G'}(f')$. We set 
$f=\ti t_0f'=f'\ti t_0$. Then the $\Ad(G)$-orbit of $f$ is well defined by $M$.
Note that we have not used $\di$ in this definition. Thus we have the following
result.

\proclaim{Theorem 10.4} Assume that $v\in\bc^*$ is either $1$ or is not a root
of $1$. Let $r_0\in\bc$ be as in 9.1. There is a bijection 
$$\Irr_{v_0}H(G\ta,G_J,\cc,\cf)\lra\fR(G\ta,G_J,\cc,\cf)\tag a$$
depending only on $r_0$, which for any $\di$ as in 9.1(a) coincides with the
bijection 9.1(b).
\endproclaim
It is likely that the bijection (a) is independent of the choice of $r_0$. 
(Some evidence is given in 10.7.)

\proclaim{Theorem 10.5} We preserve the setup of 10.4. 
Assume that $v_0$ is not a root of $1$. Let $\z:\bc^*@>>>\br$ be
a group homomorphism such that $\z(v_0)\ne 0$. Assume that under 10.4(a), 
$M\in\Irr_{v_0}H(G\ta,G_J,\cc,\cf)$ corresponds to 
$(f,\rh)\in\fR(G\ta,G_J,\cc,\cf)$. Then

(a) $M$ is $\z$-tempered if and only if $f\in\hG^{\Ker\z}$;

(b) $M$ is $\z$-square integrable if and only if any torus in $Z_G(f)$ is
$\{1\}$.
\endproclaim
Assume first that $\card(K)=1$. Then the unique 
$M\in\Irr_{v_0}H(G\ta,G_J,\cc,\cf)$ is obviously $\z$-square integrable. The 
unique element of $\fR(G\ta,G_J,\cc,\cf)$ may be represented in the form 
$(f,\rh)$ where $f=su$ with $s$ being the unique element of $X_J$ and 
$u\in\dcc$. Assume that $T_1$ is a torus in $Z_G(f)=Z_{G_J}(u)$. In our case, 
$G_J$ is a semisimple group and $u$ is a distinguished unipotent element of 
$G_J$. Hence any torus in $Z_{G_J}(u)$ is trivial. Thus, $T=1$ hence (b) holds
in this case. We show that $f\in\hG^{\Ker\z}$. It is enough to show that 
$s\in\hG^{\Ker\z}$. Since $Z_G(s)$ is semisimple, $s$ is of finite order. 
Since $\Ker\z$ contains all roots of $1$, we have $s\in\hG^{\Ker\z}$. Thus, 
(a) holds in this case. 

In the remainder of the proof we assume that $\card(K)\ge 2$. Let $r_0,\spa$ be
as in 9.1.

{\it Case 1.} We assume that $\z(\bc^*)\sub\bq$. In this case, $\z$ restricts
to an isomorphism $\ex(\spa)@>\sim>>\bq$. In particular we have
$\bc^*=\Ker\z\opl\ex(\spa)$. Let $\di=\{a\in\bc;\z(\ex(a))=1\}$. Then $\di$ 
satisfies 9.1(a) and $\Ker\z=\ex(\di)$. We can use the definition of the 
bijection 9.1(b) in terms of this $\di$.

We prove (a). Using Lemmas 3.4, 4.3, 5.18 and the definitions we are reduced to
verifying that for $f\in G\ta$ the following two conditions are equivalent:

(i) $f_\spa\in Z_G(f_\di)^{\Ker\z}$;

(ii) $f\in\hG^{\Ker\z}$.
\nl
(When applying 4.3, we can choose $t_0$ in 4.1 so that $t_0\in\ct_\di$; then it
is unique. The assumption $t_0\in\ct^{\Ker\z}$ of 4.3 is automatically verified
since $\ct^{\Ker\z}=\ct_\di$.)

Now (i) is equivalent to the condition $f_\spa\in\hG^{\Ker\z}$. This is 
equivalent to (ii) since we have automatically $f_\di\in\hG^{\Ker\z}$ (since
$\Ker\z=\ex(\di)$).

We prove (b). Using Lemmas 3.5, 4.4, 5.16 and the definitions we are reduced to
verifying that for $f\in G$ the following two conditions are equivalent:

(iii) $Z_G(f_\di)$ is semisimple and $f_\spa$ is a distinguished unipotent 
element of $Z_G(f_\di)$;

(iv) any torus in $Z_G(f)$ is $\{1\}$.
\nl
This is immediate. 

{\it General case.} As in 1.3, let $\cl^+$ be the set of all $x\in\cl$ such 
that $\lan x,\hah_k\r\ge 0$ for all $\hg_k\in\Pi$. We can find $x_1,\do,x_N$ in
$\cl^+-\{0\}$ such that $\cl^+=\sum_{k=1}^N\bn x_k$. There exists a finite
subset $\ct_0$ of $\ct$ such that, for $t\in\ct$, the weight space $M_t$ is 
zero unless $t\in\ct_0$. Let $A$ be the (finite) subset of $\bc^*$ consisting 
of all numbers of the form $x_k(t)$ with $k\in [1,N],t\in\ct_0$. Let $B$ be the
(finite) subset of $\bc^*$ consisting of the eigenvalues of $f$ in a fixed 
faithful $V\in\ci_{\hG}$. Then (a),(b) can be restated as (c),(d) below:

(c) we have $\z(a)/\z(v_0)\ge 0$ for all $a\in A$ if and only if $\z(b)=0$ for
all $b\in B$;

(d) we have $\z(a)/\z(v_0)>0$ for all $a\in A$ if and only (iv) holds.

Assume first that $\z(a)/\z(v_0)\ge 0$ for all $a\in A$ and $\z(b)\ne 0$ for
some $b\in B$. We can find a $\bq$-linear form $u:\br@>>>\bq$ such that
$u(\z(v_0))\ne 0,u(\z(b))\ne 0$ and $u(\z(a))/u(\z(v_0))\ge 0$ for all 
$a\in A$. Applying Case 1 to $u\z:\bc^*@>>>\bq$ instead of $\z$, we see that
$u(\z(b))=0$, a contradiction.

Assume next that $\z(a)/\z(v_0)<0$ for some $a\in A$ and $\z(b)=0$ for all
$b\in B$. We can find a $\bq$-linear form $u:\br@>>>\bq$ such that
$u(\z(v_0))\ne 0,u(\z(a))/u(\z(v_0))<0$. We have $u(\z(b))=0$ for all $b\in B$.
Applying Case 1 to $u\z:\bc^*@>>>\bq$ instead of $\z$, we see that
$u(\z(a))/u(\z(v_0))\ge 0$, a contradiction. Thus, (c) holds.

Assume now that $\z(a)/\z(v_0)>0$ for all $a\in A$. We can find a $\bq$-linear
form $u:\br@>>>\bq$ such that $u(\z(v_0))\ne 0$ and $u(\z(a))/u(\z(v_0))>0$ for
all $a\in A$. Applying Case 1 to $u\z:\bc^*@>>>\bq$ instead of $\z$, we see 
that (iv) holds.

Conversely, assume that $\z(a)/\z(v_0)\le 0$ for some $a\in A$. We can find a
$\bq$-linear form $u:\br@>>>\bq$ such that $u(\z(v_0))\ne 0$ and 
$u(\z(a))/u(\z(v_0))\le 0$. Applying Case 1 to $u\z:\bc^*@>>>\bq$ instead of
$\z$, we see that (iv) does not hold. Thus, (d) holds. The theorem is proved.

\proclaim{Corollary 10.6} Let $v_0,\z,M$ be as in 10.5. 

(a) $M$ is $\z$-tempered if and only if the following holds: for any $t\in\ct$
such that $M_t\ne 0$ and any $x\in\cl^+$ we have
$\z(x(t))/\z(v_0)\in\bq_{\ge 0}$.

(b) $M$ is $\z$-square integrable if and only if the following holds: for any
$t\in\ct$ such that $M_t\ne 0$ and any $x\in\cl^+-\{0\}$ we have $x(t)=a v_0^n$
for some $n\in\bz_{>0}$ and some $a\in\bc^*$, a root of $1$.
\endproclaim
Assume that $M$ is $\z$-tempered and there exists $t\in\ct$ and $x\in\cl^+$ 
such that $M_t\ne 0$ and $\z(x(t))/\z(v_0)\n\bq_{\ge 0}$. Note that 
$\z(x(t))/\z(v_0)\in\br_{>0}$. Since $\z(x(t)),\z(v_0)$ are non-zero real 
numbers of the same sign and one is not a rational multiple of the other, we 
can find a $\bq$-linear form $u:\br@>>>\bq$ such that $u(\z(v_0))\ne 0$ and
$u(\z(x(t)))/u(\z(v_0))\in\bq_{<0}$. Hence $M$ is not $u\z$-tempered. Let $f$
correspond to $M$ as in 10.5. By 10.5(a) we have $f\in\hG^{\Ker\z}$. It follows
that $f\in\hG^{\Ker(u\z)}$. Using again 10.5(a) (for $u\z$ instead of $\z$) we
see that $M$ is $u\z$-tempered, a contradiction. This proves (a).

We prove (b). By the arguments in the proof of 3.4, 4.4, we are reduced to the
analogous statement for the algebras considered in \S5, which is proved in 
\cite{\TE, 1.22}. The corollary is proved.

\subhead 10.7\endsubhead
In the setup of 9.1 (with $\card(K)\ge 2$) we consider an element $f\in G\ta$ 
such that $(f,\rh)\in\fR(G\ta,G_J,\cc,\cf)$ for some $\rh$. We can write 
uniquely $f=f_\di f_\spa'f_u$ (three commuting factors) where $f_\di\in\hG_\di$
is semisimple, $f'_\spa\in G_\spa$ is semisimple, $f_u\in G$ is unipotent. 
Replacing $(f,\rh)$ by a $G$-conjugate, we may assume that $f_\di\in X_J$ (see
9.5) so that $G_J\sub Z_G(f_\di)$. Let $t\in\ct_\di$ be such that $\io(\cw t)$
is the $\cn$-orbit of $f_\di$. (See 10.2.) Let $\cw_t$ be the stabilizer of $t$
in $\cw$. Let $\ph\in\Hom(SL_2(\bc),Z_G(f_\di f_\spa'))$ be such that 
$\ph\left(\sm 1&1\\0&1\esm\right)=f_u$. Let
$$\ti f=f_\spa'\ph\left(\sm v_0&0\\0&v_0\i\esm\right)\in Z_G(f_\di).$$ 
Let $\ph^0\in\Hom(SL_2(\bc),G_J)$ be such that
$\ph^0\left(\sm 1&1\\0&1\esm\right)\in\dcc$. Using 5.6, we see that there 
exists $z\in Z_G(f_\di)$ such that 
$$z\ti fz\i\ph^0\left(\sm v_0\i&0\\0&v_0\esm\right)\in(Z_J^0)_\spa\tag a$$
and that the orbit of the element (a) under the normalizer of $Z_J^0$ in
$Z_G(f_\di)$ does not depend on the choice of $z$. Since 
$(Z_J^0)_\spa=\ct_\spa$ (both may be identified with $(\fz_J)_\spa$ using
$\ex_{Z_J^0},\ex_\ct$) we may regard this orbit as a $\cw_t$-orbit $c$ in 
$\ct_\spa$. Let $\Si$ be the $\cw$-orbit in $\ct$ that contains $tc$. Then 
$\Si$ depends only on the $G$-conjugacy class of $(f,\rh)$ and $(f,\rh)\m\Si$
is a map
$$\fR(G\ta,G_J,\cc,\cf)@>>>\ct/\cw.\tag b$$

We now give a second definition of the map (b). Consider $(f,\rh)$ as above. We
can write uniquely $f=f_sf_u$ (two commuting factors) where $f_s\in\hG$ is 
semisimple, $f_u\in G$ is unipotent. Replacing $(f,\rh)$ by a $G$-conjugate, we
may assume that $f_s\in X_J$ (see 9.5) so that $G_J\sub Z_G(f_s)$. Let 
$\ph\in\Hom(SL_2(\bc),Z_G(f_s))$ be such that 
$\ph\left(\sm 1&1\\0&1\esm\right)=f_u$. Let
$$\hat f=f_s\ph\left(\sm v_0&0\\0&v_0\i\esm\right)\in Z_G(f_s).$$ 
Let $\ph^0\in\Hom(SL_2(\bc),G_J)$ be such that
$\ph^0\left(\sm 1&1\\0&1\esm\right)\in\dcc$. Using 5.6, we see that there 
exists $z\in Z_G(f_s)$ such that 
$$a=:z\hat fz\i\ph^0\left(\sm v_0\i&0\\0&v_0\esm\right)\in X_J.$$
Let $\Si$ be the $\cw$-orbit in $\ct$ such that $\io(\Si)$ is the $\cn$-orbit
of $a$ (see 10.2). Then $\Si$ depends only on the $G$-conjugacy class of 
$(f,\rh)$ and $(f,\rh)\m\Si$ coincides with the map (b). In the setup of 10.4,
each $M\in\Irr_{v_0}H^{\la,\la^*}_{\car,\cl}$ belongs to
$\Irr_{\Si,v_0}H^{\la,\la^*}_{\car,\cl}$ for a unique $\cw$-orbit $\Si$ in 
$\ct$ (as in 1.3(a)). Then $M\m\Si$ is a well defined map
$$\Irr_{v_0}H^{\la,\la^*}_{\car,\cl}@>>>\ct/\cw.\tag c$$ 
Composing this with the bijection 10.4(a) we obtain a map 
$\fR(G\ta,G_J,\cc,\cf)@>>>\ct/\cw$. This coincides with the map (b) (in its 
first form) as one sees using the definitions.

The fact that the map (b) (in its second form) is independent of the choice of
$r_0,\di$ and that the same is obviously true for the map (c), suggests that 
the bijection 10.4(a) is also independent of $r_0$.

\subhead 10.8\endsubhead
We fix $v_0\in\bc^*$ which is either $1$ or is not a root of $1$. Let $\fJ$ be
the set of all triples $(J,\cc,\cf)$ as in 7.4. (Here $\cf$ is given up to 
isomorphism.) Putting together the bijections 10.4(a) for various 
$(J,\cc,\cf)\in\fJ$, we obtain a bijection 
$$\sqc_{(J,\cc,\cf)\in\fJ}\Irr_{v_0}H(G\ta,G_J,\cc,\cf)\lra
\sqc_{(J,\cc,\cf)\in\fJ}\fR(G\ta,G_J,\cc,\cf).\tag a$$
Let $\fR(G\ta)$ be the set of all $(f,\rh)$ (modulo the $\Ad$-action of $G$) 
where $f\in G\ta$ and $\rh\in\Irr\bZ_G(f)$ (the set of isomorphism classes of 
irreducible representations of $\bZ_G(f)$). We will show below that 
$$\sqc_{(J,\cc,\cf)\in\fJ}\fR(G\ta,G_J,\cc,\cf)=\fR(G\ta).\tag b$$
Combined with (a), this gives a bijection
$$\sqc_{(J,\cc,\cf)\in\fJ}\Irr_{v_0}H(G\ta,G_J,\cc,\cf)\lra\fR(G\ta).\tag c$$
We prove (b). We fix $f\in G\ta$. We may assume that $f_\di\in p(C_S^\di)$ 
where $S\sub I,S\ne\e$. Then $Z_G(f_\di)=G_{I-S}$. Let $\rh\in\Irr\bZ_G(f)$. We
must show that there is a unique $(J,\cc,\cf)\in\fJ$ such that $\rh$ belongs to
$\Irr_1\bZ_G(f)$ (defined in terms of $G,G_J,\cc,\cf$) or equivalently (see 
9.5(a)) to $\Irr_0\bZ_{G_{I-S}}(f_\spa)$ (defined in terms of 
$G_{I-S},G_J,\cc,\cf$; we have necessarily $J\sub I-S$). Recall that 
$\bZ_G(f)=\bZ_{G_{I-S}}(f_\spa)$. Define $\xi\in(\fg_{I-S})_\spa$ by 
$\exp(\xi)=f_\spa$. We are reduced to verifying the following statement:

For any $\rh\in\Irr\bZ_{G_{I-S}}(\xi)$ there exists a unique
$(J,\cc,\cf)\in\fJ$ such that $J\sub I-S$ and $\rh$ is in 
$\Irr_0\bZ_{G_{I-S}}(\xi)$ (defined in terms of $G_{I-S},G_J,\cc,\cf$).

This follows from \cite{\II, \S8}.

\subhead 10.9\endsubhead
Let $\z:\bc^*@>>>\br$ be a homomorphism such that $\z(v_0)\ne 0$. By 10.5, the
bijection 10.8(c) restricts to a 

{\it bijection between the set of $\z$-tempered representations in the left 
hand side of 10.8(c) and $\{(f,\rh)\in\fR(G\ta);f\in\hG^{\Ker\z}\}$}

and to a

{\it bijection between the set of $\z$-square integrable representations in the
left hand side of 10.8(c) and the set of all $(f,\rh)\in\fR(G\ta)$ such that
any torus in $Z_G(f)$ is $\{1\}$.} 

Special cases of this result can be found in \cite{\KL},\cite{\RE},\cite{\WA}.

\subhead 10.10\endsubhead
Let $\bk,q$ be as in 1.1. Let $\ti{\bk}$ be a maximal unramified extension of
$\bk$. Let $\bg$ be a connected, adjoint simple algebraic group defined over 
$\bk$ which is split over $\ti{\bk}$. We identify $\bg$ with $\bg(\ti{\bk})$.
Assume that $\bg$ is of type dual in the sense of Langlands to $G$. Define 
$\bi$ as in \cite{\IM, 1.10}. This is the set of vertices of the affine Dynkin
graph of $\bg$. Let $\ti S(\bi)$ be the set of bijections $\bi@>\sim>>\bi$ that
preserve the graph structure. We have a canonical (surjective) homomorphism
from $\ti S(\bi)$ to the group of automorphisms of $G$ modulo the group of
inner automorphisms of $G$ (see \cite{\IM, 8.1}). Let $\ti S(\bi)_\ta$ be the 
fibre of this map over the coset of $\Ad(\ta):G@>>>G$. For $u\in\ti S(\bi)_\ta$
we can find a $\bk$-rational structure on $\bg$ (compatible with the 
$\ti{\bk}$-rational structure) with Frobenius map $F_u$ (see \cite{\IM, 1.1})
such that the permutation of $\bi$ induced by $F_u$ (as in \cite{\IM, 1.12}) is
equal to $u$. Then $\bg(\bk)$ coincides with the fixed point set $\bg^{F_u}$.
Let $\cu(\bg^{F_u})$ be the set of isomorphism classes of unipotent 
representations of $\bg^{F_u}$ (see \cite{\IM, 1.21}). 

\proclaim{Theorem 10.11} There is a natural bijection 
$\sqc_{u\in\ti S(\bi)_\ta}\cu(\bg^{F_u})\lra\fR(G\ta)$.
\endproclaim
By \cite{\IM, 1.22} we have a natural bijection between
$\sqc_{u\in\ti S(\bi)_\ta}\cu(\bg^{F_u})$ and the disjoint union of the sets of
irreducible representations (up to isomorphism) of a finite collection of 
affine Hecke algebras $\ch'(\bi,\bj,u,\be)$ given by a presentation of 
Iwahori-Matsumoto type with explicitly known parameters and with the
indeterminate $v$ being specialized to $\sqrt{q}$. (Here $u\in\ti S(\bi)_\ta$,
$\bj$ is a proper $u$-stable subset of $\bi$, $E$ is a unipotent cuspidal 
representation of the $F_u$-fixed points of the parahoric subgroup attached to
$\bj$). The various $\bj,u,\be$ are listed in the tables in \S11 as 
"arithmetic diagrams"; the corresponding affine Hecke algebras 
$\ch'(\bi,\bj,u,\be)$ are listed in the same tables as "H.A.". Rather 
surprisingly, it turns out out that these affine Hecke algebras are exactly the
same as the geometric affine Hecke algebras attached to $G\ta$ (which are also
described in the tables of \S11). Therefore, the theorem follows from 10.8(c)
with $v_0=\sqrt{q}$.

\subhead 10.12\endsubhead
For any homomorphism $\chi:Z_G@>>>\bc^*$, let $\fR(G\ta)_\chi$ be the subset of
$\fR(G\ta)$ consisting of all $(f,\rh)$ such that via the obvious map 
$Z_G@>>>\bZ_G(f)$, $Z_G$ acts on $\rh$ through the character $\chi$. This gives
us a partition $\fR(G\ta)=\sqc_\chi\fR(G\ta)_\chi$. On the other hand, the
bijection in 10.11 induces a partition of $\fR(G\ta)$ into subsets indexed by 
the elements of $\ti S(\bi)_\ta$. This coincides with the previous partition of
$\fR(G\ta)$.

\subhead 10.13\endsubhead
Let $\bk,q,\ti{\bk},\bg$ be as in 10.10 except that $\bg$ is no longer assumed
to be split over $\ti{\bk}$. We identify $\bg$ with $\bg(\ti{\bk})$. One can 
still define the set $\bi$ which indexes the maximal parahoric subgroups of 
$\bg$ (see \cite{\TI}), a Dynkin graph with set of vertices $\bi$ and a
bijection $u:\bi@>>>\bi$ (preserving the graph structure, including the
orientation of the double or triple edges) which specifies the $\bk$-rational
structure on $\bg$. Then $\cu(\bg(\bk))$ can be defined in the same way as in
\cite{\IM, 1.21}. We can find another connected adjoint simple algebraic group
$\bg'$ which is split over $\ti{\bk}$ whose associated $\bi,u:\bi@>>>\bi$ (as 
in 10.10) is the same as the $\bi,u$ associated to $\bg$ and such that the
corresponding Dynkin graph (for $\bg'$) is the same as that for $\bg$, except
possibly for the orientation of the double or triple edges. Then 
$\cu(\bg'(\bk))$ is defined. Moreover we have a natural bijection
$$\cu(\bg(\bk))\lra\cu(\bg'(\bk)).\tag a$$
Indeed, each side of (a) is naturally in bijection with the disjoint union of
the sets of irreducible representations of a finite collection of affine Hecke
algebras given by a presentation of Iwahori-Matsumoto type at $v=\sqrt{q}$. But
the affine Hecke algebras associated to the two sides of (a) are the same, 
since the recipe that describes them is not sensitive to the orientation of the
double or triple edges. (This is analogous to the known statement that the sets
of unipotent representations of the finite groups $SO_{2n+1}(F_q)$ and 
$Sp_{2n}(F_q)$ are in bijection.) Since 10.11, 10.12 are applicable to 
$\cu(\bg'(\bk))$, they also provide, via (a), a parametrization of 
$\cu(\bg(\bk))$.

We describe the various pairs $(\bg,\bg')$ using the names in the tables of
\cite{\TI}.

$(CB_n,C_n),(BC_n,B_n),(CBC_n,C_n),(G_2^I,G_2),(F_4^I,F_4),({}^2BC_n,{}^2B_n),
({}^2CB_n,{}^2C_n)$.

\subhead 10.14\endsubhead
It is likely that our results can be extended to the case where the assumption
that $\bg$ is adjoint simple is weakened to the assumption that $\bg$ is 
semisimple. Indeed, our main technique, that of reducing to the case of graded
Hecke algebras is still available in this more general case (see \cite{\LG}).

\head 11. Tables\endhead
\subhead 11.1\endsubhead
In this section we list the various possibilities for $G\ta$ and $J$ as in 7.4
assuming that $d\ge 2$. (The cases where $d=1$ are listed in \cite{\IM}.) In 
each case we describe the affine Dynkin graph associated to the affine Cartan 
matrix $(a_{i_1,i_2})$ (resp. $({}'a_{i_1,i_2})$) in 6.21; we call this the 
$(\g_i)$-graph (resp. the $(\b_i)$-graph). Both these graphs have vertices in
bijection with $I$. The vertices of the $(\b_i)$-graph that are inside a box 
correspond to the subset $J$ of $I$. The full subgraph with vertices $J$ is the
Dynkin graph of $G_J$ (see 6.23). We also describe the affine Dynkin graph
associated to the affine Cartan matrix $(\hg_k(\hah_{k'}))_{k,k'\in K}$. We 
call this the $\fla-\sh$ diagram; its vertices are in bijection with $K$ and we
attach to any vertex the symbol $\fla$ or $\sh$ according to whether the
corresponding element of $K$ is in $K^\fla$ or $K^\sha$. For any vertex 
correponding to $k\in K$ we specify some data of the form $a\T b\T c$ where
$a=\uc_k,b=\baz_k,c=d_k$. From the $\fla-\sh$ diagram one obtains an affine 
Hecke algebra as in 8.2, 8.3. This affine Hecke algebra (in a presentation of 
Iwahori-Matsumoto type) is denoted by H.A. It turns out to be the same as the 
affine Hecke algebra attached to the arithmetic diagram (see 10.11) which is 
also given in each case. The notation for affine Hecke algebras follows the
conventions of \cite{\IM, 6.9, 6.11}.

\subhead 11.2\endsubhead
$G$ is of type $A_n$, $n$ even, $d=2$.

$a\in 4\bz,b\in 1+4\bz$, $n+1=2s+2+a(a+1)/2+b(b+1)/2$, $s\ge 1$.

$(\g_i)$-graph:    
$$\s\Ra\s\h\s\h\s\h\s\h\s\h\do\h\s\h\s\h\s\Lar\s$$
$(\b_i)$-graph: 
$$\us{C_q}\to{\x{\sm\s&\Ra&\s&\s&\do&\s\esm}}\h\s_1\h\s_2\h\s_3\h\do\h\s_s\h
\us{B_p}\to{\x{\sm\s&\do&\s&\s&\Ra&\s\esm}}$$
if $a-b\ne-1,|a+b+1|\ne 2$; here $2p+1=(a+b+1)^2/4,2q=(a-b-1)(a-b+1)/4$; 
$$\s_1\Ra\s_2\h\s_3\h\do\h\s_s\h\us{B_p}\to{\x{\sm\s&\do&\s&\s&\Ra&\s\esm}}$$
if $a-b=-1,|a+b+1|\ne 2$; here  $2p+1=(a+b+1)^2/4$;
$$\us{C_q}\to{\x{\sm\s&\Ra&\s&\s&\do&\s\esm}}\h\s_1\h\s_2\h\s_3\h\do\h\s_{s-1}
\Lar\s_s$$
if $a-b\ne-1,|a+b+1|=2$; here $2q=(a-b-1)(a-b+1)/4$;
$$\s_1\Ra\s_2\h\s_3\h\do\h\s_{s-1}\Ra\s_s$$
if $a-b=-1,|a+b+1|=2$.

For $s\ge 2$, the $\fla-\sh$-diagram is
$$\fla_s^{|a+b+1|\T 1\T 2}\Lar\sh_{s-1}^{2\T 1\T 2}\h\sh_{s-2}^{2\T 1\T 2}\h
\do\h\sh_2^{2\T 1\T 2}\Ra\fla_1^{|a-b|\T 1\T 2}$$
if $a-b\ne-1,|a+b+1|\ne 2$; 
$$\fla_s^{|a+b+1|\T 1\T 2}\Lar\sh_{s-1}^{2\T 1\T 2}\h\sh_{s-2}^{2\T 1\T 2}\h
\do\h\sh_2^{2\T 1\T 2}\Ra\fla_1^{2\T\fra{1}{2}\T 2}$$
if $a-b=-1,|a+b+1|\ne 2$;
$$\fla_s^{2\T\fra{1}{2}\T 4}\Lar\sh_{s-1}^{2\T 1\T 2}\h\sh_{s-2}^{2\T 1\T 2}\h
\do\h\sh_2^{2\T 1\T 2}\Ra\fla_1^{|a-b|\T 1\T 2}$$
if $a-b\ne-1,|a+b+1|=2$;
$$\fla_s^{2\T\fra{1}{2}\T 4}\Lar\sh_{s-1}^{2\T 1\T 2}\h\sh_{s-2}^{2\T 1\T 2}\h
\do\h\sh_2^{2\T 1\T 2}\Ra\fla_1^{2\T\fra{1}{2}\T 2}$$
if $a-b=-1,|a+b+1|=2$.

For $s=1$, the $\fla-\sh$-diagram is $\e$.

H.A: $\ti C_{s-1}^{sc}[{}_{|2a+1|}2_{|2b+1|}]$ if $s\ge 2$ and $\e$ if $s=1$.

Arithmetic diagram: $\ti A_n$, $u^2=1,u\ne 1$, $\bj$ of type 
$A_{p'-1}\T A_{q'-1}$ (both components are $u$-stable), 
$p'=a(a+1)/2,q'=b(b+1)/2$.

\subhead 11.3\endsubhead
$G$ is of type $A_n$, $n$ odd, $d=2$.

Either $a\in 4\bz,b\in 3+4\bz$ or $a\in 2+4\bz,b\in 1+4\bz$.

$n+1=2s+2+a(a+1)/2+b(b+1)/2$, $s\ge 1$.

$(\g_i)$-graph:
$$\matrix{}&{}&\s&{}&{}&{}&{}&{}&{}&{}&{}\\{}&{}&\vert&{}&{}&{}&{}&{}&{}&{}&{}
\\ \s&\hor&\s&\hor&\s&\hor&\do&\hor&\s&\Ra&\s\endmatrix$$
$(\b_i)$-graph:
$$\us{D_p}\to{\x{\sm\s&{}&{}&{}\\ \s&\s&\do&\s\\ \s&{}&{}&{}\esm}}\h\s_1\h\s_2
\h\s_3\h\do\h\s_s\h\us{C_q}\to{\x{\sm\s&\do&\s&\s&\Lar&\s\esm}}$$
if $a-b\ne 1$, $a+b\ne-1$; here $2p=(a+b+1)^2/4,2q=(a-b-1)(a-b+1)/4$;
$$\us{D_p}\to{\x{\sm\s&{}&{}&{}\\ \s&\s&\do&\s\\ \s&{}&{}&{}\esm}}\h\s_1\h\s_2
\h\s_3\h\do\h\s_{s-1}\Lar\s_s$$
if $a-b=1$, $a+b\ne-1$; here $2p=(a+b+1)^2/4$;
$$\matrix{}&{}&\s_1&{}&{}&{}&{}&{}&{}&{}&{}\\{}&{}&\vert&{}&{}&{}&{}&{}&{}&{}&
{}\\\s_2&\hor&\s_3&\hor&\s_4&\hor&\do&\hor&\s_s&\hor&
\us{C_q}\to{\x{\sm\s&\do&\s&\s&\Lar&\s\esm}}\endmatrix$$
if $a-b\ne 1$, $a+b=-1$; here $2q=(a-b-1)(a-b+1)/4$; 
$$\matrix{}&{}&\s_1&{}&{}&{}&{}&{}&{}&{}&{}\\{}&{}&\vert&{}&{}&{}&{}&{}&{}&{}&
{}\\\s_2&\hor&\s_3&\hor&\s_4&\hor&\do&\hor&\s_{s-1}&\Lar&\s_s\endmatrix$$
if $a-b=1$, $a+b=-1$.

For $s\ge 2$, the $\fla-\sh$-diagram is
$$\fla_s^{|a+b+1|\T 1\T 2}\Lar\sh_{s-1}^{2\T 1\T 2}\h\sh_{s-2}^{2\T 1\T 2}\h
\do\h\sh_2^{2\T 1\T 2}\Ra\fla_1^{|a-b|\T 1\T 2}$$
if $a-b\ne 1$, $a+b\ne-1$;
$$\fla_1^{|a+b+1|\T 1\T 2}\Lar\sh_2^{2\T 1\T 2}\h\sh_3^{2\T 1\T 2}\h\do\h
\sh_{s-1}^{2\T 1\T 2}\Ra\fla_s^{2\T 1\T 1}$$
if $a-b=1$, $a+b\ne-1$;
$$\matrix\sh^{|a-b|\T 1\T 2}_s&\Lar&\sh^{2\T 1\T 2}_{s-1}&\hor&
\sh^{2\T 1\T 2}_{s-2}&\hor&\do&\hor&\sh^{2\T 1\T 2}_3&\hor&\sh^{2\T 1\T 2}_2\\
{}&               {}&          {}&   {}&       {}&    {}&  {}&  {}&\vert&{}&\\
{}&               {}&          {}&   {}&       {}&    {}&  {}&  {}&
\sh^{2\T 1\T 2}_1&{}&\endmatrix$$
if $a-b\ne 1$, $a+b=-1$;
$$\matrix\sh^{2\T 1\T 1}_s&\Lar&\sh^{2\T 1\T 2}_{s-1}&\hor&
\sh^{2\T 1\T 2}_{s-2}&\hor&\do&\hor&\sh^{2\T 1\T 2}_3&\hor&\sh^{2\T 1\T 2}_2\\
{}&               {}&          {}&   {}&       {}&    {}&  {}&  {}&\vert&{}&\\
{}&               {}&          {}&   {}&       {}&    {}&  {}&  {}&
\sh^{2\T 1\T 2}_1&{}&\endmatrix$$
if $a-b=1$, $a+b=-1$.

For $s=1$, the $\fla-\sh$-diagram is $\e$.

H.A.: $\ti C_{s-1}^{sc}[{}_{|2a+1|}2_{|2b+1|}]$ if $s\ge 2$ and $a+b\ne-1$;
$\ti C_{s-1}[{}_{|2a+1|}2_{|2b+1|}]$ if $s\ge 2$ and $a+b=-1$; $\e$ if $s=1$.

Arithmetic diagram: $\ti A_n$, $u^2=1,u\ne 1$, $\bj$ of type 
$A_{p'-1}\T A_{q'-1}$ (both components are $u$-stable), 
$p'=a(a+1)/2,q'=b(b+1)/2$.

\subhead 11.4\endsubhead
$G$ is of type $D_n$, $d=2$.

$a\ge 1$ odd; $b\ge 0$ even; $n+1=s+a^2+b^2$, $s\ge 1$. 

$(\g_i)$-graph:    
$$\s\Ra\s\h\s\h\s\h\s\h\s\h\do\h\s\h\s\h\s\Lar\s$$
$(\b_i)$-graph: 
$$\us{B_p}\to{\x{\sm\s&\Lar&\s&\s&\do&\s\esm}}\h\s_1\h\s_2\h\s_3\h\do\h\s_s\h
\us{B_q}\to{\x{\sm\s&\do&\s&\s&\Ra&\s\esm}}$$
if $a+b\ne 1$,$|a-b|\ne 1$; here $2p+1=(a+b)^2,2q+1=(a-b)^2$;
$$\us{B_p}\to{\x{\sm\s&\Lar&\s&\s&\do&\s\esm}}\h\s_1\h\s_2\h\s_3\h\do\h\s_{s-1}
\Ra\s_s$$
if $a+b\ne 1$,$|a-b|=1$; here $2p+1=(a+b)^2$;
$$\s_1\Lar\s_2\h\s_3\h\do\h\s_{s-1}\Ra\s_s$$
if $a+b=1$,$|a-b|=1$.

For $s\ge 2$, the $\fla-\sh$-diagram is
$$\fla_1^{2(a+b)\T 1\T 1}\Lar\sh_2^{2\T 1\T 1}\h
\sh_3^{2\T 1\T 1}\h\do\h\sh_{s-1}^{2\T 1\T 1}\Ra\fla_s^{2|a-b|\T 1\T 1}$$
if $a+b\ne 1$,$|a-b|\ne 1$; 
$$\fla_1^{2(a+b)\T 1\T 1}\Lar\sh_2^{2\T 1\T 1}\h\sh_3^{2\T 1\T 1}\h\do\h
\sh_{s-1}^{2\T 1\T 1}\Ra\fla_s^{2\T\fra{1}{2}\T 2}$$
if $a+b\ne 1$,$|a-b|=1$; 
$$\fla_s^{2\T\fra{1}{2}\T 2}\Lar\sh_{s-1}^{2\T 1\T 1}\h\sh_{s-2}^{2\T 1\T 1}\h
\do\h\sh_2^{2\T 1\T 1}\Ra\fla_1^{2\T\fra{1}{2}\T 2}$$
if $a+b=1$,$|a-b|=1$.

For $s=1$, the $\fla-\sh$-diagram is $\e$.

H.A.: $\ti C_{s-1}^{sc}[{}_{2a}1_{2b}]$ if $s\ge 2$ and $\e$ if $s=1$.

Arithmetic diagram: $\ti D_n$, $u:\bi@>>>\bi$ has exactly $n-1$ fixed points,
$\bj$ of type $D_{p'}\T D_{q'}$ ($u$ acts non-trivially on $D_{p'}$), 
$p'=a^2,q'=b^2$. 

\subhead 11.5\endsubhead
$G$ is of type $D_n$, $d=2$.

$a\ge 0$, $b\ge 0$, $a=n+1\mod 2$, $(b^2+b)/2=n\mod 2$;

$n+1=2a^2+(b^2+b)/2-1+2s$, $s\ge 1$.

$(\g_i)$-graph: 
$$\s\Ra\s\h\s\h\s\h\s\h\s\h\do\h\s\h\s\h\s\Lar\s$$
$(\b_i)$-graph: 
$$\us{B_p}\to{\x{\sm\s&\Lar&\s&\s&\do&\s\esm}}\h\s_1\h\x\s\h\s_2\h\x\s\h\do\h
\x\s\h\s_s\h\us{B_q}\to{\x{\sm\s&\do&\s&\s&\Ra&\s\esm}}$$
if $2a+b\ne 1,|4a-2b-1|\ne 3$; here 

$2p+1=(2a+b)(2a+b+1)/2,2q+1=(2a-b)(2a-b-1)/2$; 
$$\us{B_p}\to{\x{\sm\s&\Lar&\s&\s&\do&\s\esm}}\h\s_1\h\x\s\h\s_2\h\x\s\h\do\h
\x\s\Ra\s_s$$
if $2a+b\ne 1,|4a-2b-1|=3$; here $2p+1=(2a+b)(2a+b+1)/2$;
$$\s_1\Lar\x\s\h\s_2\h\x\s\h\do\h\x\s\Ra\s_s$$
if $2a+b=1,2a-b=-1$.

For $s\ge 2$, the $\fla-\sh$-diagram is
$$\fla_1^{(4a+2b+1)\T 1\T 1}\Lar\sh_2^{4\T 1\T 1}\h\sh_3^{4\T 1\T 1}\h\do\h
\sh_{s-1}^{4\T 1\T 1}\Ra\fla_s^{|4a-2b-1|\T\fra{1}{2}\T 2}$$
if $2a+b\ne 1,|4a-2b-1|\ne 3$; 
$$\fla_1^{(4a+2b+1)\T 1\T 1}\Lar\sh_2^{4\T 1\T 1}\h\sh_3^{4\T 1\T 1}\h\do\h
\sh_{s-1}^{4\T 1\T 1}\Ra\fla_s^{3\T\fra{1}{2}\T 2}$$
if $2a+b\ne 1,|4a-2b-1|=3$; 
$$\fla_1^{3\T\fra{1}{2}\T 2}\Lar\sh_2^{4\T 1\T 1}\h\sh_3^{4\T 1\T 1}\h\do\h
\sh_{s-1}^{4\T 1\T 1}\Ra\fla_s^{3\T\fra{1}{2}\T 2}$$
if $2a+b=1,2a-b=-1$.

For $s=1$, the $\fla-\sh$-diagram is $\e$.

H.A.: $\ti C_{s-1}^{sc}[{}_{4a}2_{2b+1}]$ if $s\ge 2$ and $\e$ if $s=1$.

Arithmetic diagram: $\ti D_n$, 
$u:\bi@>>>\bi$ has $<n-1$ fixed points, $\bj$ of type 
$D_{p'}\T D_{p'}\T A_{r-1}$ where $p'=a^2,r=(b^2+b)/2$.

\subhead 11.6\endsubhead
$G$ is of type $E_6$, $d=2$. ($2$ possible cuspidal local systems.)

$(\g_i)$-graph: 
$$\s\h\s\h\s\Ra\s\h\s$$
$(\b_i)$-graph: 
$$\x{\s\h\s}\h\s_1\Lar\x{\s\h\s}$$
H.A.: $\e$.

Arithmetic diagram: $\ti E_6$, $u^2=1,u\ne 1$, $\bj$ of type $E_6$ (with two
possible unipotent cuspidal representations, one the dual of the other).

\subhead 11.7\endsubhead
$G$ is of type $E_6$, $d=2$.

$(\g_i)$-graph: 
$$\s\h\s\h\s\Ra\s\h\s$$
$(\b_i)$-graph: 
$$\s_1\h\x{\s\h\s\Lar\s\h\s}$$
H.A.: $\e$.

Arithmetic diagram: $\ti E_6$, $u^2=1,u\ne 1$, $\bj$ of type $E_6$ (with a 
self-dual unipotent cuspidal representation).

\subhead 11.8\endsubhead
$G$ is of type $E_6$, $d=2$.

$(\g_i)$-graph: 
$$\s\h\s\h\s\Ra\s\h\s$$
$(\b_i)$-graph: 
$$\x\s\h\s_1\h\x\s\Lar\s_2\h\x\s$$
$\fla-\sh$-diagram:   
$$\fla_1^{4\T 1\T 2}\overset\infty\to\h\fla_2^{5\T 2\T 2}$$
H.A.:  $\bold 1\overset\infty\to\h\bold 9$.

Arithmetic diagram: $\ti E_6$, $u^2=1,u\ne 1$, $\bj$ of type $A_5$.

\subhead 11.9\endsubhead
$G$ is of type $E_6$, $d=2$.

$(\g_i)$-graph: 
$$\s\h\s\h\s\Ra\s\h\s$$
$(\b_i)$-graph: 
$$\s_1\h\s_2\h\s_3\Lar\s_4\h\s_5$$
$\fla-\sh$-diagram:    
$$\sh_1^{2\T 1\T 2}\h\sh_2^{2\T 1\T 2}\h\sh_3^{2\T 1\T 2}\Ra\sh_4^{2\T 1\T 1}
\h\sh_5^{2\T 1\T 1}$$
H.A.: $$\bold 2\h\bold 2\Lar\bold 1\h\bold 1\h\bold 1$$

Arithmetic diagram: $\ti E_6$, $u^2=1,u\ne 1$, $\bj=\e$.

\subhead 11.10\endsubhead
$G$ is of type $D_4$, $d=3$.

$(\g_i)$-graph: 
$$\s\h\s\equiv>\s$$
$(\b_i)$-graph: 
$$\s_1\h\x{\s<\equiv\s}$$
H.A.: $\e$.

Arithmetic diagram: $\ti D_4$, $u^3=1,u\ne 1$, $\bj$ of type $D_4$.

\subhead 11.11\endsubhead
$G$ is of type $D_4$, $d=3$.

$(\g_i)$-graph: 
$$\s\h\s\equiv>\s$$
$(\b_i)$-graph: 
$$\x\s\h\s_1<\equiv\x\s$$
H.A.: $\e$.

Arithmetic diagram: $\ti D_4$, $u^3=1,u\ne 1$, $\bj$ of type $D_4$ with a 
unipotent cuspidal other than that in 11.10.

\subhead 11.12\endsubhead
$G$ is of type $D_4$, $d=3$.

$(\g_i)$-graph: 
$$\s\h\s\equiv>\s$$
$(\b_i)$-graph: 
$$\s_1\h\s_2<\equiv\s_3$$
$\fla-\sh$-diagram:    
$$\sh_1^{2\T 1\T 3}\h\sh_2^{2\T 1\T 3}\equiv>\sh^{2\T 1\T 1}$$
H.A.:  $\bold 3<\equiv\bold 1\h\bold 1$.

Arithmetic diagram: $\ti D_4$, $u^3=1,u\ne 1$, $\bj=\e$.

\head Appendix. Proof of Lemma 5.5\endhead
\subhead A.1\endsubhead
We may assume that $\fg$ is simple. If $\ul=\fg$ there is nothing to prove. If
$\ul$ is a Cartan subalgebra, then $h^0=0$ and there is nothing to prove. In 
the rest of the proof we assume that $\ul\ne\fg$ and $\ul$ is not a Cartan 
subalgebra. 

Since $\ft\opl\bc h^0$ is a Cartan subalgebra of $\un{\tZ}$, it is enough to
prove the following statement:

(a) {\it Let $x,x'\in\ft,z\in\bc$ be such that $x+zh^0,x'+zh^0$ are 
$G$-conjugate in $\fg$. Then $x,x'$ are in the same $W$-orbit.}

If (a) holds for $z=1$ then it also holds for any $z\ne 0$. (We replace 
$x,x',z$ by $z\i x,z\i x',1$.) Thus it is enough to prove (a) for
$z\in\{0,1\}$.

Let $\fS=\{x+ah^0\in\ft+\bc h^0; \a(x)\ge 0\quad\f\a\in\Pi\}$. As in 6.6 we see
that $\fS$ is a fundamental domain for the action of $W$ on $\ft\opl\bc h^0$.
Hence it is enough to prove the following statement: 

{\it Let $x,x'\in\ft,z\in\{0,1\}$ be such that $x+zh^0,x'+zh^0$ belong to $\fS$
and are $G$-conjugate in $\fg$. Then $x=x'$.}
\nl
We consider the various cases separately.

For a multiset $X$ consisting of finitely many numbers in $\bc$ we denoted by
$\max X$ the complex number $x\in X$ such that $x-x'\ge 0$ for any $x'\in X$.

\subhead A.2\endsubhead
Assume that $\fg=\fs\fl_{ab}(\bc),\ul=\fs\fl_a(\bc)^b\opl\bc^{b-1}$. Here 
$a>1,b>1$.
Let $V$ be a $\bc$-vector space with basis $e_1,e_2,\do,e_{ab}$. We may assume
that $\fg=\fs\fl(V)$. We may assume that 

$x(e_{ai+l})=x_ie_{ai+l},x'(e_{ai+l})=x'_ie_{ai+l}$ for 
$i\in[0,b-1],l\in[0,a-1]$
\nl
where $x_i,x'_i\in\bc$ satisfy $\sum_ix_i=\sum_ix'_i=0$ and 
$x_i-x_{i+1}\ge 0, x'_i-x'_{i+1}\ge 0$ for $i\in[0,b-2]$ and that 

$h^0(e_{ai+l})=(a-1-2l)e_{ai+l}$ for $i\in[0,b-1],l\in[0,a-1]$. 
\nl
Since $x+zh^0,x'+xh^0$ are conjugate under $SL(V)$, they must have the same 
eigenvalues in $V$. Thus, the multisets

$X=\{x_i+z(a-1-2l)\}_{i\in[0,b-1],l\in[0,a-1]},
X'=\{x'_i+z(a-1-2l)\}_{i\in[0,b-1],l\in[0,a-1]}$
\nl
coincide. Clearly, $\max X=x_0+z(a-1)$ and $\max X'=x'_0+z(a-1)$. Since $X=X'$
we have $x_0+z(a-1)=x'_0+z(a-1)$. Hence $x_0=x'_0$. Removing 

$x_0+z(a-1),x_0+z(a-3),\do,x_0+z(-a+1)$ 

(resp. $x'_0+z(a-1),x'_0+z(a-3),\do,x'_0+z(-a+1)$) 
\nl
from $X$ (resp. $X'$) we obtain a multiset $X_1$ (resp. $X'_1$). We have 
$X_1=X'_1$. Clearly, $\max X_1=x_1+z(a-1)$ and $\max X'_1=x'_1+z(a-1)$. Since 
$X_1=X'_1$ we have $x_1+z(a-1)=x'_1+z(a-1)$. Hence $x_1=x'_1$. Continuing in 
this way we find $x_i=x'_i$ for $i\in[0,b-1]$. Hence $x=x'$.

\subhead A.3\endsubhead
Assume that $\fg=\fs\fp_{2n+2p}(\bc),\ul=\fs\fp_{2n}(\bc)\opl\bc^p$. Here 
$n>1,p>1$ and $n=(m^2+m)/2$.

Let $V$ be a $\bc$-vector space with basis 
$e_1,e_2,\do,e_{n+p},e'_{n+p},\do,e'_2,e'_1$ and with a symplectic form 
$(,):V\T V@>>>\bc$ such that $(e_i,e'_j)=\de_{ij},(e_i,e_j)=(e'_i,e'_j)=0$ for
$i,j\in[1,n+p]$. We may assume that $\fg=\fs\fp(V)$ and that

$x(e_i)=x_ie_i,x(e'_i)=-x_ie'_i,x'(e_i)=x'_ie_i,x'(e'_i)=-x'_ie'_i$ for 
$i\in[1,p]$, 

$x(e_i)=0,x(e'_i)=0,x'(e_i)=0,x'(e'_i)=0$ for $i\in[p+1,p+n]$,
\nl
where $x_i,x'_i\in\bc$ satisfy $x_i-x_{i+1}\ge 0, x'_i-x'_{i+1}\ge 0$ for 
$i\in[1,p-1]$, $x_p\ge 0,x'_p\ge 0$. We may also assume that

$h^0(e_i)=0,h^0(e'_i)=0$ for $i\in[1,p]$, 

$h^0(e_i)=c_ie_i,h^0(e'_i)=-c_ie_i$ for $i\in[p+1,p+n]$ 
\nl 
where $c_i\in\bz$. Since $x+zh^0,x'+xh^0$ are conjugate under $Sp(V)$, they 
must have the same eigenvalues in $V$. Thus, the multisets

$Y=\{x_i,-x_i (i\in[1,p]),c_i,-c_i (i\in[p+1,p+n])\}$,

$Y'=\{x'_i,-x'_i (i\in[1,p]),c_i,-c_i (i\in[p+1,p+n])\}$
\nl
coincide. Removing $\{c_i,-c_i (i\in[p+1,p+n])\}$ from $Y$ (resp. $Y'$) we 
obtain a multiset $X$ (resp. $X'$). We have $X=X'$. Clearly, $\max X=x_1$ and 
$\max X'=x'_1$. Since $X=X'$ we have $x_1=x'_1$. Removing $x_1$ (resp. $x'_1$)
from $X$ (resp. $X'$) we obtain a multiset $X_1$ (resp. $X'_1$). We have
$X_1=X'_1$. Clearly, $\max X_1=x_2$ and $\max X'_1=x'_2$. Since $X_1=X'_1$ we
have $x_2=x'_2$. Continuing in this way we find $x_i=x'_i$ for $i\in[1,p]$.
Hence $x=x'$.

\subhead A.4\endsubhead
Assume that $\fg=\fs\fo_{n+2p}(\bc),\ul=\fs\fo_n(\bc)\opl\bc^p$. Here $n>2,p>1$
and $n=m^2$. This case is completely similar to that in A.3.

\subhead A.5\endsubhead
Assume that 
$\fg=\fs\fo_{2n+4p}(\bc),\ul=\fs\fo_{2n}(\bc)\opl\fs\fl_2^p\opl\bc^p$. Here 
$n>0,p>0$ and $2n=(m^2+m)/2$.

Let $V$ be a $\bc$-vector space with basis 
$e_1,e_2,\do,e_{n+2p},e'_{n+2p},\do,e'_2,e'_1$ and with a symmetric bilinear
form $(,):V\T V@>>>\bc$ such that $(e_i,e'_j)=\de_{ij},(e_i,e_j)=(e'_i,e'_j)=0$
for $i,j\in[1,n+2p]$. We may assume that $\fg=\fs\fo(V)$ and that

$x(e_{2i-1})=x_ie_{2i-1},x(e_{2i})=x_ie_{2i-1},x(e'_{2i-1})=-x_ie_{2i-1},
x(e'_{2i})=-x_ie_{2i-1},$

$x'(e_{2i-1})=x'_ie_{2i-1},x'(e_{2i})=x'_ie'_{2i-1},
x'(e'_{2i-1})=-x'_ie_{2i-1},x'(e'_{2i})=-x'_ie_{2i-1}$
\nl
for $i\in[1,p]$,

$x(e_i)=0,x(e'_i)=0,x'(e_i)=0,x'(e'_i)=0$ for $i\in[2p+1,2p+n]$,
\nl
where $x_i,x'_i\in\bc$ satisfy 

$x_i-x_{i+1}\ge 0, x'_i-x'_{i+1}\ge 0$ for $i\in[1,p-1]$, $x_p\ge 0,x'_p\ge 0$.
\nl
We may also assume that

$h^0(e_{2i-1})=e_{2i-1},h^0(e_{2i})=-e_{2i},h^0(e'_{2i-1})=-e_{2i-1},
h^0(e'_{2i})=e_{2i}$
\nl
for $i\in[1,p]$ and 

$h^0(e_i)=c_ie_i,h^0(e'_i)=-c_ie_i$ for $i\in[2p+1,2p+n]$ 
\nl
where $c_i\in\bz$. Since $x+zh^0,x'+xh^0$ are conjugate under $SO(V)$, they 
must have the same eigenvalues in $V$. Thus, the multisets

$Y=\{x_i+z,x_i-z,-x_i+z,-x_i-z (i\in[1,p]), c_i,-c_i (i\in[2p+1,2p+n])\}$,

$Y'=\{x'_i+z,x'_i-z,-x'_i+z,-x'_i-z (i\in[1,p]), c_i,-c_i (i\in[2p+1,2p+n])\}$
\nl
coincide. Removing $\{c_i,-c_i (i\in[2p+1,2p+n])\}$ from $Y$ (resp. $Y'$) we 
obtain a multiset $X$ (resp. $X'$). We have $X=X'$.

Clearly, $\max X=x_1+z$ and $\max X'=x'_1+z$. Since $X=X'$ we have 
$x_1+z=x'_1+z$ hence $x_1=x'_1$. Removing $x_1+z,x_1-z,-x_1+z,-x_1-z$ (resp. 
$x'_1+z,x'_1-z,-x'_1+z,-x'_1-z$) from $X$ (resp. $X'$) we obtain a multiset 
$X_1$ (resp. $X'_1$). We have $X_1=X'_1$. Clearly, $\max X_1=x_2+z$ and 
$\max X'_1=x'_2+z$. Since $X_1=X'_1$ we have $x_2+z=x'_2+z$ hence $x_2=x'_2$. 
Continuing in this way we find $x_i=x'_i$ for $i\in[1,p]$. Hence $x=x'$.

\subhead A.6\endsubhead
Assume that 
$\fg=\fs\fo_{2n+1+4p}(\bc),\ul=\fs\fo_{2n+1}(\bc)\opl\fs\fl_2^p\opl\bc^p$.
Here $n\ge 0,p>0$, $2n+1=(m^2+m)/2$. This case is completely similar to that in
A.5.

\subhead A.7\endsubhead
Assume that $\fg$ is of type $E_6$ and $\ul\cong\fs\fl_3(\bc)^3\opl\bc^2$. We
number the vertices of the Coxeter diagram by $1,2,\do,6$ where the edges are
$1-2-3-4-5$ and $3-6$. For $i\in[1,6]$ define $u_i\in\ft$ by 
$\a_j(u_i)=\de_{ij}$ for all $j\in[1,6]$. Then $\ft$ is spanned by $u_3,u_6$. 
We may assume that $h^0=2\cha_1+2\cha_2+2\cha_4+2\cha_5$.

We have $x=au_3+bu_6$ $x'=a'u_3+b'u_6$ where $a,b,a',b'\in\bc$ are $\ge 0$. Let
$Y$ be the multiset consisting of the numbers

$2a+b+2z,2a+b,2a+b-2z,a+b+2z,a+b,a+b-2z,a+2z,a,a-2z,4z,2z,2z$,
\nl
their negatives, and of $0,0,0$. Let $Y'$ be the multiset consisting of the 
numbers

$2a'+b'+2z,2a'+b',2a'+b'-2z,a'+b'+2z,a'+b',a'+b'-2z$,

$a'+2z,a',a'-2z,4z,2z,2z$,
\nl
their negatives, and of $0,0,0$. The eigenvalues of $x+zh^0$ (resp. $x'+zh^0$)
on a minuscule $\fg$-module $V$ are the $27$ numbers in the multiset $Y$ (resp.
$Y'$). Since $x+zh^0,x'+zh^0$ are in the same $G$-orbit, they have the same 
eigenvalues on $V$. Thus, $Y=Y'$. 

Removing $4z,2z,2z$ from $Y$ (resp. $Y'$) we obtain a multiset $X$ (resp. $X'$)
with $24$ elements. We have $X=X'$. Clearly, $\max X=2a+b+2z$ and 
$\max X'=2a'+b'+2z$. It follows that

(a) $2a+b+2z=2a'+b'+2z$.
\nl
Removing from  $X$ (resp. $X'$) the  numbers $2a+b+2z,2a+b,2a+b-2z$ (resp. 
$2a'+b'+2z,2a'+b',2a'+b'-2z$) we obtain a multiset $X_1$ (resp. $X'_1$). By 
(a), we have $X_1=X'_1$. Clearly, $\max X_1=a+b+2z$ and $\max X'_1=a'+b'+2z$.
It follows that 

(b) $a+b+2z=a'+b'+2z$.
\nl
From (a),(b) we deduce that $a=a',b=b'$. Thus $x=x'$ as required.

\subhead A.8\endsubhead
Assume that $\fg$ is of type $E_7$ and $\ul\cong\fs\fl_2(\bc)^3\opl\bc^4$. We
number the vertices of the Coxeter diagram by $1,2,\do,7$ where the edges are
$1-2-3-4-5-6$ and $3-7$. For $i\in[1,7]$ define $u_i\in\ft$ by 
$\a_j(u_i)=\de_{ij}$ for all $j\in[1,7]$. Then $\ft$ is spanned by 
$u_1,u_2,u_3,u_5$. We may assume that $h^0=\cha_4+cha_6+\cha_7$.

We have $x=au_1+bu_2+cu_3+du_5$ $x'=a'u_1+b'u_2+c'u_3+d'u_5$ where
$a,b,c,d,a',b',c',d'\in\bc$ are $\ge 0$. Let $Y$ be the multiset consisting of
$$\align&a+2b+3c+2d+z,a+2b+3c+2d-z,a+2b+3c+d+z,a+2b+3c+d-z,\\&
a+2b+2c+d+z, a+2b+2c+d-z,a+b+2c+d+z,a+b+2c+d-z,\\&
a+b+c+d+z,a+b+c+d-z,b+2c+d+z,b+2c+d-z,\\&
b+c+d+z,b+c+d-z,a+b+c+z,a+b+c-z,\\&
c+d+z,c+d-z,b+c+z,b+c-z,d+z,d-z,c+z,c-z,3z,z,z,z\endalign$$
and their negatives. Let $Y'$ be the multiset obtained from $Y$ by replacing 
$a,b,c,d,z$ by $a',b',c',d',z$.

The eigenvalues of $x+zh^0$ (resp. $x'+zh^0$) on the minuscule $\fg$-module $V$
are the $56$ numbers in the multiset $Y$ (resp. $Y'$). Since $x+zh^0,x'+zh^0$ 
are in the same $G$-orbit, they have the same eigenvalues on $V$. Thus, $Y=Y'$.
Removing $3z,z,z,z$ from $Y$ (resp. $Y'$) we obtain a multiset $X$ (resp. $X'$)
with $52$ elements. We have $X=X'$. Clearly, $\max X=a+2b+3c+2d+z$ and 
$\max X'=a'+2b'+3c'+2d'+z$. It follows that

(a) $a+2b+3c+2d+z=a'+2b'+3c'+2d'+z$.
\nl
Removing from  $X$ (resp. $X'$) the  numbers $a+2b+3c+2d+z,a+2b+3c+2d-z$ (resp.
$a'+2b'+3c'+2d'+z,a'+2b'+3c'+2d'-z$) we obtain a multiset $X_1$ (resp. $X'_1$).
By (a), we have $X_1=X'_1$. Clearly, $\max X_1=a+2b+3c+d+z$ and
$\max X'_1=a'+2b'+3c'+d'+z$. It follows that 

(b) $a+2b+3c+d+z=a'+2b'+3c'+d'+z$.
\nl
Removing from  $X_1$ (resp. $X'_1$) the  numbers $a+2b+3c+d+z,a+2b+3c+d-z$ 
(resp. $a'+2b'+3c'+d'+z,a'+2b'+3c'+d'-z$) we obtain a multiset $X_2$ (resp. 
$X'_2$). By (b), we have $X_2=X'_2$. Clearly, $\max X_2=a+2b+2c+d+z$ and
$\max X'_2=a'+2b'+2c'+d'+z$. It follows that 

(c) $a+2b+2c+d+z=a'+2b'+2c'+d'+z$.
\nl
Removing from  $X_2$ (resp. $X'_2$) the  numbers $a+2b+2c+d+z,a+2b+2c+d-z$ 
(resp. $a'+2b'+2c'+d'+z,a'+2b'+2c'+d'-z$) we obtain a multiset $X_3$ (resp. 
$X'_3$). By (c), we have $X_3=X'_3$. Clearly, $\max X_3=a+b+2c+d+z$ and
$\max X'_3=a'+b'+2c'+d'+z$. It follows that 

(d) $a+b+2c+d+z=a'+b'+2c'+d'+z$.
\nl
From (a),(b),(c),(d) we deduce that $a=a',b=b',c=c',d=d'$. Thus $x=x'$ as 
required.

Lemma 5.5 is proved.

\enddocument